\magnification=1200

\def\proof{\medskip\noindent{\bf Proof. }}
\def\remark{\medskip\noindent{\bf Remark }}
\def\example{\noindent{\bf Example }}
\def\claim{\noindent{\bf Claim }}

\def\definition {\noindent{\bf Definition}}

\def\today{\ifcase\month\or January \or February\or
March\or April\or May\or June\or July\or August\or
September\or October\or November\or December\fi
\space\number\day, \number\year}

\def\O{{\cal O}}

\def\E{{\cal E}}

\def\I{{\cal I}}
\def\J{{\cal J}}
\def\L{{\cal L}}

\def\P{{\bf P}}
\def\Q{{\bf Q}}
\def\P{{\bf P}}
\def\Z{{\bf Z}}
\def\C{{\bf C}}
\def\F{{\bf F}}
\def\S{{\bf S}}
\def\R{{\bf R}}

\def\f{\varphi}
\def\ra{\rightarrow}
\def\raa{\longrightarrow}
\def\iso{\simeq}

\def\har#1{\smash{\mathop{\hbox to .8 cm{\rightarrowfill}}
\limits^{\scriptstyle#1}_{}}}

\font\ninerm=cmr10 at 10truept
\font\ninebf=cmbx10 at 10truept
\font\nineit=cmti10 at 10truept
\font\bigrm=cmr12 at 12 pt
\font\bigbf=cmbx12 at 12 pt
\font\bigit=cmti12 at 12 pt

\def\smalltype{\let\rm=\ninerm \let\bf=\ninebf
\let\it=\nineit  \baselineskip=9.5pt minus .75pt
\rm}
\def\bigtype{\let\rm=\bigrm \let\bf=\bigbf\let\it=\bigit
\baselineskip=12 pt minus 1 pt \rm}

\phantom{.}

\vskip 1 cm

{\bigtype{
\centerline{{\bf On contractions of smooth varieties.}$^*$}
\bigskip
\centerline{{\sl Marco Andreatta}$^1$ and {\sl Jaros{\l}aw A. Wi\'sniewski}$^2$}
\medskip
}}

\medskip
\item{$^1$}      Dipartimento di Matematica, Universit\'a di Trento,
38050 Povo (TN), Italy
\par
e-mail :  {\tt andreatt@science.unitn.it}
\par
\item{$^2$}     Instytut Matematyki UW, Banacha 2, 02-097 Warszawa, Poland
\par
e-mail:   {\tt jarekw@mimuw.edu.pl}
\smallskip
\item{$^*$} Dedicated to Silvia  and Ma{\l}gosia.
\medskip
\noindent{\bf Abstract.}
Let $\f: X \ra Z$ be a proper surjective map from a smooth complex manifold
$X$ onto a normal variety $Z$. If $\f$ has connected fibers and $-K_X$
is $\f$-ample then $\f$ is called a good contraction.
In the present paper we study good contractions, fibers of which have 
dimension less or equal than two:
after describing possible two dimensional isolated fibers
we discuss their scheme theoretic structure and the geometry of 
$\f:X\ra Z$ nearby such a fiber. 
If $dimX=4$ and $\f$ is birational with an isolated
2 dimensional fiber then we obtain a complete description of $\f$.
We provide also a description of a 4 dimensional 
conic fibration with an isolated
fiber which is either a plane or a quadric.
We construct pertinent examples.


\smallskip

\beginsection Introduction.

The present paper is about maps of complex algebraic varieties.  A
{\sl contraction} $\f:X\ra Z$ is a proper surjective map of normal
varieties with connected fibers. The contraction $\f$ is called {\sl
good} if the anticanonical divisor $-K_X$ of $X$ is $\f$-ample.  Good
contractions occur naturally in classification theories of algebraic
varieties. This is most apparent in the Minimal Model Program where
the point is to use good contractions and birational transformation
related to them to produce a minimal model of a given variety.  In
what follows we study the local structure of a good contraction of a
smooth variety.

\par
As they appear naturally in classification of algebraic varieties,
good contractions have been studied on many occasions. First, let us
note that a Fano manifolds $X$ admits a good contraction to a point
and from this position the theory of good contractions is just a very
broad generalization of the theory of Fano manifolds. Next let us
notice that the classical and the contemporary adjunction theory is
about the adjunction morphism which is nothing else but a special
good contraction.  In particular, the classical works of Enriques and
Castelnuovo provide a complete description of good contractions of
smooth surfaces.  However, it had to come to the case of threefolds
to establish the fundamental role of good contractions in
classification of algebraic varieties. A complete description of
(elementary) good contractions of smooth 3-folds was given by S.~Mori
[Mo1] as the base of the Minimal Model Program in dimension 3. Also,
the objective of the Program in dimension 3, that is the existence of
minimal models, was achieved by studying good contractions of
varieties admitting good (terminal) singularities, see [Mo2].
Subsequently, the theory of good contractions of terminal threefolds
was extended in [Ko-Mo].
\par

In this paper we look at good contractions of smooth varieties of
dimension $n \geq 3$ with low dimensional fibers (i.e. of dimension
$\leq 2$). In particular we provide a complete list of isolated two
dimensional fibers. Then, after reviewing the 3 dimensional case we
focus on the case $n= 4$ and $\f$ birational. In this range some
results were already proved by T. Ando [An], M. Beltrametti [Be]
and Y.Kawamata [Ka1].
\par

In order to understand the local structure of a good contraction
$\f:X\ra Z$ we will assume that the target $Z$ is affine and $z\in Z$
is a fixed (geometric) point. First we will describe the geometric
structure of the fiber $F=\f^{-1}(z)$. Subsequently we will discuss
its normal sheaf of $F$ and the fiber scheme structure on $F$.
Finally we will provide a description of the singularity of $Z$ at
$z$ and the description of $\f$ around $F$.
\par

The paper is divided into six sections. The first two sections are
preparatory.  In Section 1 we recall first the definition of good and
crepant contractions and subsequently we present our main technical
tools. This includes theorems about vanishing, relative base point
freeness and deformation of rational curves.  In Section 2 we collect
pertinent results about blow-ups and about rank two vector bundles on
Fano varieties.  In particular we discuss in detail the blow-up of a
smooth variety along a codimension two locally complete intersection
subvariety and we prove a generalized Castelnuovo's contraction
theorem (2.4).

In Section 3 we present several examples of good contractions.
They are constructed with various algebraic geometry techniques:
as complete intersections in projective bundles, or as blow-ups and
blow-downs of special varieties, or as double coverings of special
varieties, or as toric varieties. Although the constructions are not
particularly difficult we found the results of some of them rather
surprising.  Actually, the list of examples grew together with our
understanding of the classification of good contractions. In the
subsequent section we show that the list of examples covers (almost)
all possible two dimensional isolated fiber of a good contraction.
In Section 4 we recall also the case of one dimensional fiber and at
the end of the section we reprove the theorems for good contractions
of a $3$-fold (see (4.13)).  The classification of 2 dimensional 
fibers is achieved in a two step argument.  First we prove that a
two dimensional fiber is normal and it has Fujita's $\Delta$-genus
equal to zero (see (4.2.1))  Then, to get
the final list of possible fibers, we compare deformations of
rational curves inside the fiber with the deformations of such curves
inside the ambient variety. The main results are summarized in Table II
and (4.5), (4.7) and (4.11).

In Section 5 we discuss the scheme theoretic structure of a fiber of
a good contractions. We prove that if the fiber $F$ is a locally
complete intersection and the blow-up along this fiber has good
singularities then the conormal bundle of the fiber provides a lot of
information about the contraction $\f$ around $F$. In particular, if
the conormal bundle is spanned by global sections then the fiber
structure on $F$ is trivial and the contraction can be factored
through the blow-up of $Z$ along the maximal ideal of $z$ (see
(5.5)).  Subsequently, we focus on the case $n=4$ and we prove that
the good situation mentioned above occurs for any two dimensional
isolated fiber of a birational good contraction.  The strategy in
this section is as follows: we take a general smooth divisor
$X'\in|-K_X|$ and then we consider the map $\f':=\f_{|X'}: X' \ra Z'$
which is a {\sl crepant contraction} of a smooth 3-fold.  The
structure of such a map is rather well understood by results of the
Minimal Model Program in dimension 3 (we discuss it in detail in
(5.6)).  Then we apply an {\sl ascending lemma} (5.7.2)
which gives the spannedness of the
conormal of the fiber.

Section 6 concludes the paper with a geometric description of a good
birational contraction $\f$ of a smooth 4 fold.  In particular, we
show that $\f$ can be {\sl resolved} in terms of some special
blow-ups and blow-downs.  Next we describe the singularities of the
target $Z$ and of the image of the exceptional locus $\f(E) = S$.
Here we apply the classification of spanned rank two vector bundles
on Fano manifolds presented in Section 2 as well as arguments
involving Hilbert scheme.

The following is a summary of the 4 dimensional birational result:

\bigskip
{\bf Theorem.} {\sl Let $\f:X\ra Z$ be a birational good contraction
from a smooth variety $X$ of dimension $4$ onto a normal variety $Z$
(possibly affine).  Let $F=\f^{-1}(z)$ be a (geometric)
fiber of $\f$ such that $\hbox{dim}F=2$. Assume that all other
fibers of $\f$ have dimension $<2$ and all components of the
exceptional locus $E$ of $\f$ meet $F$ (this may be achieved by
shrinking $Z$ to an affine neighbourhood of $z$ and restricting $\f$
to its inverse image, if necessary).

If $\f$ is not divisorial then $E=F\iso\P^{2}$
and its normal is $N_{F/X} = \O(-1) \oplus \O(-1)$.
In this situation the {\it flip} of $\f$
exists. (This was the situation studied in [Ka1]).

If $E$ is a divisor then $Z$ as well as $S:=\f(E)$ are
smooth outside of $z$. Moreover, outside of $F$ the map $\f$ is a
simple blow-down of the divisor $E$ to the surface $S\subset Z$.  The
scheme theoretic fiber structure over $F$ is trivial, that is the
ideal $\I_F$ of $F$ is equal to the inverse image of the maximal
ideal of $z$, that is $\I_F=\f^{-1}(m_z)\cdot\O_X$.

The fiber $F$ and its conormal $\I_F/\I_F^2$ as well as
the singularity of $Z$ and $S$ at $z$ can be described as follows

\medskip
\settabs\+$\P^1 \times \P^1$\ \ \ \ \
&$T_{\P^2} (-1) \cup \O \oplus \O(1)$  \ \ \ \ \ \
&quadratic sing. \ \ \ &cubic singularity \ \cr
\+ $F$  &$N^*_{F/X}$ &$Sing Z$ &$Sing S$\cr
\hrule
\smallskip
\+$\P^2$ &$T(-1)\oplus\O(1)/\O $
&cone over $\Q^3$&smooth \cr
\+$\P^2$ &$\O^{\oplus 4}/\O(-1)^{\oplus 2}$
&smooth &cone over rational twisted cubic\cr
\+ quadric &spinor bundle from $\Q^4$ &smooth&non-normal\cr
\smallskip
\hrule
\bigskip
The quadric fiber can be singular, even reducible, and in the subsequent
table we present a refined description of its conormal bundle.
The last entry in the table provides information
about the ideal of a suitable surface $S$
which is computed via degeneracy locus technique in Example 3.2.
\medskip
\settabs\+ quadric cone \ \ \ \ \ & the normal bundle\ \ \ \ \ \
 \ \ \ \ \ \ \ \ \  &
the ideal of $S$ in $\C[[x,y,z,t]]$\cr
\+ quadric & conormal bundle & the ideal of $S$ in $\C[[x,y,z,t]]$\cr
\smallskip
\hrule
\smallskip
\+$\P^1 \times \P^1$ &$\O(1,0)\oplus\O(0,1)$
&$(xz,xt,yz,yt)$\cr
\+ quadric cone & $0\ra\O\ra N^*\ra \J_{line}\ra 0$
 & generated by 5 cubics \cr
\+$\P^2 \cup \P^2$ &$T_{\P^2} (-1) \cup (\O \oplus \O(1))$
& generated by 6 quartics\cr
\smallskip
\hrule
}

\bigskip
A similar classification is expected for contractions of fiber type
from a smooth 4 dimensional projective variety with an isolated two
dimensional fiber.  The list
of possible fibers is given in Section 4 (see (4.11)) and many
examples are discussed in Section 3.  In this case similar results
were proved also by Kachi [Kac]. The understanding of the normal of
such a fiber requires a new approach, i.e. the {\sl ascending
lemma} (5.7.2) has to be replaced by a
{\sl trace argument} (see 5.9). As the result we obtain a structure theorem
for 4 dimensional conic fibrations with an isolated 2 dimensional 
fiber which is either $\P^2$ or a quadric.

Due to the relative base point freeness proved in [A-W], the
4-dimensional results can be extended for adjoint contractions of
varieties of higher dimension (see the corollary (5.8.1))

\medskip

In the course of the work on the present paper we enjoyed hospitality
of Universit\'a degli Studi di Trento, Uniwersytet Warszawski and
Max-Planck-Institut f\"ur Mathematik in Bonn.  We would like to thank
these institutions for providing support as well as excellent working
conditions and stimulating atmosphere.  The first named author was
moreover partially supported by MURST and GNSAGA.  The second named
author was partially supported by an Italian CNR grant (GNSAGA) and
by a Polish grant KBN GR564 (2P03A01208).

\bigskip

\beginsection 0. Notation and assumptions.

\par
We work with schemes defined over complex numbers.
In particular, a variety is a puredimensional
reduced separated scheme of finite type over $\C$,
a curve (a surface) is a variety of pure dimension 1 (2, respectively)
(thus it does not have to be irreducible).

\par
On a variety $X$ by $K_X$ we will denote its canonical divisor.
If $K_X$ is Cartier the associated line bundle will be denoted by the
same name. More generally: we will confuse Cartier divisors and line
bundles whenever it makes sense. We will also identify vector bundles
and locally free sheaves. Whenever possible,
the Chern or Segre classes will be identified with integers.
If $\E$ is a vector bundle over a variety $X$ then
$$\P(\E)=Proj_X(\bigoplus_{m>0} S^m\E)$$ is a {\sl projective bundle}
with a relatively ample line bundle $\O_{\P(\E)}(1)$.
If $\I$ is a coherent sheaf of ideals on $X$ then
$$\hat X = Proj_X(\bigoplus_{m\geq 0} \I^m)$$
is the {\sl blowing up of $X$ with respect to the coherent sheaf of ideal
$\I$}. If $F$ is the closed subscheme of $X$ corresponding to $\I$
we also call $\hat X$ the {\sl blowing up of $X$ along $F$}.

A Hirzebruch surface $\F_r$ is a $\P^1$-bundle $\P(\O(r)\oplus\O)$
over projective line $\P^1$ with a unique section $C_0\subset \F_r$
(isomorphic to $\P^1$) such that $C_0^2=-r\leq 0$.
A fiber of the projection $\F_r\ra\P^1$ will be denoted by $f$.
A (normal) cone $\S_r$ is defined by contracting $C_0\subset \F_r$ to
a normal point; in terms of projective geometry $\S_r$ is a cone over
$\P^1\hookrightarrow\P^r$ embedded via Veronese map ($r$-uple embedding).
The restriction of the hyperplane section line bundle
from $\P^{r+1}$ to $\S_r$ will be denoted by $\O_{\S_r}(1)$;
the pull-back of this bundle to $\F_r$ is $\O(C_0+rf)$.

Let $\Q_4 \iso Gr(1,3)$ be the smooth 4-dimensional quadric,
identified with the Grassmaniann of lines in $\P^3$, and let $\cal S$ be
the universal bundle which we call also the spinor bundle on
$\Q_4$. Consider a codimension 2 linear section
$i: V \hookrightarrow Gr(1,3)$. The surface
$V$ is again a quadric, either $\Q_2=\F_0=\P^1\times\P^1$, or
the quadric cone $\S_2$, or a reducible quadric, i.e. $V \iso \P^2 \cup \P^2$
where the two $\P^2$ intersect along a line $l$.
We will denote by $\O_V(1)$ (or just $\O(1)$) the restriction of
$\O(1)$ from $\Q_4$.
The {\sl spinor bundle} over $V$, which we denote again by ${\cal S}$,
is the rank two vector bundle
defined as $i^*(\cal S)$.
We note that if $V\iso \P^1 \times\P^1$
is a smooth quadric then ${\cal S}= \O(-1,0) \oplus \O(0,-1)$;
if $\pi: \F_2\ra \S_2=V$ is the resolution of the vertex of the singular
quadric, then the pull-back of ${\cal S}$ is in the non-splitting extension
$0\ra\O(-f)\ra f^*{\cal S}\ra \O(-f-C_0)\ra 0$.
If $V \iso \P^2 \cup \P^2 = P_1 \cup P_2$, then ${\cal S}$ restricted
to $P_1$ is $T\P^2(-2)$ while restricted to $P_2$ is $\O\oplus\O(-1)$.

\medskip
Let $\cal E$ be a rank $r$
vector bundle on a projective variety $X$ and
let $C \cong \P^1 \subset X$ be a rational curve in $X$.
The splitting type of $E$ on $C$ is a sequence of $r$ numbers,
$(a_1,\dots,a_{r})$ such that
$E_{|C} = \O(a_1) \oplus ...\oplus \O(a_{r})$;
we assume that $a_1 \leq a_2 ...\leq a_r$.

\medskip
For other definitions and notations connected with the theory of
minimal models that we will use in the paper we refer in general the reader
to [K-M-M]; for the main ones see also the next section.

\bigskip
\beginsection 1. Good contractions: definitions and fundamentals tools.

(1.0) A {\sl contraction} is a proper map
$\f:X\ra Z$ of normal irreducible varieties
with connected fibers. We assume that a contraction is not an isomorphism.
The map $\f$ is birational or otherwise
$dimZ< dimX$, in the latter case we say that $\f$ is of fiber type.
The exceptional locus $E(\f)$ of a birational contraction $\f$
is equal to the smallest subset of $X$ such that
$\f$ is an isomorphism on $X\setminus E(\f)$.
\par
Throughout the paper we will assume that $X$ is smooth.
In this situation the contraction $\f$ is called {\sl good} if
the anti-canonical divisor $-K_X$ is $\f$-ample.
If the map $\f$ is birational and if $K_X$ is a pull-back
of a line bundle from $Y$ then we say that $\f$ is crepant.
We say that $\f$ is {\sl elementary} if $PicX/\f^*(PicZ)\iso\Z$.
An elementary contraction is called {\sl small} if its exceptional locus
is of codimension $\geq 2$.
\par

In the present paper we are interested in the local description
of a contraction: we would like to know a structure of the target
$Z$ and of the fiber of the map $\f$.
Thus we choose a point $z\in Z$, we assume that the target $Z$ is affine,
and we consider the topological,
(or set theoretical) fiber of $\f$ over $z$, that is $\f^{-1} (z)$.
The set $\f^{-1}(z)$ may be reducible, however we will usually assume that
(unless otherwise specified) $\f^{-1} (z)$ is equidimensional
(this is because we deal with low dimensional fibers, see (4.1)).
We have two natural scheme structures on $\f^{-1} (z)$. One is
the scheme theoretic fiber structure, which we denote by $\tilde F$,
which is the closed subscheme of $X$
defined by the ideal $\I_{\tilde F}:=\f^{-1}(m_z)\cdot\O_X$.
Since $Z$ is affine and normal, and $\f$ is proper with connected fiber,
the ideal $\f^{-1}(m_z)\cdot\O_X$ is generated by
global functions on $X$ vanishing along $\f^{-1}(z)$.
The other is what we can call the geometric structure, which we denote
by $F$; this is the smallest scheme structure on $\f^{-1}(z)$. With this
structure $F$ is a variety, that is it is reduced and it has
no embedded component.

For a good contraction $\f$ we
will also consider a $\f$-ample line bundle $L$ such that
$K_X+L$ is a pullback of a line bundle from $Z$, if $Z$ is affine
then $K_X+L=\f^*(\O_Z)$.

\medskip \noindent
(1.1) Let us begin with a classical example.
The following is the list of all possible good contractions $\f:X\ra Z$ of
smooth surfaces:
\item{(a)} $Z$ is a point and $X$ is a  del Pezzo surface;
\item{(b)} $Z$ is a smooth curve and $\f: X\ra Z$ is a conic
or $\P^1$-bundle, in particular every fiber $\tilde F$ is reduced
and isomorphic to $\P^1$ or to union of two $\P^1$'s meeting transversally;
\item{(c)} $Z$ is  a smooth surface (thus $\f$ is birational) and
the exceptional locus consists of disjoint smooth rational curves
with normal bundle $\O(-1)$, thus $\f$ is a composition of blow-downs
of disjoint rational curves to smooth points on $Z$.
\item{}
\par Similarly one can describe crepant birational contractions of
smooth surfaces (see (1.5.2).
\par
The description of 2-dimensional contractions was known classically.
To understand them it is enough to apply adjunction formula, Grauert
criterion and the theory of divisors on surfaces.
\par
To understand higher
dimensional contractions one has to use some other properties
of good contractions. The fundaments of the theory were set in the 80's
by S. Mori, Y. Kawamata, J. Koll\'ar, M. Reid.
In particular, in the famous paper [Mo1], S. Mori produced
the list of all possible good (elementary) contractions for smooth three
dimensional projective varieties (which will be reproved
here in (4.1) and (4.14)).

\bigskip
The aim of this section is to recall to the reader some properties of
good contractions and to state them in the form which is convenient
for our purposes.
\par
The chief tool is the vanishing theorem due to Y.Kawamata, E.Viehweg and
J. Koll\'ar (see [K-M-M], section (1-2), or [E-W], corollary 6.11):

\proclaim Theorem (1.2). {\bf (Vanishing theorem)}
Let $\f:X\ra Z$ be a good or crepant contraction with target $Z$
being affine. Assume that $L$ is a $\f$-ample line bundle.
Then for any non-negative integer $t$ we have
$$H^{i}(X,tL)=0 \hbox{ for } i > 0.$$
If $\f$ is a good contraction and $K_X+L=\f^*(\O_Z)$ then also
$$H^{i}(X,-L)=0 \hbox{ for } i > dimX - dimZ.$$

\par
Let us note that although in the present paper we discuss only
the case of smooth $X$, the vanishing theorem and many of its consequences
remain true if we allow that $X$ has log terminal singularities.
\par
In the case of good contractions the fiber structure scheme
$\tilde F$ has nice properties.
Namely, its structural sheaf admits all the vanishings which hold for the
ambient space.
The following lemma is used very often; for the proof we
refer the reader to the following papers:
([Mo1];(3.20) and (3.25.1)),  ([Fu];(11.3)), ([An]; (2.2)),
([Y-Z], lemma 4).

\proclaim Lemma (1.2.1). Let $\f:X\ra Z$ be a good contraction
and let $L$, $z\in Z$, $F$ and $\tilde F$ be as in (1.0).
Moreover, let $\hat F$ be a scheme structure on $\f^{-1}(z)$
defined by an ideal $\f^{-1}(\I)\cdot\O_X$, where $\I$ is an ideal of
a zero dimensional subscheme of $Z$ supported at $z$ (in particular
$\hat F=\tilde F$ if $\I=m_z$).
If either $t \geq 0$ and $i\geq 1$ or $t=-1$ and $i > dimX-dimZ$ then
$$H^{i}(\hat F,tL_{|\hat F})=0.\leqno(a)$$
Let $F'$ be any subscheme of $X$ whose support is contained in $F$
so that $\f(F') = z$.
If either $t\geq 0$ and $r = dimF$ or  $t =-1$ and $r \geq
max\{dimF, dimX-dimZ +1\}$ then
$$H^{r}(F',tL_{|F'})=0.\leqno(b)$$

\proof Suppose that $\I$ is generated by functions
$\{f_1,\dots,f_r\}$. Then let us consider a sequence
$X=S_0\supset S_1\supset\dots \supset S_r=\hat F$ of
subschemes of $X$, each $S_k$ defined in $X$ by functions
$\{f_1\circ\f,\dots,f_k\circ\f\}$.
Now vanishing (a) is proved for all schemes $S_k$
by induction on $k$ if we consider
sequences
$$0\raa tL_{|S_k}\har{\cdot f_{k+1}} tL_{|S_k}\raa tL_{|S_{k+1}}\raa 0$$
and we start with the vanishing (1.2) --- see also [An] or [Fu].
To prove (b) let us note that any $F'$ supported on $\f^{-1}(z)$
is contained in a subscheme $\hat F$ which is supported on $\f^{-1}(z)$
and defined by global functions.
Since $F' \subset \hat F$ the restriction map $tL_{|\hat F} \ra tL_{|F'}$
is surjective and thus also
$H^{r}(\hat F,tL_{|\hat F}) \ra H^{r}(F',tL_{|F'})$ is
surjective for $r \geq dimF$. This proves the part (b).

\proclaim Lemma (1.2.2). Let $\f:X\ra Z$ be a good contraction
and let $L$, $\hat F$ and $F'$ be as in (1.2.1).
Let also $X'\in |L|$ be the zero locus of
a non-trivial section of $L$. Then we have
$$H^{i}(\hat F \cap X',tL_{|\hat F\cap X'})=0\leqno(a)$$
if either $t \geq 1$ and $i\geq 1$
or $t=0$ and  $i \geq max\{dimX-dimZ, 1\}$; we also have
$$H^{r}(F' \cap X',tL_{|F'\cap X'})=0 \leqno(b)$$
if either $t \geq 1$ and $r = dimF\cap X'$,
or $t=0$ and $r \geq max\{dimX-dimZ, dimF\cap X'\}$.

\proof The first part of the lemma follows from the previous one if we
consider the cohomology sequence associated to the exact sequence
$$0 \ra -L_{|\hat F} \ra \O_{\hat F} \ra \O_{\hat F\cap X'}\ra 0$$
tensored by $tL$.
The part (b) is proved similarly as in the previous lemma.

\proclaim Lemma (1.2.3).
Let $\f:X\ra Z$ be a crepant contraction
and let $\hat F$, $F$ and $F'$ be as in (1.2.1). Then
$$H^{i}(\hat F,\O_{\hat F})=0  \hbox{ \ for \ }i\geq 1 \hbox {\ and \
}\leqno(a)$$
$$H^{r}(F',\O_{F'})=0 \hbox{ \ for \ } r\geq dimF. \leqno(b)$$

\proof The same as of the above.

\medskip
We also note the following important consequence of the vanishing:

\proclaim Theorem (1.2.4). Let $\f: X\ra Z$ be a good contraction of a 
smooth variety $X$ of dimension $n$ onto a normal variety $Z$
of dimension $m$. Then $Z$ has rational singularities and
$R^{n-m}\f_*(K_X)=\omega_Z$, where $\omega_Z$ is the dualizing sheaf of $Z$.

\proof The rationality of singularities follows immediately from the
vanishing (see, [Ko1, Cor.~7.4] for $\f$ of fiber type. The descending
property of the canonical sheaf was noted first by Kempf 
[Ke, pp.~49--51]
in case $m=n$ and for $m\leq n$ it was proved by Koll\'ar [Ko1, Prop.~7.6] 
in case $Z$ is smooth. The general case is obtained by compilation of these
two result and application of Grothendieck spectral sequence.
Namely, let $\alpha:\tilde Z\ra Z$ be a desingularisation of $Z$ and 
let $\tilde X$ be a desingularisation of the fiber product
$X\times_Z\tilde Z$ with the induced morphisms
$\beta: \tilde X\ra X$ and $\tilde\f:\tilde X\ra\tilde Z$,
such that $\f\circ\beta=\alpha\circ\tilde\f$. Thus, by Grothendieck
spectral sequence, the sequences $R^i\f_*(R^j\beta_*(K_{\tilde X}))$
and $R^i\alpha_*(R^j\tilde\f_*(K_{\tilde X}))$ have the same limit.
But, because of Grauert-Riemenschneider vanishing 
$R^j\tilde\f_*(K_{\tilde X})=0$ for $j>0$ and $\tilde\f_*(K_{\tilde X})
=K_{\tilde Z}$ and by Kempf's result 
$R^{n-m}(\alpha\circ\tilde\f)_*(K_{\tilde X})=\omega_Z$.
On the other hand $R^i\f_*(R^j\beta_*(K_{\tilde X}))=0$ for $i>0$, 
[Ko1, Thm.~(3.8.i)], and $\beta_*(K_{\tilde X})=K_X$ so that the other 
sequence degenerates to $R^{i}(\f\circ\beta)_*(K_{\tilde X})
=R^i\f_*(K_X)$, and we are done.

\bigskip

Another feature of good contractions is the special
behaviour of some divisors, we will use it to apply
some inductive arguments. Namely, in the set up of (1.0)
we will choose a good section of $(K_X+L)$ or of $L$ and
we restrict to this section. We call this procedure
{\sl vertical}, respectively {\sl horizontal},
{\sl slicing}; in order to do this
we need the following (for
a proof see [A-W], (2.5) and (2.6)).
(We note that vertical slicing was already used in the proof of (1.2.1).)

\medskip \noindent
{\bf Lemma (1.3).}
{\sl Let $\f:X\ra Z$ be a good contraction of a smooth variety,
assume moreover that $Z$ is affine and $K_X+L=\O_X$.
\item{(1.3.1)} {\bf (Vertical slicing)} Assume that $X''\subset X$
is a non-trivial
divisor defined by a global function $h\in H^0(X,K_X+L)=H^0(X,\O_X)$.
Then for a general choice of $h$, $X''$ is smooth
and any section of $L$ on $X''$ extends to $X$.

\item{(1.3.2)} {\bf (Horizontal slicing)} Let $X'$ be a general divisor
from the linear system $\vert L \vert$. Then ,
outside of the base point locus of $\vert L\vert$, $X'$ is smooth and
any section of $L$ on $X'$ extends to $X$.
\itemitem{(i)} If $\f':=\f_{\vert X'}$, then $K_{X'}$ is
$\f'$-trivial
\itemitem{(ii)} Let $Z':=Spec(X',\O_{X'})$. If $\f$ is birational
then the induced map $Z'\ra Z$ is a
closed immersion. Therefore
the map $\f$ restricted to $X'$ has connected fibers.
\itemitem{}}

\medskip
The above lemma on horizontal slicing is particularly effective if we
can prove a relative base point free theorem for the line bundle $L$,
which means that the evaluation $\f^*\f_*L\ra L$ is surjective. The
next result is in this direction and it is a special case of the main
theorem of [A-W] (i.e. Theorem (5.1) in [A-W]).

\proclaim Proposition (1.3.3). {\bf Relative spannedness. }
Let $\f:X\ra Z$ be a good birational contraction
of a smooth $n$-fold. Assume that a fiber $F$ of $\f$ is of dimension
$\leq 2$. Then the evaluation morphism
$\f^*\f_*L\ra L$ is surjective at every point of $F$
(we say that $L$ is $\f$-spanned).

\medskip
\proclaim Proposition (1.3.4). In the same hypothesis of the proposition
(1.3.3) (or, more generally, of theorem 5.1 from [A-W])
the bundle $L$ is $\f$-very ample which means that there exists
an embedding $X\ra Z\times\P^N$ over $Z$ such that $L$ is the pull-back of
$\O(1)$.

\proof The proof is the same as the one of the theorem (5.1)
in [A-W]: in the hypothesis one can slice horizontally until fibers of
$\f$ are zero dimensional. In this case the $\f$-very ampleness is clear.
By (1.3.2) the sections of $L$ extend up from a divisor from $|L|$
so does the map to $Z\times \P^N$.

\medskip
A useful consequence of the relative very ampleness is an estimate on
the normal of a linear subspace of a fiber of a contraction:

\proclaim Lemma (1.3.5). {\rm (c.f.~[Ei])} Let $\f:X\ra Z$ be a good
contraction with a $\f$-very ample line bundle $L$. Assume that
$S\iso\P^r$ is contained in a fiber of $\f$ and $L_{|S}\iso\O(1)$.
Then the twisted conormal bundle $N^*_{S/X}(1)$ is spanned by global
sections.

\proof The inclusions $S\subset X\subset \P^N$ yield a surjection
$N^*_{S/\P^N}\ra N^*_{S/X}$. Since $N^*_{S/\P^N}\iso\O(-1)^{N-r}$ the lemma
follows.

\bigskip
Another fundamental property of good contractions we will use is the
existence of rational curves in fibers: through any point of the
contracted locus there passes a rational curve. A proof of this fine
property requires deformation theory. We only note that the vanishing
(1.2.2) implies that any 1-dimensional component of a fiber must be
$\P^1$ (see (2.1)).
\par

We will study the deformations of rational curves as well as
their chains and for this we need the following estimate on the
dimension of a component of the Hilbert scheme containing the class
of a curve $C \subset X$.

\definition \ (1.4). A proper curve $C$ is called smoothable if
there is an irreducible pointed variety $0\in T$ and a proper flat
family of curves $g: W \ra T$ such that $C= g^{-1}(0)$ and the generic
fiber of $g$ is smooth.
\medskip
There are two natural examples of reducible smoothable curves of genus 0
a tree and a bunch of rational curves.
A curve $C=\cup_i R_i$ is a (connected) tree of rational curves if:
\item {(i)}any $R_i$ is a smooth rational curve
\item {(ii)} $R_i$ intersects $\sum_{i-1}^{j = 1} R_J$ in
a single point which is an ordinary node of $C$.
\par\noindent The smoothing of a tree of rational curves is obtained
by modifying (blowing-up) the special fiber in the trivial family 
$\P^1\times T$.
\par
A bunch of $m$ rational curves is a projective cone over $m$ generic points
in $\P^{m-1}$. In other words such a bunch is a section of
$\S_m$ by a hyperplane which passes through the vertex of $\S_m$.
The smoothing of the bunch of rational curves is obtained by considering
a generic pencil of sections of $\S_m$ which contains the section in question.
Let us also note that if we attach to a one of the stems of a bunch of 
rational curves a tree (so that the meeting point is an ordinary node)
then the resulting curve is again smoothable.

\proclaim Proposition (1.4.1). Let $C$ be a proper curve
without embedded points.
Suppose that $f:C\ra X$ is an immersion
of $C$ into a smooth variety $X$ and that $C$ is smoothable.
Then any component of the Hilbert scheme containing $f(C)$ has dimension
$-K_X.C + (n-3) \chi (\O_C)$ at least.

We will apply the proposition in the case $C$ is a tree or a bunch of
rational curves; in particular we will have $\chi (\O_C) = 1$ (see the
section 4, in particular (4.5) and the next).

The proposition is proved
in the book of J. Koll\'ar (see [Ko],ch. II, theorem (1.14)).
Another
version of this result
for irreducible rational curves was used by Mori to prove
the existence of rational curves in fibers of good
contractions (see [Mo1]). It was also used to make relation between the
dimension
of a fiber and the dimension of the exceptional locus of $\f$ (see [Wi],
Theorem (1.1) and [Io], Theorem 0.4).

\bigskip \noindent
(1.5). Our inductive proofs will frequently lead from good contractions
to crepant contractions. Therefore we will use some results on crepant
contraction.  As an immediate application of the Lemma (1.2.1) we obtain

\proclaim Corollary (1.5.1). Any one dimensional component $F'$ of a fiber
of a good
or crepant contraction is a smooth rational curve and
$H^1(F',\J/\J^2)=0$, where $\J$ is the sheaf of ideals of $F'$.  If
$F$ is one dimensional fiber of a good or crepant contraction then
$H^1(F,\O_F)=0$ and thus the graph of $F$, with edges representing
its components and vertexes representing their incidence points, is
simply connected.

\remark (1.5.2). A complete description
of the one dimensional fibers of a good contraction of a smooth variety
$X$ will be given in (4.1). If $\f$ is creapant contraction of a smooth
variety $X$
then a more detailed and refined description
of the configuartion of the curves in $F$ can be given
if $n = 2$ or $n=3$;
some more subtle arguments are needed.
For $n=2$ the incidence of the curves is described by a dual graph which is
isomorphic to one of the following Dynkin diagrams: $A_n$, $D_n$, $E_6$,
$E_7$, $E_8$ (see for instance [B-P-V]).
If $n=3$ the incidence of the curves is described by a dual graph which is
isomorphic to one obtained by contracting
any subset (possibly the empty subset) of the $(-2)$-curves of
a Dynkin diagram ($A_n, D_n,E_6,E_7,E_8$) (this was proved in [Re]).

\medskip
We restrict our attention now to the case of a crepant contraction
of a smooth $3$-fold; this was studied first by M. Reid (see [Re]) and
subsequently by J. Koll\'ar (see [C-K-M]), D.Morrison-S.Katz (see [Ka-M]
and Y.Kawamata
(see [Ka2]).

In particular M.Reid (see [Re], section 1,
Theorem (1.4 iii) and its proof) proved, among other results,
the following:

\proclaim Proposition (1.5.3).  Let $\f: X \ra Z$
be a small crepant contraction  of a smooth $3$-fold $X$. Then
$Z$ has only a terminal-Gorenstein singularity or,
equivalently, $Z$ has only cDV singularity. If $p$  is a  singular
point of $Z$ and $H$ is a generic divisor through $p$
(therefore $p$ is a rational double point on $H$
by the definition of cDV singularity)
then $\f^{-1}H$ is non singular at the general point
of any curve
of $\f^{-1} (p)$. In particular $\f^{-1}H$ is normal and
$\f^{-1}H \ra H$ is a partial resolution of $H$.

\bigskip

\beginsection 2. Generalities on blow-ups and vector bundles.

In the present section we collect some general results concerning
blow-ups and vector bundles which we will use in the sequel.
Our notation is consistent with that of Hartshorne's book [Ha, Sect.II.7].

The first results are about blow-ups of locally complete intersections.

\proclaim Lemma (2.1). Let $S\subset X$ be a locally complete intersection
subvariety in a smooth variety $X$ defined by a sheaf of ideals
$\I_S$. Suppose that $\beta: \hat X\ra X$ is the blow-up of $X$ along
the subvariety $S$. Then the exceptional set $E=E(\beta)$ is a
Cartier divisor on $\hat X$ and $\beta_E:E\ra S$ a projective bundle
over $S$ isomorphic to $\P(\I_S/\I_S^2)$. Moreover $\beta_*\O_{\hat
X}(-E)=\I_S$.  If $S$ is connected then $E$ generates the Picard
group of $\hat X$ over $PicX$. Moreover $\hat X$ is Gorenstein and
$K_{\hat X}=\beta^*K_X+(dimX-dimS-1)E$.

\proof The variety $\hat X$ can be locally embedded into $X\times \P^{k-1}$,
where $k$ is the codimension of $S$. Indeed, if $f_1,\dots f_k$ are
the functions defining locally the ideal $\I_S$ then $\hat X$ is
defined in $X\times\P^{k-1}$ by equations $f_it_j=f_jt_i$, where
$t_i$ are coordinates in $\P^{k-1}$. Using this observation one may
verify all the assertions. We refer the reader to [EGA] and [Ha] for
this and further properties of blow-ups.

\bigskip
Let $X$ be a smooth variety of dimension $n\geq 3$.
Assume that $S_1$ and $S_2$ are two
codimension 2 smooth subvarieties of $X$ which meet transversaly along
a set $\Delta$ of dimension $n-3$. That is, in local coordinates
$(z_0,z_1,\dots,z_{n-1})$ in a neighbourhood of a given point $x\in \Delta$,
the subvarieties $S_1$ and $S_2$ are defined by, respectively,
$z_0=z_1=0$ and $z_0=z_2=0$.
Therefore the reducible subvariety $S:=S_1\cup S_2$ is defined locally by
two functions and in these coordinate system its equations are
$z_0=z_1z_2=0$.
In particular, the sheaf $N^*_S:=\I_{S_1\cup S_2}/\I^2_{S_1\cup S_2}$
is locally free rank 2 over $S_1\cup S_2$.

\proclaim Lemma (2.2).
There is an exact sequence of $\O_{S_1}$-modules:
$$0\raa (N^*_{S/X})_{|S_1}\raa N^*_{S_1/X}\raa N^*_{\Delta/S_2}\raa 0.$$

\proof
We have an embedding of sheaves of ideals $\I_S\hookrightarrow \I_{S_1}$
which is the identity outside of $S_2$.
For $x\in \Delta$ and local coordinates as above, this inclusion
can be expressed as $(z_0,z_1z_2)\subset (z_0,z_1)$.
Therefore over $\Delta$ the quotient $\I_{S_1}/(\I_S+\I^2_{S_1})$
is generated by the function $z_1$. On the other hand,
the inclusion $\Delta\subset S_2$  is related to
$\I_{S_2}\hookrightarrow \I_{\Delta}$
and thus to the inclusion $(z_0,z_2)\subset (z_0,z_1,z_2)$.
Therefore $\I_{\Delta}/(\I^2_\Delta+\I_{S_2})$ is generated by $z_1$ as well.

\medskip

The geometric meaning of the above sequence can be described as follows.
Let $\pi_1:\hat X_1\ra X$ be the blow-up of $X$ along
$S_1$ with the exceptional divisor $\hat E_1$. Then
$\hat E_1=\P(N^*_{S_1})$ and $N^*_{S_1}=(\pi_1)_*\O_{\hat E_1}(-\hat E_1)$.
Let $\hat S_2\subset \hat X_1$  be the strict transform
of $S_2$. We note that $\hat S_2\iso S_2$ and $\hat S_2$ meets $E_1$
along a section of the $\P^1$-bundle
$\pi_1:\hat\Delta=\pi_1^{-1}(\Delta)\ra
\Delta$. Let us call the section $\Delta_1$. We note that the section
$\Delta_1$ is associated to the surjective morphism
$(N^*_{S_1})_{|\Delta}\ra N^*_{\Delta /S_2}\ra 0$
of $\O_{\Delta}$-modules (dually: $\Delta_1$ parametrizes normal vectors
along which $S_2$ enters into $\Delta$).
Now we blow-up $\hat X_1$ along $\hat S_2$ and we call the result
$\tilde X_2$ (and the exceptional divisor $\tilde E_2$). The strict
transform of $\hat \Delta$, call it
$\tilde\Delta$, has now the normal whose restriction to any fiber of the
ruling $\tilde\Delta\ra\Delta$ is $\O(-1)\oplus\O(-1)$. Thus, we
see that $\tilde\Delta$ can be contracted now in $\tilde X_2$.
That is, we
have a contraction $\tilde X_2\ra\bar X$ over $X$
which is isomorphism outside of $\tilde\Delta$ and which contracts
$\tilde\Delta$ to the set isomorphic to $\Delta$ and any point of this
set is a quadric cone singularity on $\bar X$. The map $\bar\pi: \bar X\ra X$
is the blow-up of the ideal of $S=S_1\cup S_2$. The strict transform of
$\hat E_1$, call it $\bar E_1$, is again a $\P^1$-bundle over $S_1$
and ${N^*_{S}}_{|S_1}$ is the pushforward of
$\O_{\bar E_1}(-\bar E_1-\bar E_2)$.

Over $S_1$, however, the whole operation can be described as the blowing-up
of the section $\Delta_1$ in $E_1$ and then contracting
of the strict transform of $\hat\Delta$.
The birational map $\bar X\ra\hat X_1$ is associated to the inclusion
$\I_S\subset\I_{S_1}$ and also to the injection
$0\ra N^*_{S|S_1}\ra N^*_{S_1}$.
The cokernel of the injection is $\O_{\Delta_1}(-E_1)=N^*_{\Delta/S_2}$.
The details of
this geometric interpretation of ``vector bundle surgery'' are explained
in general by Maruyama in [Ma, pp.116--117].

\medskip

The above blow-up argument is useful to understand the subsequent situation.

\proclaim Lemma (2.3).
With the above notation
$$N^*_{\hat S_2/\hat X_1}\iso (N^*_S)_{|S_2}\otimes \O_{S_2}(\Delta).$$

\proof
Indeed, the exceptional set $\tilde E_2$ of the blow-up $\tilde X_2\ra
\hat X_1$ is not influenced by the contraction of $\tilde\Delta$, that is
its strict tranform $\bar E_2$ is the same as $\tilde E_2$. But $\tilde E_2$
and $\bar E_2$ are projectivisations of, respectively, left-hand and
right-hand side of the above equality; thus the equality is proved up to
the twisting. The twist follows by comparing the 1st Chern classes of
$N^*_{\hat S_2/\hat X_1}$ and ${N^*_S}_{|S_2}$; e.g. using the above lemma.
Namely, in the above argument we noticed that (changing the indices
appropriately)
$(N^*_{S/ X})_{|S_2}=(\bar\pi)_*\O_{\bar E_2}(-\bar E_1-\bar E_2)$
while $N^*_{\hat S_2/\hat X}=(\tilde\pi_2)_*\O_{\tilde E_2}(-\tilde E_2)$,
where $\bar\pi$ and $\tilde\pi_2$ are the appropriate blow-downs.
Another proof of the above equality can be obtained by calculation in
local coordinates of the blow-up $\hat X_1$. That is, we can write the
equations of all the needed sets in the
local (mixed: affine-homogenous) coordinates $((z_0,z_1,\dots, z_{n-1}),[\hat
z_0,\hat z_1])$ around each point of $\Delta$ and we can compute
the equality directly.

\bigskip
We conclude the part of the section devoted to blow-ups with a version of
Castelnuovo theorem, see [Ha, V.5.7].

\proclaim Proposition (2.4).
Let $\f :X \ra Z$ be a projective morphism from
a smooth variety $X$ onto a normal variety $Z$ with
connected fibers (a contraction). Suppose that $z \in Z$ is a point of $Z$ and
$F = \f^{-1} (z)$ is the geometric fiber over $z$ (see (1.0)).
which is locally complete intersection in $X$. Let
$\I$ denote the ideal of $F$ and $r := dimH^0(F,\I/\I^2)$.
Assume that for any positive integer $k$
$$\matrix{ H^1( F, S^k(\I/ \I^2)) = 0& \hbox{ and }&
H^0(F ,S^k(\I/ \I^2)) = S^k H^0(F ,(\I/ \I^2))}$$
then $z$ is a smooth point of $Z$ and $dimZ=r$.

\proof Since $\f_* \O_X = \O_Z$ we can apply the theorem on formal function
(see [Ha, III.11.1]) to describe the completion $\hat\O_z$ of the
local ring of $z\in Z$ $$\hat {\O}_z = \mathop{lim}_{\longleftarrow}
H^0(\hat F_k,\O_{F_k}),$$ where $\hat F_k$ is the closed subscheme of
$X$ defined by $\f^{-1}(m_z^k)\cdot\O_X$.  Since $F = \f^{-1} (z)$ it
follows that the sequence of ideals $\f^{-1}(m_z^k)\cdot\O_X$ is
cofinal (in the sense of [Ha, p. 194]) with the sequence of ideals
$\I^k$, so we may use schemes $F_k$ defined by $\I^k$ instead of
$\hat F_k$.
\par
We will prove for each $k$ that $H^0(F_k,\O_{F_k})$ is isomorphic to
a truncated power series ring $A_k =\C[[x_1,...,x_r]]/(x_1,...,x_r)^k$.
This will imply that $\hat {\O}_z\cong \C [[x_1,...,x_r]]$
and therefore that $z$ is a smooth point of $Z$.
\par
For $k=1$ we have that $H^0(F,\O_F)= \C$.
For $k>1$ we note that $\I^k/\I^{k+1}=S^k(\I/\I^2)$,
because $F$ is a local complete intersection.
Thus the claim can be proved by induction, using the cohomology
sequence associated to the exact sequence
$$0 \ra \I^k/\I^{k+1} \ra \O_{F_{k+1}} \ra \O_{F_k} \ra 0$$
and the cohomology hypothesis (see [Ha], proof of (V.5.7)).

\bigskip

Now we pass to the discussion of vector bundles. We start with a very simple
result for vector bundles on a tree of rational curves.

\proclaim Lemma (2.5). Let $C=\bigcup C_i$ be a
connected (possibly reducible) curve which is a tree
of rational curves, $C_i\iso\P^1$.
Suppose that $\E$ is a vector bundle over $C$. Then the following
conditions are equivalent:
\par{\sl
\item{(a)} The bundle $\E$ is nef; that is,
$\O_{\P(\E)}(1)$ is nef on $\P(\E)$
(in particular, for any $i$ and any surjective
morphism over the $i$-th component $\E_{|C_i}\ra \E'_i$ we have
$deg(\E'_i)\geq 0$).
\item{(b)} The bundle $\E$ is spanned by global sections at a generic
point of any component of C; that is, for any $i$ and a generic $x\in C_i$
the evaluation $\Gamma(C,\E)\ra\E_x$ is surjective.
\item{(c)} The bundle $\E$ is spanned by global sections at every
point of C; that is, for every $x\in C$
the evaluation $\Gamma(C,\E)\ra\E_x$ is surjective.
\item{}    }

\medskip
\proof The implications (c)$\Rightarrow$(b)$\Rightarrow$(a)
are clear so we shall prove (a)$\Rightarrow$(c).
We will proceed by induction with respect to the number of
irreducible components of the curve $C$.
The lemma is clearly true if $C\iso\P^1$.
Now the curve $C$ which contains $n\geq 2$ irreducible components
can be presented as the union $C'\cup C_n$ where $C'$ has $n-1$ components,
$C_n\iso \P^1$ and $C'$ and $C_n$ meet at one point, say $p$.
By the inductive assumption $\E_{C'}$ is globally generated by
$\Gamma(C',\E_{C'})$ and $\E_{C_n}$ is globally generated by
$\Gamma(C_n,\E_{C_n})$. The bundle $\E$ is obtained from
$\E_{C_n}$ and $\E_{C'}$ by the identification of their fibers
over $p$. Moreover $\Gamma(C,\E)$ is obtained from
$\Gamma(C',\E_{C'})$ and $\Gamma(C_n,\E_{C_n})$ by identification
of the values of sections over $p$. Now it is clear that the spannedness
of the restriction of $\E$ to both $C'$ and $C_n$ implies the spannedness
of $\E$. For example: to get a section with a prescribed value $v_x$ over
$x\in C'$ we choose first a section $s'\in\Gamma(C',\E_{C'})$
such that $s'(x)=v_x$ and then $s_n\in \Gamma(C_n,\E_{C_n})$
such that $s_n(p)=s'(p)$. Now $s'$ and $s_n$ glue to a global section
$s\in\Gamma(C,\E_{C})$ such that $s(x)=v_x$.

\bigskip
The next result is from [S-W1]:

\proclaim Lemma (2.6). Let $\E$ be a rank-2 vector bundle over $\P^2$.
Suppose that $\E$ is spanned by global sections and $0\leq c_1(\E)\leq 2$.
Then $\E$ is isomorphic to one of the bundles listed in
Table I in which we provide also its Chern classes together with
a description of the map $\P(\E)\ra \P(H^0(\P^2,\E)^*)$
associated to the evaluation of sections.

\midinsert{{\bf Table I}
{\sl{\baselineskip=18pt
\settabs\+ $c_1$ & $c_1$ \ \ \ &
$0\ra \O\ra T\P^2(-1)\oplus\O(1)\ra\E\ra 0$ \ \ \ &
description of $\P(\E)\ra \P(H^0(\P^2,\E)^*)$\cr
\+  $c_1$ & $c_2$   & {\sl description of bundle} $\E$&
{\sl description of} $\P(\E)\ra \P(H^0(\P^2,\E))$\cr
\hrule
\medskip
\+  0&0&$\E=\O\oplus\O$&$\P^2 \times \P^1\ra\P^1$\cr
\+  1&0&$\E=\O\oplus\O(1)$&blow-up of a point at  $\P^3$\cr
\+  1&1&$\E=T\P^2(-1)$& $\P^1$-bundle over $\P^2$  \cr
\+  2&0&$\E=\O\oplus\O(2)$&blow-up of vertex of a cone over
$\P^2\subset\P^5$  \cr
\+  2&2&$0\ra \O\ra T\P^2(-1)\oplus\O(1)\ra\E\ra 0$&blow-up of a line in
a smooth quadric  \cr
\+  2&3&$0\ra\O(-1)^{\oplus 2}\ra\O^{\oplus 4}\ra \E\ra 0$&
blow-up of a twisted rational curve in $\P^3$  \cr
\+  2&4&$0\ra\O(-2)\ra\O^{\oplus 3}\ra \E\ra 0$&conic bundle over $\P^2$\cr
\medskip
\hrule
}}
}\endinsert

\proclaim Lemma (2.7). Let $\E$ be a rank-2 vector bundle over $\P^2$.
If the splitting type of $\E$ on every line is the same, then $\E$ is
either a twist of the tangent bundle or it is decomposable.  If the
splitting type of $\E$ on each line is either $(-1,1)$ or $(-2,2)$
then either $\E$ is decomposable or it is not semi-stable (which
means that $\E(-1)$ has a section) with $c_2(\E)=0$. In the latter
case the section of $\E(-1)$ vanishes at one point only and thus $\E$
is in the following sequence $$0\raa \O(1)\raa\E\raa\I_x(-1)\raa
0.\eqno(2.7.0)$$

\proof The first part of the lemma is just a theorem of Van de Ven.
The second part follows from a theorem of Grauert---Muellich.
Namely, because of its general splitting type the bundle
$\E(-1)$ has a section which (again by the splitting type)
vanishes at one point at most (see [O-S-S]).

\proclaim Corollary (2.7.1). Suppose that $\E$ is a rank-2 vector bundle
over $\P^2$ such that $c_1(\E)=2$ and for any line $l\subset\P^2$ we have
$H^1(l,\E_l)=0$. Then either the general splitting of $\E$ is
$\O(1)\oplus\O(1)$ or $\E$ is decomposable, or it is not 
semistable as in (2.7.0).
\medskip

Now we pass to the discussion of vector bundles over quadrics.

\proclaim Lemma (2.8). Let $\E$ be a rank-2 vector bundle over a 2-dimensional
quadric $V$ which is either a smooth quadric $\Q_2 = \F_0$,
or a quadric cone $\S_2$, or $\P^2\cup\P^2$.
Suppose that $det(\E)$ is the restriction of $\O(1)$
from $\Q^4$ and that $\E$ is spanned by global sections.
Then, up to an automorphism of the quadric, we have either
$\O\oplus\O(1)$, or $\E\iso{\cal S}(1)$ or $\E$ is a pull-back of
$T\P^2(-1)$ via a double covering $V\ra\P^2$.

\proof The case of $\F_0$ is discussed in [S-W2]. The case of $\S_2$ can be done
similarly. We note only that since $\E$ is spanned then its
second Chern class $c_2$ as well as its Segre class $s_2=c_1^2-c_2=2-c_2$ have
to be non-negative and thus $c_2=0,\ 1,\ 2$. Because $\P(\E)$ is a Fano
variety one can apply vanishing and get $h^0(\E)=5-c_2$.

The case of $\P^2\cup\P^2$ can be easily done with the help of Van de Ven
theorem.
Namely, the restriction of $\E$ to each of the components is either
$\O\oplus\O(1)$ or $T\P^2(-1)$ and the glueing of these two bundles
along $\P^2\cap\P^2$ is unique as it follows from an easy argument concerning
extensions of automorphisms of these bundles from a line to $\P^2$.

\remark (2.9). Let $F$ be a 2-dimensional
quadric which is either a smooth quadric $\Q_2 = \F_0$, or a quadric
cone $\S_2$, or a reducible quadric $\P^2\cup\P^2$, and let $\E =
{\cal S}(1)$.  Then the map $\f: \P(\E)\ra \P(H^0(\P^2,\E)^*)$,
associated to the evaluation of sections, is a birational map to
$\P^3$. If $F = \F_0$ then $\f$ is the blow-up of two disjoint lines
in $\P^3$: up to a change of coordinates $[x_0,x_1,x_2,x_3]$ in
$\P^3$ the morphism $\f$ is the blow-up of the ideal
$(x_0x_2,x_0x_3,x_1x_2,x_1x_3)$. If $F = \S_2$ then
$\f$ is the blow-up of a double line associated to the ideal
$(x_0^2,x_1^2,x_0x_1,x_0x_2+x_1x_3)$.  If $F$ is the reducible
quadric then $\P(\E)$ is reducible and it consists of two components:
$\P(\E) = \P (\E_{|F_1}) \cup \P (\E_{|F_2}) = \P (T_{\P^2} (-1))
\cup \P (\O \oplus \O(1))$ and $\f_{|\P (T_{\P^2} (-1))}$
is a $\P^1$-bundle over $\P^2 \subset \P^3$ while $\f_{|\P (\O \oplus
\O(1))}$ is the blow up of a point $x \in \P^2\subset \P^3$.
We note also that, since each of the quadrics is a linear section of the
Grassmann variety of lines on $\P^3$, each of the projective bundles
has a natural interpretation in terms of incidence of lines in $\P^3$
(c.f.~3.4.0).

\medskip

Suppose now that $V$ is an irreducible quadric, that is either $\F_0$
or the quadric cone $\S_2$ with the resolution $\pi: \F_2\ra \S_2$ of
its vertex.  Line bundles over $\F_0$ are of the form
$\O(a_1,a_2)=\O(a_1f_1+a_2f_2)$, where $(f_1,f_2)$ are fibers of two
different ruling over $\P^1$ and $a_1$, $a_2$ are integers. A conic
on a quadric is an element of the $|\O(1)|$, or a linear section of
the surface embeded into $\P^3$. If $V=\F_0$ then $\O(1)=\O(1,1)$
and if $V=\S_2$ then $\pi^*\O(1)=\O(C_0+2f)$ where $C_0$ is the
exceptional divisor of $\pi$ and $f$ is a fiber of the ruling
$\F_2\ra \P^1$.

\proclaim Lemma (2.10.1). Let $\E$ be a rank 2-vector bundle over $\F_0$.
Suppose that the restriction of $\E$ to a generic conic splits into a
direct sum $\O(a)\oplus\O(b)$, where $a\geq b$. If $a-b\geq 2$ then
there exists a zero-dimensional subscheme $Z\subset\F_0$ such that
the bundle $\E$ is in the following exact sequence
$$0\raa\O(a_1,a_2)\raa \E\raa\I_Z\otimes \O(b_1,b_2)\raa 0,$$
where $a_1+a_2=a$ and $b_1+b_2=b$.

\proclaim Lemma (2.10.2). Let $\E$ be a rank 2-vector bundle over $\S_2$.
Suppose that the restriction of $\E$ to a generic conic splits into a
direct sum $\O(a)\oplus\O(b)$, where $a\geq b$. If $a-b\geq 2$ then
there exists a zero-dimensional subscheme $Z\subset\F_2$ such that
the bundle $\pi^*\E$ is in the following exact sequence
$$0\raa\O(a_0C_0+af)\raa\pi^*\E\raa\I_Z\otimes \O(b_0C_0+bf)\raa 0,$$
where $2(a_0+b_0)=a+b$.

\proof The proof of both lemmata depends on a criterion from [F-H-S],
Theorem (4.2).  Although the quoted theorem is formulated for smooth
varieties, it is not hard to see that the construction of a rank 1
subsheaf of $\E$ holds also in the case of the resolution of the
singular irreducible quadric. We note that the condition of Lemma 4.1
of [ibid] is satisfied. This is because the incidence variety $Z\ra
V$ of conics on the quadric $V$ is obtained from the flag manifold of
planes in $\P^3\supset V$ and one computes easily that the
restriction of $T_{Z/V}$ to the unique section $\hat C$ of $Z\ra V$
over a smooth conic $C$ is equal to $\O(-1)\oplus\O(-1)$.

\proclaim Lemma (2.11). Let $\E$ be a rank 2 vector bundle over an
irreducible quadric $V$ such that $det\E=\O(1)$.  Assume that the first
cohomology of the restriction of the bundle $\E$ to any conic is
zero.  Then either the restriction of $\E$ to a generic smooth conic
is $\O(1)\oplus\O(1)$ or one of the following is true:
\item{(a)} $V=\F_0$ and either
\itemitem{(i)} $\E$ is decomposable and (up to the change of the order
in the product $\P^1\times\P^1$) isomorphic to $\O(1,1)\oplus\O$ or
$\O(2,0)\oplus\O(-1,1)$, or $\O(2,1)\oplus\O(-1,0)$, or
\itemitem{(ii)} there exists a point $x\in V$ such that $\E$ fits in the
exact sequence $$0\raa\O(1,1)\raa\E\raa\I_x\ra 0.$$
\item{(b)} $V=\S_2$ and if $\E'=\pi^*\E$ then either
\itemitem{(i)} $\E'$ is decomposable and  isomorphic to $\O(C_0+2f)\oplus\O$ or
\itemitem{(i')} $\E'$ fits in one of the sequences
$$0\raa\O(2C_0+2f)\raa\E'\raa\O(-C_0)\raa 0\eqno(2.11.0)$$
 $$0\raa\O(2C_0+3f)\raa\E'\raa\O(-C_0-f)\raa 0$$
\itemitem{(ii)} there exists a point $x\in \F_2\setminus C_0$
such that $\E'$ fits in the exact sequence
$$0\raa\O(C_0+2f)\raa\E'\raa\I_x\ra 0.$$

\proof The assumption on vanishing of the cohomology implies that
the splitting type of $\E$ ony any smooth conic is either $(1,1)$ or
$(2,1)$, or $(3,-1)$ while the splitting type on any line is either
$(1,0)$ or $(2,-1)$. If the generic spliting type is not $(1,1)$ then
we can use the previous two lemmata to place $\E$ (or $\E'$) as a
middle term of an exact sequence. We moreover know that the
intersection of the invertible subsheaf $\L_1$ of $\E$ (resp. $\E'$),
provided by this sequence, with a conic is 2 or 3, while the intersection
with any line is $\leq 2$. In case of $\S_2$ we also know that
$\L_1.C_0\leq 0$. This gives a list of possible $\L_1$. Similarly we
find the condition on the degeneracy of the map $\L_1\ra\E$, that is on the
0-cycle $Z$ from the previous lemmata. If $V=\F_0$ then we note that
$length(Z)\leq 1$ because otherwise we would find a smooth conic containing
at least two of the points from $Z$ and the spliting type of $\E$ on such
a conic would not be admissible. Considering the restriction of $\E$ to lines
passing through $Z$ we find out that $Z\ne 0$ only if $\L_1=\O(1,1)$.
The argument for $V=\S_2$ is similar: considering lines passing through
$Z$ we find out that $Z\ne 0$ only if $\L_1=\O(C_0+2f)$, in such a case
however $Z\cap C_0=\emptyset$ because $\L_1.C_0=0$. Also, $Z$ consits of one
(reduced) point because we can consider the splitting type of $\E'$
on a conic passing through two points from $Z$.

\proclaim Corollary (2.11.1). Suppose that $\E$ is as in the above lemma and
its generic splitting type on conics is neither $(1,1)$ nor $(2,0)$, then
$h^0(V,\E\otimes\O(-1))-h^1(V,\E\otimes\O(-1))>0$.

\proof It follows immediately from the previous lemma.

\beginsection 3. Good contractions: examples.

In this section we will present several examples of good contractions of
a smooth projective variety $X$; we look especially
at the case in which all fibers have dimension $\leq 2$.

\medskip \noindent
(3.1) {\bf Projective bundles.}

\smallskip \noindent
Let $(F,\E)$ be a pair consisting of a smooth variety $F$ and a
numerically effective (for example spanned) vector bundle $\E$ such
that $-K_F-det\E$ is ample.  The bundle $\E$ is called Fano since its
projectivisation is a Fano manifold, see [S-W1, S-W2] and Table I.
Let then $X :=\P(\cal E\oplus\O)$ and let $\xi$ denote the
tautological line bundle on $X$, that is $\O_{\P(\E\oplus\O)}(1)$.
Let us consider the section of the projective bundle $X\ra F$
determined by the surjection ${\cal E}\oplus\O \ra \O \ra 0$, we will
call it again $F$. It is easy to check that $N^*_{F/X} = {\cal E}$.
Since $\xi$ is nef and $\xi-K_{\P(\E\oplus\O)}$ is ample, thus by
Kawamata-Shokurov contraction theorem it follows that $m\xi$ is base
point free for $m\gg 0$ and it defines a good contraction $\f: X \ra
Z$.  In particular $Z=Proj(\bigoplus_{m\geq 0}
H^0(F,S^m(\E\oplus\O)))$ and if $\E$ is spanned then $\f$ is the
connected part of the Stein factorsation of the map given by the
linear system $|\xi|$.  The map $\f$ contracts $F$ to a point.
\par

\example (3.1.1).
To get 4-dimensional birational contractions we consider pairs
$$\matrix{(\P^2,\ (T\P^2(-1)\oplus\O(1))/\O), &(\P^2,\ \O^{\oplus
4}/\O(-1)^{\oplus2})\cr (\P^1\times \P^1,\
\O(1,0)\oplus\O(0,1)),&(\P^2,\  \O(1)\oplus\O(1)).}$$
By the theorem of Leray and Hirsch one gets that $\xi^4 = c_1^2({\cal
E}\oplus \O) - c_2({\cal E}\oplus \O) > 0$.  This implies that the
divisor $\xi$ is big so that the map $\f$ is birational.  In each of
the above cases the map $\f$ will have exactly one 2-dimensional
fiber, namely $F$. If $(F,{\cal E})=(\P^2,\ \O(1)\oplus\O(1))$ then
the contraction is small while in each of the remaining cases it will
contract a divisor (see also (6.1) and (6.2)).  For the pair $(F,\cal
E)$ different from $(\P^1\times \P^1,\ \O(1,0)\oplus\O(0,1))$ the map
$\f$ is elementary. Let us also note that $X$ can be described
in terms of the rational map $\P(H^0(\E))\supset\f(X)- \ra F$. 
For example, if $\E=\O^{\oplus 4}/\O(-1)^{\oplus 2}$ then the map
is defined by a net of quadrics which contain the rational twisted
cubic in $\P^3$.
\par
Similarly, considering other bundles from (2.6) and (2.8) we get a series
of good contractions of fiber type with the special fiber $F$ equal to
$\P^2$ or $\F_0$. More precisely, to get an isolated 2-dimensional fiber 
one has to consider $T\P^2(-1)$ and $\O^{\oplus 3}/\O(-2)$ over $\P^2$
and the pullback of $T\P^2(-1)$ to $\F_0$.
\smallskip

\example (3.1.2). Using the results of [S-W2] one can produce examples
of contractions of a manifold $X=\P(\E\oplus\O)$, $dimX=n\geq 5$,
with an isolated 2-dimensional fiber $F\iso\P^2$. Birational
contractions can be obtained for $n=5,\ 6$ if we take $\E=\O^{\oplus
n}/\O(-1)^{\oplus 2}$ (which means that $\E(1)$ is the quotient of
$\O(1)^{\oplus n}$ divided by two general sections) or
$\E=T\P^2(-1)\oplus\O(1)$ for $n=5$.  Fiber type contractions (conic
fibrations with an isolated $\P^2$ fiber) will be produced this way
if we take $\E=\O^{\oplus n-1}/\O(-2)$ for $n\leq 7$.
Similarly, if we take $\E=\O^{\oplus 4}/\O(-1,-1)$ over
$\F_0=\P^1\times\P^1$ then $\P(\E\oplus\O)$ has a structure of a conic
fibration over $\P^4$ with an isolated 2 dimensional fiber equal
to $\F_0$.

\medskip \noindent
(3.2) {\bf Complete intersections in projective bundles}.

\smallskip \noindent
\example (3.2.1). Let $\Q_3$ be a smooth 3-dimensional quadric and let
${\cal S}$ be the spinor bundle over $\Q_3$.  (As described in
section 0 the bundle ${\cal S}$ is a rank two vector bundle which is
the pull back of the universal bundle under the inclusion $i: \Q_3
\ra Gr(1,3)$.) Let us consider the projectivization $M := 
\P(\O\oplus{\cal S}(1))$ with the projection $p: M\ra \Q_3$.  
By $\xi$ let us denote the line bundle $\O_{\P(\O\oplus {\cal
S}(1))}(1)$.  The complete linear system $|\xi|$ gives a contraction
of $M$ onto $\P^4$ such that the section corresponding to the trivial
factor of $\O\oplus{\cal S}(1) $ is contracted to a point by $\f$ and
all the other fibers of $\f$ are $\P^1$'s.

Let then $X$ be a smooth divisor in the
linear system $|\xi+p^*\O_{\Q^3}(1)|$. Since this linear system
is ample $Pic(X)$ is the same as $Pic(M)$ and thus $\f$ restricted
to $X$ --- call it again $\f$ --- is an elementary contraction.
Moreover, it is easy to see that $\f: X\ra \P^4$ is birational
and divisorial. Also, it will have a 2-dimensional fiber
$F$ which comes from a 3-dimensional fiber of $M\ra \P^4$,
thus $F$ is a quadric. For a general choice of $X$ the fiber
$F$ is a smooth quadric $\P^1\times\P^1$. Thus we can
realise the pair $(F, N_{F/X}) =(\P^1\times \P^1,\ \O(-1,0)\oplus\O(0,-1))$
in an elementary contraction (c.f.~the previous series of examples).

Moreover, we claim that there exists a smooth $X$ such that $F$ is
the quadric cone $\S^2$. Indeed, let us choose a section $\alpha\in
H^0(\Q^3,\O(1))$ which vanishes along a quadric cone with a vertex
at $w$ and a section $\beta\in H^0(\Q^3,{\cal S}(2))$ which does not
vanish at $w$. Then the section of $\xi+p^*\O(1)$ associated to the
section $(\alpha,\beta)$ of $\O(1)\oplus{\cal S}(2)$ is smooth along
the 2-dimensional fiber of the contraction. (This can be verified
locally --- as in the subsequent example.) Thus we can produce an
example of the a fiber and its normal being $(\S_2, {\cal
S}_{|\S_2})$.

\medskip
\example (3.2.2). Now we extend the above example to the codimension 2 
complete intersection in a projective bundle over the smooth 4
dimensional quadric $\Q_4$. This time we consider a spinor bundle
${\cal S}$ over the quadric $\Q_4$ and the projective bundle
$p: M=\P(\O\oplus{\cal S}(1))\ra \Q_4$.  Again, let $\xi$ denote the
relative hyperplane bundle --- the evaluation map associated to the
system $|\xi|$ is onto $\P^4$ and contract a section $Q_0$ to a point
$v\in \P^4$.  The quadric $\Q^4$ parametrizes planes in $\P^4$
containing $v$ and $M$ is the incidence variety of points on the
planes.

Again, we consider the linear system $\Lambda=|\xi+p^*\O(1)|$ over
$M$ and $X$ will be a complete intersection of two divisors from
$\Lambda$. The intersection $F:=X\cap Q_0$ is a linear section of the
4-dimesional quadric. The codimension 2 linear section of
$Q_0\iso\Q^4$ can be either a smooth quadric $\F_0=\P^1\times\P^1$ or
a quadric cone $\S_2$, or a union of two planes meeting along a line.
We will show that each of these cases can be realised so that the
complete intersection $X$ is smooth. The first two cases were dealt
with in the preceeding examples so we focus on the case of the two
$\P^2$s.  We distinguish the two planes calling them $P_1$ and
$P_2$, so that ${\cal S}(1)_{|P_1}=T\P^2(-1)$ and ${\cal
S}(1)_{|P_2}=\O\oplus\O(1)$.  In particular, over the line
$l:=P_1\cap P_2$ the bundle $\S(1)$ decomposes into $\O\oplus\O(1)$.

Sections of $\xi+p^*\O(1)$ are associated to sections of
$\O(1)\oplus{\cal S}(2)$; the decomposition of the later bundle
enables to write the sections in the form $(\alpha,\beta)$, where
$\alpha\in H^0(\Q^4,\O(1))$ and $\beta\in H^0(\Q^4,{\cal S}(2))$.
For a suitable trivialisation of ${\cal S}(2)$ (and thus of $\P({\cal
S}(1)\oplus\O)$) the zero locus of the section of $\xi+p^*\O(1)$
associated to $(\alpha,\beta)$ is given by the equation
$$\alpha(x)z_0+\beta^1(x)z_1+\beta^2(x)z_2=0$$ where $[z_0,z_1,z_2]$
are the associated uniform coordinates in the fiber of $p$ and $x$ is
a coordinate in $\Q^4$.  Moreover, we may assume that the
trivialisation of ${\cal S}(2)$ over $l$ coincides with the splitting
${\cal S}(2)_l=\O(1)\oplus\O(2)$ and therefore ${\beta^k}_{|l}\in
H^0(l,\O(k))$ for $k=1,\ 2$.  In these coordinates the section $Q_0$
is where $z_1=z_2=0$ and thus in affine coordinates around $Q_0$ the
zero locus is described by the equation
$$\alpha(x)+\beta^1(x)z_1+\beta^2(x)z_2=0.$$

Let $\alpha_1$ and $\alpha_2$ be the two sections of $\O_{\Q^4}(1)$
such that their common zero is the union of two $\P^2$'s,
i.e.~$P_1\cup P_2$. Then each of $\alpha_i$'s is a quadric cone with
the vertex on the line $P_1\cap P_2$.  The partial derivative
$\partial\alpha_i/\partial x_j$ vanishes only at the vertex of the
cone defined by the section $\alpha_i$.  Therefore, the Jacobi matrix
$(\partial\alpha_i/\partial x_j)$ (where $(x_j)$ are local
coordinates on $\Q^4$) is of rank one along $l=P_1\cap P_2$.  Now we
choose two sections, $\beta_i$ with $i=1,\ 2$, of ${\cal S}(2)$ such
that over $l$ they coincide with the splitting, that is
$\beta_{1|l}=(\beta^1_1,0)$ and $\beta_{2|l}=(0,\beta^2_2)$, where
$\beta^i_i\in H^0(l,\O(i))$.  (We can make such a choice because the
restriction map $H^0(\Q^4,{\cal S}(2))\ra H^0(l,{\cal S}(2)_{|l})$ is
surjective --- which is easy to verify e.g.~by considering
cohomology.) Moreover, we may assume that the zero set of the section
$\beta^i_i$ does not include the vertex of the cone defined by
$\alpha_i$.  And we consider two sections $f_i\in
H^0(M,\xi+p^*\O(1))\iso H^0(Q^4,\O(1)\oplus{\cal S}(2))$ which
locally are of the form $$\matrix{f_1:=\alpha_1(x)z_0+\beta_1z_1&\ \
& f_2:=\alpha_2(x)z_0+\beta_2z_2}$$ Now we can compute the matrix of
derivatives $(\partial f_i/\partial x_j,\partial f_i/\partial z_k)$
in local coordinates $(x_j,z_1,z_2)$ around $l\subset Q_0$.  The
result evaluated on $l$ is the matrix: $$\Big(\matrix{ {\partial
\alpha_1\over \partial x_j}&{\beta^1_1}&0\cr
{\partial \alpha_2\over \partial x_j}&0&{\beta^2_2}}
\Big)
$$
Because of our assumption on the zero locus of $\beta_i^i$ it follows that
the above matrix is of rank 2 everywhere on $l$. Therefore, the complete
intersection of divisors defined by $f_i$'s is smooth.

\medskip
{\bf Remark 1. }{\sl Discussion of normal bundles of $P_i$'s in $X$. }

The conormal of the total fiber $\I_{P_1\cup P_2}/\I^2_{P_1\cup P_2}$
is the restriction of the spinor bundle; in particular, its
restriction to $P_1$ is $T\P^2(-1)$ and to $P_2$ is $\O\oplus\O(1)$.
Using Lemma 2.2 one proves that the conormal of $P_1$ and $P_2$ is,
respectively, a stable bundle with $c_1=2$ $c_2=4$ which is spanned
outside of $l$ (see [S-W1]) and, respectively, an unstable,
semistable bundle with $c_1=2$, $c_2=3$. In both cases $l$ is the
unique jumping line with the splitting type $\O(-1)\oplus\O(3)$. (For
further discussion see Section 6.)

\medskip
{\bf Remark 2.}
{\sl Digression on the geometry of surfaces in $\P^4$.}

The above series of examples with special fiber a two dimensional
quadric can be described from the point of view of the geometry of
surfaces in $\P^4$. Namely, in each of the cases $X$ is a blow up of
$\P^4$ along a surface $S$ with an isolated singularity at the point
$v$.  Equivalently, $X$ is a graph of a rational map
$\P^4-\ra\Q^s\subset\P^{s+1}$ for $s=2,\ 3,\ 4$ (respectively for
(3.1.1), (3.2.1), (3.2.2)).  The resolution of the sheaf of ideals
$\I_S$ is as follows: $$0\raa \O(-s-2)\raa \O(-s-1)^{\oplus
(s+2)}\raa \O(-s)^{\oplus (s+2)}\raa\I_S\raa 0.$$ If $s=2$ then $S$
consists of two planes meeting transversally at $v$ and the above
sequence can be computed directly. In the other two cases this can be
computed by looking at $S$ as degeneracy locus of sections of
sheaves.  Namely, the projective bundle $M$ (over $\Q^3$ or $\Q^4$,
respectively) in which we consider the complete intersection has a
projective bundle structure over $\P^4$ outside of the point $v$. If
we pull the sheaf $\O(1)$ from the quadric to $M$ and then we push it
forward to $\P^4$ then the resulting sheaf, call it ${\cal F}$, will
have an isolated singularity at $v$ and outside of $v$ it will be
locally free. The sheaf ${\cal F}$ is of rank 2 or 3 and $S$ is cut
out by its 1 or 2 sections of ${\cal F}\otimes\O(1)$, respectively.
Moreover, since $M$ comes from a flag variety, one finds out easily
that the restriction of ${\cal F}$ to any linear $\P^3\subset\P^4$
which does not contain $v$ is equal to ${\cal N}(1)$ (${\cal N}$
denoting the well known null-correlation bundle on $\P^3$) or
$\Omega_{\P^3}(2)$, respectively. This provides enough information to
get the above resolution and more. In particular one finds that the
degree of $S$ is 5 and 9, and the ideal $I_S$ is defined by five
cubics and six quartics, respectively. Using {\sl CoCoa} program [CC]
for symbolic computations we found the following example of the ideal
$I_{S}=
(txy-x^2y-txz+x^2z-t^2u, tx^2-xyz+xz^2-xyu+tzu, t^3-ty^2xy^2-xyz+tz^2,
t^2x-xyz+xz^2-tyu, txy+t^2z-y^2z+z^3-y^2u)$,
where $[x,y,z,t,u]$ are homogeneus coordinates and $v=[0,0,0,0,1]$.
If we blow up the ideal $I_S$ in $\P^4$ we obtain a special fiber
of dimension two which is the quadric cone $\S_2$. The reader may easily
verify that the tangent cone of $S$ at $v$ gives a double line scheme.

\medskip \noindent
(3.3) {\bf A conic bundle with a special fiber ${\bf S}_3$.}
\par\noindent
Let  $\beta: V\ra\P^3$  be the blow-up of $\P^3$
at a point $x_0$ with the exceptional divisor $V_0$ and
another projection $\alpha: V\ra \P^2$ which makes
$V$ the $\P^1$ bundle $\P(\O\oplus\O(1))$.
Over $V$ we consider the pull-back bundle $\E:=\alpha^*(\O^4/\O(-1)^2)$
and its projectivization $\pi : W=\P(\E)\ra V$. Over $W$ we have the relative
hyperplane section bundle $\O_{\P(\E)}(1)$, which we denote by $\xi$.
The bundle $\xi$ generates $PicW$ together with $H:=(\alpha\pi)^{*}(\O(1))$
and $D:=(\beta\pi)^{*}(\O(1))$. Let us note that all the above three
line bundles are spanned and they define maps onto $\P^4$, $\P^2$ and
$\P^3$, respectively.
Moreover, we can compute that $-K_W=3\xi+2D$.

Let $W_0:=(\beta\pi)^{-1}(x_0)$,
$W_0$ is the unique effective divisor in $|D-H|$.
We note that $W_0$ is a 4-fold whose case was discussed in a previous series
of examples: it has a structure of a $\P^2$ bundle over $\P^2$
and another good contraction supported by the divisor $\xi$ which contracts
an irreducible divisor $Y_0\subset W_0$ to a cone ${\bf S}_3\subset\P^4$.
The divisor $Y_0$ is the unique effective divisor in $W_0$
equivalent to the restriction of $2\xi-H$.
A section of $W_0\ra\P^2$ is contracted to the vertex of ${\bf S}_3$,
call it $\Pi_0$. Let us note that $Y_0$ is smooth outside of $\Pi_0$.
The set $Y_0\subset W_0$ is contracted by the birational map supported by
$a\xi+bD$ where $a,\ b>0$.
Moreover, let us note that $N^*_{\Pi_0/W}=(\O^4/\O(-1)^{\oplus 2})\oplus\O(1)$.

\medskip

\claim.  $Y_0$ extends to a smooth divisor $Y\subset W$ such
that $Y\in |2\xi-H+mD|$ where $m\gg 0$.

\proof Let us note that the linear system $|2\xi-H+mD|$ is base point free
outside of $W_0$ for $m\gg 0$. On the other hand, from the sequence
$$0\ra (2\xi+(m-1)D)\ra(2\xi-H+mD)\ra(2\xi-h+mD)_{|W_0}\ra 0$$
and from Kawamata-Viehweg vanishing, it follows that $Y_0$ extends.
If $Y$ is an extension of $Y_0$ then since $Y_0=Y\cap W_0$
we have the smoothness
of $Y$ at the smooth points of $Y_0$, that is outside of $\Pi_0$.
Thus, since smoothnes is an open property, we will be done if we show
a divisor in $|2\xi-H+mD|$ which is smooth along $\Pi_0$.
To this end, note that it is enough to  take a reducible
$Y'=W_0+Y''$ such that the divisor $Y''\in |2\xi+(m-1)D|$ does not meet $\Pi_0$.
The existence of such $Y''$ is clear since $2\xi+(m-1)D$ is trivial on $\Pi_0$.

\medskip

Thus we can produce a smooth conic fibration
$\psi=(\beta\circ\pi)_{|Y}: Y\ra \P^3$
with a non-normal divisorial fiber $Y_0$ over $x_0$
such that $-K_Y=\xi+H+m'D$.
Now a good supporting divisor $\xi+D$ defines a divisorial
contraction of $\rho: Y\ra X\subset \P^4\times\P^3$ of $Y$ over $\P^3$.
The contraction $\rho$ contracts the divisor $Y_0$ to the cone
${\bf S}_3= X\cap(\P^4\times\{x_0\})$. Thus $X$ is smooth and the
resulting contraction $\f: X\ra \P^3$ is a conic fibration with
an exceptional fiber ${\bf S}_3$.

\medskip\noindent
{\bf Remark.} A different, simpler and very nice construction of a conic
bundle with the special fiber $\S_3$ was given by N.~Shepherd-Barron. The
example is reported in a paper of Y.~Kachi (see [Kac], example (11.6)).

\bigskip \noindent
(3.4) {\bf Blow-ups, blow-downs.} 
\par\noindent
A very convenient way to produce a conic
fibration is to alter another fibration. We will use this method to produce
non-elementary fibrations also with reducible fibers. The fibration
which will be the base of the construction is either a simple $\P^1$-bundle
or  a non-equidimensional 4-dimensional scroll  $\psi : Y\ra Z$
with an exceptional fiber
$V\iso \P^2$, as dealt with in [A-W, Remark (4.12)], and [B-W].

\example (3.4.0) A 4 dimensional scroll with an exceptional fiber $\P^2$.
\par\noindent
Our favorite example of such a scroll comes from incidence
construction (c.f.~3.2.2)): we set $Z:=\P^3$ and we fix a point
$v\in \P^3$. Then we consider the incidence variety $$Y:=
\{(z,\Pi)\in \P^3\times Grass(\P^2,\P^3)\colon z\in \Pi \hbox{ and }
v\in\Pi\}$$ with the projections $\psi\colon Y\ra Z=\P^3$ and $\pi\colon Y\ra\P^2$.
The $\pi$ makes $Y$ a projective bundle $\P(\O\oplus T\P^2(-1))$
(c.f.[A-W],(4.12))
while $\psi$ is the scroll with a unique 2-dimensional fiber $V\iso\P^2$
over $v$. If $(z_0,z_1,z_2)$ are coordinates in the affine neighbourhood
of $v=(0,0,0)$ and $[t_0,t_1,t_2]$ are homogenous coordinates in the
$\P^2$ then the equation of $Y$ in $\C^3\times\P^2$ is 
$t_0z_0+t_1z_1+t_2z_2=0$ (this is just a duality pairing between the
plane at infinity of $\C^3$ with its dual).

\medskip

\example (3.4.1) A conic fibration with exceptional fiber $\F_1$.
\par\noindent
Let $S\subset Y$ be a smooth surface meeting $V$ transversaly at one point
and the other fibers of $\psi$ at one point (transversaly!) at most.
Then the blow-up of $Y$ along $S$: $\alpha: X\ra Y$ with the morphism
$$\f:=\psi\circ\alpha\colon X\raa Z$$
is a conic bundle with discriminant divisor
$\Delta=\psi(S)$ and a special fiber over $v$ isomorphic to $\F_1$.
Let us note that $X$ admits another good contraction $\beta$ which
is a simple blow-down map of $X$ to a $\P^1$-bundle over $Z$.
In local coordinates $S=\{([t_0,t_1,t_2],(z_0,z_1,z_2)): z_0=t_1=t_2=0\}$,
$\beta\circ\alpha^{-1}$ is the projection from $S$ and its inverse 
$\alpha\circ\beta^{-1}$ can be described as follows:
$$\P^1\times\C^3\ni([\tau_1,\tau_2],(z_0,z_1,z_2))\mapsto
([-\tau_1z_1-\tau_2z_2,\tau_1z_0,\tau_2z_0],(z_0,z_1,z_2))\in V.$$
The exceptional set of $\alpha\circ\beta^{-1}$ is the simple blow-up
of the surface $\Delta$ at $v$, that is $\{([\tau_1,\tau_2],(z_0,z_1,z_2):
z_0=\tau_1z_1+\tau_2z_2=0\}$.
\medskip
\example (3.4.2) A conic fibration with exceptional fiber $\F_0$.
\par\noindent
Let us consider a rational map $\gamma: \P^1\times\C^3-\ra\P^1\times\C^3$
given by the formula
$$\gamma([\tau_1,\tau_2],(z_0,z_1,z_2))=
([\tau_1z_0+\tau_2z_2,\ \tau_1z_1+\tau_2z_0],(z_0,z_1,z_2)).$$
The inverse of $\gamma$ is $([\tau_1,\tau_2],(z_0,z_1,z_2))\mapsto
([\tau_1z_0-\tau_2z_2,\ -\tau_1z_1+\tau_2z_0],(z_0,z_1,z_2))$ and
the exceptional set of each of these two maps is a resolution
of the quadric singularity $\{z_0^2=z_1z_2\}$. If $X$ is the resolution
of $\gamma$ then it is smooth and the map $X\ra\C^3$ is a conic fibration
with the discriminant $\Delta=\{(z_0,z_1,z_2): z_0^2=z_1z_2\}$ 
and an exceptional 2-dimensional fiber equal to $\F_0$.
\medskip

\example (3.4.3) A conic fibration with exceptional fiber $\F_1\cup\P^2$.
\par\noindent
Let $S\subset Y$ be a smooth surface meeting $V$ along a line and
such that $\psi_{|S}$ is a blow-down of a $(-1)$ curve.  For example,
in the local coordinates $([t_0,t_1,t_2],(z_0,z_1,z_2))$, which we
introduced above, we can take $S$ given by equations $t_0=z_0=0$.
Then the blow-up of $Y$ along $S$, $\alpha: X\ra Y$, with morphism
$\f:=\psi\circ\alpha$ is a conic fibration with discriminant divisor
$\Delta=\psi(S)$ and special fiber over $v$ isomorphic to $\F_1\cup\P^2$.

Alternatively, $X$ is the closure of the graph of a rational map of $Y$:
$$Y\ni ([t_0,t_1,t_2],(z_0,z_1,z_2))\mapsto 
([t_0,t_1z_0,t_2z_0],(z_0,z_1,z_2))\in W$$
where $W\subset \P^2\times\C^3$ is given by an equation 
$t_0z_0^2+t_1z_1+t_2z_2=0$ and thus it has an isolated quadric cone 
singularity at $([1,0,0],(0,0,0))$. The contraction
$X\ra W$ has an isolated two dimensional fiber $\iso \P^2$ with normal 
$(T\P^2(-1)\oplus\O(1))/\O$. The exceptional set of the inversed rational
map is the smooth  surface $\{z_0=t_1=t_2=0\}\subset W\subset \P^2\times\C^3$.

\medskip

\example (3.4.4) A conic fibration with exceptional fiber $\F_0\cup\P^2$.
\par\noindent
A similar argument as in the previous example leads to a conic
fibration with an exceptional fiber  $\F_0\cup\P^2$. In this case, however,
we choose $\Delta$ having a quadric cone singularity at $v$ so that
the map $\psi_{|S}:S \ra \Delta$ is a contraction of a $(-2)$ curve.
Namely, let us consider $S$ given in $Y\subset\P^2\times\C^3$
by equations $t_1=t_0z_0+t_2z_2=t_2z_0+t_0z_1=0$ (c.f.~(3.4.2)). 
Then $S$ is a resolution of $\Delta=\{(z_0,z_1,z_2): z^2_0=z_1z_2\}$ 
and the blow-up of $Y$ along $S$
is a conic fibration over $Z$ with an exceptional fiber 
$\iso \F_0\cup\P^2$.

Similarly as before one can describe $X$ as a graph of a rational map.
Let us consider $\gamma: Y-\ra \P^1\times\C^3$ such that
$\gamma([t_0,t_1,t_2],(z_0,z_1,z_2))=
([t_0z_1+t_2z_0,t_1],(z_0,z_1,z_2))$. Then the inverse of $\gamma$ 
is as follows
$$([\tau_0,\tau_1], (z_0,z_1,z_2))\mapsto 
([\tau_0z_2+\tau_1z_0z_1, \tau_1(z_1z_2-z_0^2), -\tau_0z_0-\tau_1z_1^2],
(z_0,z_1,z_2)).$$
The exceptional set of $\gamma$ is the smooth surface $S$ which we 
have just described above, while the exceptional set of $\gamma^{-1}$ is 
a surface with an isolated singularity at $([0,1],(0,0,0))$ which is
of the cubic cone type. Thus $X$ has two elementary contractions:
a simple blow-down to $Y$ and the birational contraction to
$\P^1\times\C^3$ with an exceptional $\P^2$ whose conormal is
$\O^{\oplus 4}/\O(-1)^{\oplus 2}$.

\medskip
\example (3.4.5) Conic fibrations with (non-isolated)
2-dim fibers equal to $\P^2\cup\F_2$ and $\P^2\cup\F_1\cup\P^2$.
\par\noindent
In $Y$ or $W$, which we have discussed above,
we consider a smooth surface $S=\{z_0=t_1=z_2=0\}$.
Then the blow-up of either $Y$, or respectively, 
$W$ along $S$ is smooth and the
induced map $\f_Y:X_Y\ra\C^3$ (resp. $\f_W:X_W\ra\C^3$)
is a good contraction. The fibers of
these contractions are of the following type: 
$\f^{-1}_Y((0,0,0))=\P^2\cup\F_2$,
$\f^{-1}_W((0,0,0))=\P^2\cup\F_1\cup\P^2$,
$\f^{-1}_Y(0,z_1,0))=\f^{-1}_W(0,z_1,0)=\F_0$ for $z_1\ne 0$ 
and all the other fibers are $\P^1$. 
We note that the exceptional fiber is a limit of two dimensional 
fibers and apart of the strict tranform of the special fiber of the
initial scroll, it contains the specialization of the pair $(\F_0,C_0+2f)$.
We also note that in both cases the contraction  factors through
a small contraction which contracts the strict tranform of
the special fiber of the initial scroll.

\bigskip \noindent
(3.5) {\bf Double coverings.} Let $\psi: Y\ra Z'$ be a good contraction
of a smooth variety $Y$. Assume that $L$ is a $\psi$-ample line bundle
and $-K_Y-2L=\psi^*(L')$ for some line bundle $L'$ over $Z'$. If $B\in |2L|$
is a smooth divisor then we can construct a double covering $\pi:X\ra Y$
which is branched along $B$ (see for instance [B-P-V], pp.42--43);
the variety $X$ is then smooth
and $-K_X=\pi^*(L+\psi^*(L'))$ so that $-K_X$ is $\psi\circ\pi$-ample.
\par
If $\psi$ is of fiber type then fibers of $\f:=\psi\circ\pi$
are connected and thus $\f:X\ra Z'$ is a good contraction. If $\psi$
is birational then $\psi\circ\pi$ is generically $2:1$.
In this case however, the connected
part of the Stein factorisation of $\psi\circ\pi$, which we denote
by $\f:X\ra Z$, is a good contraction. The finite part of the Stein
factorisation $\pi':Z\ra Z'$ is a double covering branched along $\psi(B)$.
\par
Let us note that a similar construction with $L$ such that $K_X+L=\psi^*(L')$
leads from a good birational contraction $\psi$ to a crepant contraction
$\f:X\ra Z$. The first example concerns this case.

\example (3.5.1) A crepant elementary divisorial contraction of
a smooth 3-fold.

Let us consider a product $Y=\F_1\times\C$ with projections
$p_1: Y\ra \F_1$, $p_2:Y\ra\C$. $Y$ admits a good contraction
$\psi:Y\ra\P^2\times\C$ supported by $p_1^*(C_0+f)$ which is a simple
blow-down of the exceptional divisor $C_0\times\C$. Let $L:=p_1^*(C_0+2f)$.
We claim that there exists a smooth divisor $B\in |2L|$
 such that $B_t:=B\cap \F_1\times\{t\}$ is a smooth
curve of genus 2 for general $t$ and $B_0=C_0\cup C_1$ where
$C_1\in |C_0+4f|$ is a smooth rational curve meeting $C_0$ transversally at
3 points. This follows from a general

\proclaim Lemma (3.5.2). Let $\Lambda$ be a base-point-free linear system
on a smooth variety $X$. Let $D_0\in \Lambda$ be a divisor which is smooth
except finite number of points. Then there exists a linear pencil
of divisors $\{D_\lambda:\lambda\in\P^1\subset\Lambda\}$ which contains
$D_0$ and such that the divisor $D_{\P^1}:=\bigcup_\lambda D_\lambda$
in $X\times\P^1$ is smooth.

\proof Let $D_\Lambda\subset X\times\Lambda$ be the universal
divisor (incidence variety); locally $D_\Lambda$ is defined
by a function $f(x,\lambda)$ with variables $x$ and $\lambda$
being coordinates in $X$ and $\Lambda$, respectively.
Since $\Lambda$ is base-point-free it follows that $D_\lambda$ is non-singular
(it has projective bundle structure over $X$) and thus the vector
of partial derivatives $(\partial f/\partial x,\partial f/\partial\lambda)$
is nowhere zero on $D_\Lambda$. If $D_{\lambda_0}$ is singular at
$x_0$ then $\partial f/\partial x$ vanishes at $(x_0,\lambda_0)$. However
in such a situation
$\partial f/\partial\lambda$ does not vanish so thus if we choose
a linear pencil $D_\lambda$ which contains $\lambda_0$ but is not
contained in the kernel of $\partial f/\partial\lambda$ then the resulting
divisor $\bigcup_\lambda D_\lambda$ will be smooth at $(x_0,\lambda_0)$.
Now the lemma follows easily.
\medskip

\noindent{\sl Example (3.5.1), continued.}
Let $\f: X\ra Z$ be a crepant contraction obtained from
$\psi: Y\ra \P^2\times\C$ with a double covering
$\pi:X\ra Y$ branched along $B$.
Note that $Z$ admits a morphism
onto $\C$ so that fiberwise we have a family of maps
$\f_t: X_t\ra Z_t$ parametrized by $t\in \C$.
For a general $t$ the surface $Z_t$ is a double cover of $\P^2$
branched along a quartic which has a simple double point
and $\f_t:X_t\ra Z_t$ is a resolution of the resulting $A_1$ singularity.
Let $E$ denotes the exceptional divisor of the map $\f$ and by $E_t$
let us denote the reduced fiber $E\cap X_t$. We claim that
\item{(a)} for any $t$, $E_t\iso\P^1$, $E.E_t=-2$ for a general
$t$ and $E.E_0=-1$,
\item{(b)} $E$ is non-normal along $E_0$ and smooth elsewhere,
\item{(c)} the normal bundle of $E_t$ in $X$ is $\O\oplus\O(-2)$
for a general $t$ and $\O(1)\oplus\O(-3)$ for $t=0$.
\item{}

\par
The first two properties follow easily from the construction while the third
one comes from

\proclaim Lemma (3.5.3). Let $\f:X\ra Z$ be a crepant divisorial contraction of
a smooth 3-fold with the exceptional irreducible divisor $E$
which is contracted to a curve $C\subset Z$. Let $0$ be a fixed point on $C$.
For any point $t\in C$ let $E_t$ denote the fiber $\f^{-1}(t)$ with the reduced
structure. Then the property (a) above yields property (c).

\proof We are to prove that the normal of $E_0$ in $X$ is $\O(1)\oplus\O(-3)$.
By (5.6.2) it is enough to exclude possibility that the normal is either
$\O(-1)\oplus\O(-1)$ or $\O\oplus\O(-2)$. Let us blow-up $X$ along
$E_0$ and call the resulting variety $\hat X$, the exceptional
divisor by $A$ and the strict transform of $E$ by $\hat E$. In $\hat X$
we have a family of effective 1-cycles $E_t$ parametrized by $C\setminus \{0\}$.
By $\hat E_0$ let us call the limit cycle of $E_t$ as $t\ra 0$.
Then $\hat E_0$ is an effective cycle, it is supported on the set
$\hat E\cap A$ and $A.\hat E_0=0$. This last equality alone excludes
the possibility of the normal $\O(-1)\oplus\O(-1)$ --- indeed, in this case
the divisor $-A$ is ample on $A$, so it has positive
intersection with any effective
cycle on $A$. If the normal of $E_0$ were $\O\oplus\O(-2)$ then $-A$ would be
 nef on $A=F_2$ and $\hat E_0$ would be supported on the unique curve
$C_0\subset F_2$ whose intersection with $-A$ is zero. But then
$\hat E_{|A}$ would be supported on $C_0$ so that
$$\gamma C_0=\hat E_{|A}=-\alpha A+(E.E_0)f=
\alpha C_0+(2\alpha+E.E_0)\cdot f$$
where $\alpha$ and $\gamma$ would be positive integers.
Since $E.E_0=-1$ this is impossible.
\bigskip

\noindent
\example (3.5.4). Divisorial elementary contraction of a 4-fold with
quadric fibers.

This is a 4-dimensional version of the previous example: let $V_1\ra\P^3$
be the blow-up of $\P^3$ at one point, by $S_0$ let us denote the exceptional
divisor of the blow-up. $V_1$ has a $\P^1$-bundle structure,
$V_1=\P(\O\oplus\O(1))\ra \P^2$, by $H$ let us denote the pull-back of the line.
Let us consider a product $Y:=V_1\times\C$ with projections $p_1$, $p_2$
onto factors. Over $Y$ we have a line bundle $L:=p_1^*(S_0+2H)$.
The data consisting of the contraction morphism $\psi: Y\ra \P^3\times\C$,
and $\psi$-ample line bundle $L$ can be plugged to the construction
described above as soon as we provide a smooth divisor $B$ in $|2L|$.
The resulting elementary contraction is divisorial and the exceptional
divisor $E$ is contracted to $\C$. From the theory of 3-dimensional
good contration [Mo1] we know that a general fiber $E_t$ of $E\ra \C$
is either a smooth quadric or a quadric cone. The singularity of a special
fiber $E_0$ will depend on the singularity of $B_0:=B\cap(S_0\times\{0\})$.
In particular, if we find a smooth $B$ such that this intersection
is a double line then $E_0$ will be a reducible quadric, that is, a union
of two planes meeting along a line. The construction of such a $B$ can
be done as follows: first, using the notation introduced in [Wi, p.154]
an arguing as in [ibid, pp.156--157] we can prove that there exists
a divisor $B_0$ on $V_1$ such that: $B_0\in 2S_0+4H$, $B_0\cap S_0$
is a double line and $B_0$ has only isolated singularities.
Then using lemma (3.5.2) we can extend $B_0$ to a smooth $B$ on $V_1\times\C$.

\medskip

\example (3.5.5). Conic fibrations.
\par \noindent
A similar argument will enable us to construct
different conic fibration: also non-equidi\-mensio\-nal and also with
reducible 2-dimensional fibers. If we use a $\P^1$ bundle as the base
of our covering then the resulting contraction will be equidimensional 
i.e.~a conic bundle. If however, we begin with the scroll
from (3.4) then we get a contraction with an isolated
2 dimensional fiber.
Moreover, according to the choice
of the branching divisor, the special fiber will be one of the following:
\item{(i)} a smooth quadric if the branching divisor intersects the
2-dimensional
fiber (i.e.~$\P^2$) along a smooth conic,
\item{(ii)} a quadric cone if the intersection is a reducible conic,
\item{(iii)} two planes if the intersection is a double line.

It can be shown that there exists no smooth branch divisor which contains
the 2-dimensional fiber (so that the case of a "double $\P^2$" can not
be produced this way).

\bigskip\noindent
(3.6) {\bf Toric examples.} An especially nice class of examples comes with
the geometry of toric varieties. We refer the reader to Oda [Od],
Danilov [Da] or Fulton [Fl] for the language and notation of toric geometry.

\example (3.6.1).  Blow-up of two transversal planes meeting at a point
(see also (3.1) and the Remark 2 in (3.2)).

The blow-up of two transveral planes
in $\C^4$ can be realised as follows. Let $N_\R$ be a 4-dimensional
real vector space with a basis $\{ e_1,\ e_2,\ e_3,\ e_4\}$
and let $N=\Z e_1+ \Z e_2 +\Z e_3 +\Z e_4$ be a lattice in $N_\R$.
Then we can subdivide the cone spanned by the basis by adding
two vertices $u:=e_1+e_2$, $v:=e_3+e_4$ and considering the following
4 cones
$$\matrix{\langle u, v, e_1, e_3\rangle,&\langle u, v, e_1, e_4\rangle,&
\langle u, v, e_2, e_3\rangle,&\langle u, v, e_2, e_4\rangle}$$
where $\langle\dots\rangle$ denotes the cone spanned by appropriate vectors.
The resulting fan describes a toric variety $X$ together with
a contraction onto $\C^4$ which is the blow-up along two transversal planes
defined (in suitable coordinate system)
by equations $z_1=z_2=0$ and $z_3=z_4=0$.

The projection of the space $N_R$ along the linear subspace spanned on
$u$ and $v$ given by a matrix
$$\big(\matrix{1&-1&0&0\cr 0&0&1&-1}\big)$$
maps the above defined fan onto a fan giving $\P^1\times\P^1$.
Actually, the reader may want to verify the resulting
construction describes the closure of the the graph of the rational map
$\C^4 \ra \P^1\times \P^1$ defined as follows
$$(z_1,z_2,z_3,z_4) \ra ([z_1,z_2],\ [z_3,z_4])$$
The map onto $\P^1\times\P^1$ makes $X$ the total space of the bundle
$\O(-1,0)\oplus\O(0,-1)$ as described in (3.1).

\medskip

\example (3.6.2). Conic fibration which a 2-dimensional fiber
which consists of two planes meeting in one point. (A non-toric version
of this example was presented in [Kac]; recently J.~W{\l}odarczyk exdended
this example to higher dimensions.)

Again,
we consider a 4-dimensional vector space $N_\R$ with a fixed basis
$\{ e_1,\ e_2,\ e_3,\ e_4\}$ and the lattice $N$ as above.
 We take the following vectors in $N_\R$:
$$\matrix{v_1=-e_1-e_3,& v_2=-e_2-e_3-e_4}$$
and consider a fan in $N_\R$ with vertices
in $e_1,\ v_1,\ e_2,\ e_4,\ -e_4$ containing the following 5 simpleses:
$$\matrix{\langle e_1, v_1, e_2, e_4\rangle, &
\langle e_1, v_1, v_2, e_4 \rangle,  &
\langle e_1, v_1, e_2, v_2 \rangle , &
\langle e_1, e_2, v_2, -e_4\rangle ,  &
\langle v_1, e_2, v_2, -e_4\rangle .}$$
The toric variety associated to this fan is smooth.
Now we consider the projection of $N_\R$ along the
last coordinate to a 3-dimensional
vector space $N'_\R$ with a basis $\{ e'_1,\ e'_2,\ e'_3\}$.
The projection maps the above fan into a cone spanned by
$e'_1,\ e'_2,\ -e'_1-e'_3,\ -e'_2-e'_3$,
which is a fan of the affine quadric cone (i.e.~the affine cone over
$\P^1\times\P^1\subset\P^3$).
We consider the induced
map of toric varieties $X\ra X'$. All fibers of this map except the one
over the vertex of the cone are $\P^1$'s.
The fiber over the vertex of the quadric cone consits of two 2-dimensional
orbits associated to simpleses $\langle e_,\ v_1\rangle$ and
$\langle e_2,\ v_2\rangle$. Each of these two components is
isomorphic to $\P^2$ and they meet at a point which is the orbit associated to
$\langle e_1,\ v_1,\ e_2,\ v_2\rangle$.
Moreover, let us note that
the divisors $D_1$ and $D_2$ associated to cones $\langle e_4\rangle$  and
$\langle -e_4\rangle$ map to $X'$ so that we get the usual toric flop
$D_1\ra X'\leftarrow D_2$.

\beginsection 4. Geometric fiber.

(4.0) Let $\f : X \ra Z$ be a contraction of a smooth $n$-fold $X$.
We assume that $\f$ is either a good contraction or a crepant contraction.
If $\f$ is a good contraction we fix a relatively ample line bundle
$L:=-K_X$. In the present section we want to study the geometric structure
of a positive dimensional fiber $F$ of $\f$. The fiber $F$ will be considered
with a reduced structure. We will say that $F$ is an {\sl isolated} fiber
of dimension $k$
if $\hbox{dim}F=k$ and all fibers of $\f$ in some neighbourhood of $F$ have
dimension
$< k$.

Let us start with the following well known description of 1-dimensional fibers
of good contractions (see [Mo1] and [An]).

\medskip
\noindent
{\bf Proposition (4.1)}
{\sl Let $\f:X\ra Z$ be a good contraction of a smooth variety. Suppose that
a fiber $F$ of $\f$ contains an irreducible component of dimension $1$;
then $F$ is of pure dimension $1$ and all components of $F$ are smooth
rational curves.
\item{(1)} If $\f$ is birational then
$F$ is irreducible and it is a line relatively to $L$, moreover the 
scheme theoretic fiber structure on $F$ is reduced.
\item{(2)} If $\f$ is of fiber type then $F$ is a conic relatively to $L$,
that is:
either
\itemitem{(i)} $F$ is a smooth $\P^1$ and $L\cdot F=2$, or
\itemitem{(ii)} $F$ is a union of two smooth rational curves meeting
at one point and each of these curves is a line with respect to $L$, or
\itemitem{(iii)} $F$ is a smooth $\P^1$, $L\cdot F=1$ and the fiber
structure on $F$
is of multiplicity 2 (a reduced conic).
\item{} In the cases (i) and (ii) the fiber structure is reduced.
\par
\noindent In each one of the above cases the variety $Z$ is smooth at $\f(Z)$.}
\medskip

\proof
Let $F$ be the fiber in question. We are to prove first that $F$ is of pure
dimension 1. Let us  consider an irreducible 1-dimensional component
$C$ of $F$. Lemma (1.2.1) implies that
$h^1(\O_C)= 0$ and thus $C\iso \P^1$.
By deformation arguments (see (1.4.1)) we have that
$dim_{[C]}Hilb(X) \geq  n-3-K_X.C = n-3 + L.C$
and therefore the deformations of $C$
sweep at least a divisor. Moreover if $L.C \geq 2$ then
$C$ must move in a $n-1$ dimensional family and therefore $\f$ cannot be
birational.
On the other hand the deformation locus meets the other
components of $F$ along a subset of $C$ and thus
the components of $F$ which meet $C$ have to intersect the deformation locus
at a finite number of points hence they are of dimension $1$.
Because the fiber is connected this applies to any
component of $F$ and thus we have that $F$ is of pure dimension 1.
The configuration of curves in the fiber was described in the proposition
(1.5.1).
From deformation theory, as above
using (1.4.1), we know that $F$ moves in a family
of dimension $L\cdot F+n-3$ at least. Thus we see that $L\cdot F= 1$ in the
birational case
and $L\cdot F\leq 2$ in the fiber type case. The description of the
morphism (smoothness
of $Z$ fiber, structure of non-reduced fibers) requires either studying
Hilbert scheme
of fibers (we do not consider it now,
see [An]) or one can use the results in (5.6) and in (2.4).

\bigskip \noindent
(4.2) We now pass to the case where $F$ has dimension two. That is,
$\f:X\ra Z$ is a good contraction
of a smooth $n$-fold with a fiber $F$ of (pure) dimension $2$.
Since the target $Z$ may be assumed affine and it can be shrunk, if
necessary, we may assume that all fibers of $\f$ are of dimension $\leq 2$.
As in (4.0) $L: = -K_X$ and we assume it is $\f$-spanned by global sections;
this is true in the birational case by the proposition (1.3.3)
and in the conic fibration case if $n\leq 4$ by a result of Kachi-Kawamata
(see [Kac]). Let us also note that some of the subsequent results which
depend on the vanishing (1.2) remain true if we allow $X$ to have log terminal
singularities. This concerns e.g.~(4.2.1) and (4.3.3).

\proclaim Proposition (4.2.1).
In the above hypothesis any  component $F'$ of $F$ is normal.
The pair $(F',L_{|F'})$ has Fujita $\Delta$-genus $0$
and it is among the following:
\item{(1)} $(\P^2,\O(e))$, with $e = 1,2$,
\item{(2)} $(\F_r, C_0 + kf)$ with $k \geq r+1$, $r\geq 0$,
\item{(3)} $(\S_r,\O_{\S_r}(1))$ with $r\geq 2$.

\proof The line bundle $L'=L_{|F'}$ is base point free and for any section $C$
of $L_{|F'}$ we have $g(C):=h^1(C,\O_C) = 0$. This is a
consequence of lemma (1.2.2) since $C= F'\cap D$, where $D$ is is a section
of $L$ (see (1.3.2)). Thus the proposition follows
from the following non-normal version of 
a very well known characterization of projective
normal surfaces with sectional genus $0$ (see [Fu]).

\proclaim Proposition (4.3). Let $F'$ be an irreducible (reduced) variety
of dimension $2$ and let $L'$ be
an ample and spanned line bundle such that for {\bf any} $C\in |L'|$ we have
$g(C) =: h^1(\O_{|C}) = 0$. Then $F'$ is normal, $L'$ is very ample
and the pair $(F',L')$ has delta genus zero, $\Delta (F',L') = 0$
thus the pair $(F',L')$ is one of the pairs in the proposition (2.3).

\smallskip
\remark (4.3.1). The proposition is not true if we only assume that
$g(C) = 0$ for a general $C$. For instance take the variety
obtained by identifying two points on $\P^n$
(remark at p.~30 of [Fu]).
\par
A similar result holds for varieties of higher dimension
where $C$ is any curve obtained as intersection of
$(n-1)$ elements of $|L'|$.

\medskip
We need the following lemma.

\proclaim Lemma (4.3.2).
Let $F'$ be as in the proposition
and let $L'$ be a line bundle on $F'$ such that the zero locus of
a general section of $L'$ is reduced and connected.
Assume that through a point $x\in F'$ there is a section $C\in |L'|$
which is generically reduced, connected and such that $h^1(\O_C) =
h^1(O_{C'})$, where $C'$ is a general section of $|L'|$. Then $x$ is a
Cohen-Macaulay point of $F'$.

\proof Notice that $\chi(\O_C) = \chi(\O_C')$; this follows
from the exact sequence
$$0 \ra L'^{-1} \ra \O_{F'} \ra \O_C \ra 0,$$
true for every $C\in |L'|$, and the fact that $\chi(\O_{F'})$ and
$\chi(L'^{-1})$ does not depend on $C$.
This and the hypothesis imply that
$h^0(\O_C) = h^0(\O_{C'}) = 1$.
We then consider the exact sequence
$$0 \ra S \ra \O_C \ra \O_{C_{red}} \ra 0,$$
where $S$ is a skycraper sheaf supported on the non
Cohen-Macaulay points of $C$; since  $h^0(\O_C) = h^0(\O_{C_{red}}) = 1$
we have that $S= 0$ and therefore $C$ is Cohen-Macaulay. Since
$C$ is a Cartier divisor every point of $C$ is a Cohen-Macaulay
point of $F'$.

\medskip \noindent
{\bf Proof of proposition (4.2.1).} We are to prove that $F'$ is
normal then the rest will follow from Fujita's result.
Since $L'$ is ample and spanned by the Bertini
theorem we have through every point
$x\in F'$ a section $C\in |L'|$
which is generically reduced and connected. By hypothesis
all $C \in |L'|$ have genus $0$ and therefore we can apply the lemma
and say that $F'$ has Cohen-Macaulay singularities.
\par
On the other hand a general element of $|L'|$, being
irreducible and reduced and of sectional genus $0$, is a
smooth rational curve. Therefore $F'$ itself is smooth
along the point of a general $C$; since $L'$ is ample this implies that
$F'$ is smooth in codimension $2$.
\par
By Serre criterion $F'$ is therefore normal.
\par
To conclude the proof we apply results of T.~Fujita;
in the subsequnet lines we refer to his book [Fu].
First by the proposition (3.4) it follows that
$\Delta (F',L') = 0$; secondly by the corollary (4.12) we have
that $L'$ is very ample. Finally applying the theorem (5.15) and
the remark (5.16) we have the complete proposition.
\medskip
We note that the argument applied in the course of the proof of (4.3.2)
to a component $F'$ can be actually used for the whole fiber $F$ and it
yields the following
\proclaim Lemma (4.3.3). Let $F$ be a two dimensional fiber of a good
contration as in (4.2). Then $F$ is Cohen-Macauley unless the zero locus of
a general section $\in |L_F|$ is disconnected. 

The example of non-Cohen-Macaulay fiber is obtained in (3.6.2):
the meeting point of the two components of the fiber is its unique non-C-M point. 
If $\f$ is birational then the hyperplane section of $F$ is connected
which follows from the subsequent result.

\proclaim Lemma (4.4) (Horizontal slicing).
Assume that $\f$ is a good contraction as in (4.2).
Let $X'$ be a general section of $|L|$
not containing any component of the special fiber $F$.
Consider the restriction of $\f$ to $X'$,
$\f_{|X'}: X' \ra \f(X')$, and let
$$X'\har{\f'}Z'\raa Z$$
be the Stein factorization of $\f_{|X'}$.
Then $X'$ is a smooth $(n-1)$-dimensional variety and
$\f'$ is a crepant (birational) contraction with at most 1-dimensional
fibers onto a normal variety $Z'$.
If $\f$ is birational, then $\f'=\f_{|X'}$. If $\f$ is of fiber type than
$Z'\ra Z$ is a double covering.
Moreover for $n\geq 4$  non trivial fibers of $\f':X' \ra Z'$
are parametrized on $Z'$ by a subvariety of dimension $\geq (n-4)$.
If $n=4$ and $F$ is an isolated 2-dimensional fiber of $\f$ 
then $X'$ may be chosen so that
non trivial fibers of $\f'$ are isolated.

\proof The first part of the lemma follows from (1.3.2) and
adjunction formula. Let $E'$ be the exceptional locus of $\f'$.
We claim that $dim \f'(E') \geq (n-4)$.
Indeed, take any rational
curve in the fiber of $\f'$. Since
$K_{X'}.C = 0$ then from deformation theory (1.4.1)
we have that $C$ moves at least in a $(n-1)-3$
family. The last part is obtained by a standard dimension counting.
\bigskip

The deformation argument which we have just applied implies that the
hyperplane section of the fiber $F$ (or its part, if it is reducible)
moves in a nontrivial family, if only $n\geq 5$. Thus, knowing the
classification of good contractions for $n\leq 4$ we can estimate the
degree of $F$ with respect to $L$. The idea of moving the rational
curves and estimating the exceptional locus of a good or crepant contraction
was used in the proof of the length-fiber-locus inequality in [Wi, (1.1)],
(see also [Io, 0.4]) which related a general non-trivial fiber of a good
contracion with the dimension of its exceptional locus. 

A similar idea can be used to distinguish the locus of fibers of different
dimension. Namely, let us set consider the locus of fibers of dimension $k$.
More precisely, for $k\geq 1$, let $E_k(\f)$ denote the closure of the set 
$\{x\in X:\hbox{dim}\f^{-1}(\f(x))=k\}$.
Let $E_k$ be an irreducible component of $E_k(\f)$. Then either
$E_k$ is an irreducible component of the exceptional locus $E(\f)$
of $\f$, or it is contained in $E_m(\f)$ for some $m<k$. 

\proclaim Proposition (4.5). Let $E_2$ be an irreducible 
component of $E_2(\f)$. Then 
$\hbox{dim}E_2\geq n-2$ if $E_2$ is also a component of $E(\f)$. Otherwise
\item{(i)} $\hbox{dim}E_2\geq n-4$ if $\f$ is birational,
\item{(ii)} $\hbox{dim}E_2\geq n-5$ if $dimX-dimZ=1$.
\item{}(All the above estimates are best possible, see (3.1).)

\proof The first part of the lemma is just [Wi, (1.1)].
Thus we may assume that $E_2\subset E_1(\f)$. Moreover, we may
asume that $E_2$ is just a component of a 2-dimensional isolated
fiber of $\f$. Indeed, we can consider $X''$, the intersection of 
pull-back of general $\hbox{dim}\f(E_2)$ very ample divisors on
$Z$, and the restriction of $\f$ to $X''$ (vertical slicing in [A-W]).
Then the formula proved on $X''$ remains valid also on $X$.

Thus $E_2$ is equal to one of the surfaces listed in (4.2.1).
Suppose first that $d:=L^2\cdot E_2>2$. Let $C$ be a general curve 
from the linear system $|L_{E_2}|$. Then $C$ is a rational curve
and since $H^0(E_2,L_{E_2})=d+2$ it follows that it moves inside $E_2$
in a family of dimension $d+1$. On the other hand, because of (1.4.1)
$\hbox{dim}_{[C]}Hilb(X)\geq d+n-3$. Since the degree of neighbouring fibers
of $\f$ is 2 at most and, being chosen generally, $C$ can not move to another component of the fiber of $\f$, it follows that
$\hbox{dim}_{[C]}Hilb(E_2)\geq \hbox{dim}_{[C]}Hilb(X)$ which implies
$n\leq 4$.

If $L^2\cdot E_2\leq 2$ then $E_2$ is either $\P^2$ or $\F_0$, or $\S_2$. Let
$C\subset E_2$ be a general conic.
If $\f$ is birational then all neighbouring non-trivial fibers of $\f$
have degree 1 with respect to $L$. 
Thus again, deformations of $C$ should remain 
inside $E_2$.
This, however, implies the inequalities:
$$\hbox{dim}_{[C]}Hilb(E_2)\geq \hbox{dim}_{[C]}Hilb(X)
\geq L\cdot C+n-3 \eqno(4.5.1)$$
where the right-hand-side inequality comes from (1.4.1). Therefore,
using the description of $E_2$ we have the following bound:
$$n\leq \hbox{dim}_{[C]}Hilb(E_2)+3-L\cdot C\leq 5+3-2=6.\eqno(4.5.2)$$
This proves the birational case.

In the fiber type case the degree of a general fiber of $f$ is 2, so
$C$ can move out of $E_2$ and therefore the
argument has to be adjusted accordingly. 
(We note that we may repeat the above argument for
curves of degree $>2$ but the result is not satisfactory.)
Let ${\cal H}$ be an irreducible component of the Hilbert scheme $Hilb(X)$
which contains a general fiber of the conic fibration $\f:X \ra Z$.
Over ${\cal H}$ we have the incidence variety ${\cal C}$.
That is, there exist a conic bundle $p:{\cal C}\ra
{\cal H}$ and a dominant morphism $q:{\cal C}\ra X$ which maps fibers
of $p$ to conics contracted by $\f$. We note that the map $q$ is birational with 
its exceptional locus in $q^{-1}(E_2)$ and
the composition $\f\circ q$ is proper. By possible normalization of ${\cal H}$
and base change of the conic bundle we may assume that ${\cal C}$ is normal.
Over $E_2$ the component ${\cal H}$ meets other components which are associated
to conics inside $E_2$. In particular, the components of $q^{-1}(E_2)$ are contained
is some incidence varieties of conics inside $E_2$.
Suppose that the dimension of a component ${\cal H}_{E_2}$ of $Hilb(E_2)$
is smaller than $n-1$, i.e. dim$_{[C]}Hilb(E_2)< n-1$ 
for some conic $C\subset E_2$
which is the case if $n\geq 7$.
Then, because of the deformation argument, ${\cal H}_{E_2}\subset
{\cal H}$ and $A:=p^{-1}({\cal H}_{E_2})$ is a component of $q^{-1}(E_2)$
 However, in view of [Ko, VI.1.5], $A$ is a divisor in ${\cal C}$
and therefore
$$n-1=\hbox{dim}A=\hbox{dim}{\cal H}_{E_2}+1\eqno (4.5.3)$$
which proves the fiber type case.

\medskip 

The argument which we presented above can be extended to deal with
possibly reducible fibers if we note the following two easy observations

\proclaim Lemma (4.6.1).
Let $C=\bigcup C_k$ be a (reduced) connected curve (with irreducible
components $C_k$) contained in a variety $F=\bigcup F_k$, where $F_k$ are
irreducible components of $F$. Suppose that the generic
point of any irreducible component $C_k$ 
is contained in $F_k-(\bigcup_{j\ne k} F_j)$.
Then a small deformation of $C$ in $F$, 
call it $C'$, has a decomposition into irreducible 
components $C'_k$ 
and the generic point of $C'_k$ is contained in $F_k-(\bigcup_{j\ne k} F_j)$.

\proclaim Lemma (4.6.2).
Let $F$ be a (possibly reducible) surface and $L$ an ample and spanned
line bundle on it. Suppose that a general curve $C\in |L|$ is a connected
curve of genus 0. Then $H^0(S,L)\leq L^2\cdot F+2$.

The subsequent result provides a description of isolated (c.f.~(4.0))
 2-dimensional fibers of contractions of varieties of dimension $\geq 5$.

\proclaim Proposition (4.7). Let $F$ be an isolated 2 dimensional fiber
of a good contraction, as in (4.0) and (4.2).
If $\f$ is birational then $n=\hbox{dim}X\leq 6$ and moreover
$F=\P^2$ if $n=5,\ 6$.
If $\f$ is of fiber type (and $L$ is $\f$-spanned) then $n\leq 7$ and
for $n\geq 5$ the pair $(F,L_F)$ is either $(\P^2,\O(1))$
or, for $n=5$, it is a smooth or reducible quadric
(i.e.~$(\F_0,C_0+f)$, or
 a union of two planes meeting along a line and $L$
restricted to each of planes is $\O(1)$).

\proof The preceeding lemma provides the upper bound on $n$. Also
the birational case is an immediate consequence of the preceeding argument.
Namely, if $L^2\cdot F\geq 2$ then we take a general curve $C\in |L_F|$
and consider its deformations to get $n\leq 4$ (we note that $C$ is connected,
possibly reducible and we can apply (1.4.1) directly to $C$ or
possibly to a connected smoothable sub-curve
of degree $\geq 2$). This concludes the birational case.
For the fiber type we need the following

\proclaim Lemma (4.7.1). In the situation of Proposition
(4.7) assume that $F$ is not irreducible. Let $F_1, F_2$
be two intersecting irreducible components of $F$. 
If they meet along  a curve $K$ then it is a line relatively to $L$. 
If they have an isolated common point then $n=4$, $\f$ is of fiber type
and $F_1= F_2=\P^2$, $L_{F_i}\iso\O(1)$.

\proof
If the curve $K$ is not a line then, taking
a general smooth section $X' \in |L|$, the map
$\f':= \f_{|X'} : X' \ra Z'$ would be a crepant contraction
with a fiber that contains two curves
meeting in more then one point. This is in contradiction with
the proposition (1.5.1).
Suppose now that $x$ is an isolated point of the intersection
$F_1\cap F_2$. Then obviously $n\geq 4$.
Let us take a general section of $L$ whose zero locus
contains $x$. Then we get two curves 
$C_i\subset F_i$, with $d_i:=L\cdot C_i=L^2\cdot F_i$,
meeting at the point $x$. We may assume that the curve $C= C_1\cup C_2$
is smoothable. This is clear if both $C_i$ are smooth (then they have to
meet tranversaly because $H^1(C_1\cup C_2,\O)=0$) otherwise one of $F_i$ is
a cone and then we can take a proper subcurve of $C_i$ so that
$C$ is connected and has degree $\geq 3$. 
We claim that $\hbox{dim}_{[C]}Hilb(F_1\cup F_2) = d_1+d_2$.
Indeed, any deformation of $C$ is obtained by deforming each of
$C_i$s with $x$ fixed and therefore dim$_{[C]}Hilb(F)=\hbox{dim}|L_{F_1}\otimes J_x|
+\hbox{dim}|L_{F_2}\otimes J_x|$. 
If $L\cdot C=d_1+d_2>2$ then argueing like in the first part
of the proof of (4.5) we get $n\leq 3$, a contradiction. 
If $L\cdot C=2$
then the argument which led to (4.5.3) provides us with the following
inequality:
$$n\leq \hbox{dim}_{[C]}Hilb(F_1\cup F_2)+2=4$$
which concludes the proof of (4.7.1).

\medskip
\noindent{\sl Conclusion of the proof of (4.7).} 
Now, in view of the above Lemma we can repeat the computations from the proof
of (4.5). Indeed, a general $C\in |L_F|$ is now a connected, 
possibly reducible, curve
of genus 0 and therefore the computations are in fact the same.
Thus, for $n\geq 5$ we get $L^2\cdot F\leq 2$. Suppose that $L^2\cdot F=2$. 
Then as in (4.5.3) we get 
$(n-1)\leq \hbox{dim}_{[C]}Hilb(F)+1= 4$ and thus $n\leq 5$. 
To conclude the proof of the proposition
we are only to exclude the case $F=\S_2$ for $n=5$. 
In this case we consider $C\subset \S_2$ which is a union of a conic and a line.
Then it is not hard to see that $\hbox{dim}_{[C]}Hilb(\S_2)=4=L\cdot C+1$ and thus
this case can occur for $n\leq 4$ only.

\medskip
\remark In (3.1.2) we provided examples of contractions of
manifolds of dimension $\geq 5$ which show that the results
of Proposition (4.7) are best possible; the only exception
is the case of $F=\P^2\cup\P^2$ for $n=5$.

\bigskip
In the remaining part of this section we will deal with the case $n=4$ and 3.
First we prove
\proclaim Lemma (4.8). If $\f: X\ra Z$ is a good contraction of a 4-fold
with an isolated 2-dimensional fiber $F$ 
and a general section $X'\in |L|$ is disconnected
then $F$ is a union of two copies of $\P^2$ meeting at one point
(we denote it by $\P^2\bullet \P^2$).

\proof We argue as in the proof of (4.7.1). Namely, we can find a decomposition
of $F=F_1\cup F_2$ so that $F_1$ and $F_2$ have an isolated 
common point $x$. Now we take a general $C\in |L_F|$ which contains $x$
and arrive to a contradiction if $C\cdot L\geq 3$.
\medskip
 With this argument we exhausted the technique of using deformations
of $C\in |L_F|$ to get informations about $F$. (We note that e.g.~the
inequality (4.4.2), whose different versions we considered above,
becomes useless if $C\in |L_F|$ and $n\leq 4$.) From now on we will
choose the curve $C$ for each case separately.
First we discuss when a surface $S$ among these listed in (4.2.1)
can be a component of the fiber $F$.
This will be done with the help of the Table II in which
we indicate the cases which we have
to discuss together with our choice of the curve $C\subset S$ and the numerical
data which comes up from the choice.

\midinsert{
\baselineskip 18 pt

{\bf Table II.}
{\sl \par\par
\settabs\+ Number\  &$(\F_r,C_0+kf$, $r\geq 3$, $k> r$ \ \ \ \
& curve curve curve $C$ \ \ \ &
\ \ $L\cdot C$\ \ \  &
dim$Hilb(S)$  &\cr

\+ N$^o$& pair $(S,L_S)$ \ & curve $C$ \ \ &\ \ $L\cdot C$
&dim$_{[C]} Hilb(S)$
\cr

\smallskip
\hrule

\medskip
\+(1)&$(\P^2,\O(2))$& line & \hfill$2$\hfill&\hfill2\hfill&
\cr
\+(2)&$(\S_3,\O(1))$& two lines &\hfill$2$\hfill&\hfill2\hfill&
\cr
\+(3a)&$(\F_r,C_0+kf)$, $r\geq 1$, $k> r$& $C_0+f$ & \hfill$k-r+1$\hfill &
\hfill$\leq 2$\hfill &
\cr
\+(3b)&$(\F_r,C_0+(r+1)f)$, $r\geq 3$& $C_0+2f$ & \hfill$3$\hfill
&\hfill2\hfill&
\cr
\+(4)&$(\F_2,C_0+3f)$& $C\in |C_0+2f|$ & \hfill$3$\hfill &\hfill3\hfill&
\cr
\+(5a)&$(\F_1,C_0+2f)$& $C\in |C_0+f|$ &\hfill$2$\hfill &\hfill2\hfill&
\cr
\+(5b)&$(\F_1,C_0+2f)$& $C\in |C_0+2f|$ &\hfill$3$\hfill & \hfill4\hfill&
\cr
\+(6a)&$(\F_0,C_0+2f)$& $C_0$ &\hfill $2$\hfill &\hfill1\hfill&
\cr
\+(6b)&$(\F_0,C_0+2f)$& $C\in |C_0+f|$ & \hfill$3$\hfill &\hfill3\hfill&
\cr
\+(7a)&$(\F_0,C_0+f)$& $C\in |C_0+f|$ &\hfill$2$\hfill & \hfill3\hfill&
\cr
\+(7b)&$(\F_0,C_0+f)$& $C\in |C_0+2f|$ &\hfill $3$\hfill &\hfill5\hfill&
\cr
\medskip
\hrule
}}\endinsert
\bigskip

\noindent (4.9).
Using Table II we can proceed with the argument that we explained in the course
of the proof of (4.5). We use the inequality (4.5.1)
which gives a fundamental constrain on a possible component $S$ of the fiber
$F$. That is, if $L\cdot C>1$ in the birational case and, respectively,
 $L\cdot C>2$ in the fiber type case, then
$\hbox{dim}_{[C]}Hilb(S)+3- L\cdot C\geq n$.
\par
Thus using the respective entries from Table II we get:
\item{(1)} $(\P^2,\O(2))$ may be a component of an isolated 2-dimensional
fiber of a birational (respectively fiber type) contraction only if
$n\leq 3$ (resp. $n\leq 4$). Moreover $(\P^2,\O(2))$ is not a component of a reducible fiber by (4.7.1) and (4.8).
\item{(2)} $(\S_3,\O(1))$ is not a component of an isolated 2-dimensional
fiber of birational type  contraction.
Moreover we note that the vertex of the cone $(\S_r,\O(1))$ 
has embedding dimension equal to $r+1$ so for $n\leq 4$ only $\S_3$ and $\S_2$
have to be discussed.
\item{(3)} if $(\F_r,C_0+kf)$ is a component of an isolated 2-dimensional fiber
and $r>1$ then $L_{|S}=C_0+(r+1)f$ and therefore $L\cdot C_0=1$.
Indeed, either $C_0$ is not contained in another component
of the fiber and we may use the estimate from the table or it is
 contained in another component of
the fiber and thus in view of () $L\cdot C_0=1$.
Moreover $\F_r$ can not be a component of an isolated 2-dimensional fiber
if $r\geq 3$. To see this note a small deformation $C'$ of
$C=C_0+2f$ has to remain inside
$S$ because it will contain a component of degree 2 through a generic
point of $S$. Since a deformation $C'$ of $C$ will be connected,
the components of $C'$  passing through             
a generic point of $S$ (union of two disjoint lines) will define uniquely
the remaining component of degree 1, which must be the curve $C_0$.
\item{(4)} $\F_2$ may be a component of an isolated 2-dimensional fiber
only if $n\leq 3$.
\item{(5)} $\F_1$ may be a component of an isolated 2-dimensional fiber
of a  good birational (resp. fiber type) contraction only if
$n\leq 3$ (resp. $n\leq 4$).
\item{(6)} $(\F_0,C_0+2f)$ can not be a component of an isolated 2-dimensional fiber of
a birational contraction and 
it may be a component of an isolated 2-dimensional fiber
of a fiber type contraction only if $n=3$.
\item{(7)} $(\F_0,C_0+f)$ may be a component of an isolated 2-dimensional fiber
of a birational (respectively fiber type) contraction only if
$n\leq 4$ (resp. $n\leq 5$)
\bigskip

Therefore we obtained a finite list of possible components of an isolated
2-dimensional fiber.
The results of the preceeding discussion are listed in Table III.
For clarity we also included the case when $n\geq 5$.
The last column of the table indicates the number of an example
which admits a fiber having $S$ as an irreducible component in the case $n
> 3$.
\midinsert
{\bf Table III}
{\sl
\smallskip
{\baselineskip 18 pt
\settabs\+ Number   &component
$(S,L_S)$\ \ \ & $\f$ birational \ \  & $\f$ of fiber type\ \ \ &Example
(for $n>3$)
\ \ \ \ \ \ \  \cr
\+ N$^o$ \ &component $(S,L_S)$\ & $\f$ birational \ \ & $\f$ of fiber type
& example (for $n>3$)
\cr
\medskip
\hrule
\medskip
\+ (0)&$(\P^2,\O(1))$&\hfill$n \leq 6$\hfill&\hfill $n\leq 7$\hfill&(3.1) and other
 \cr
\+ (1)&$(\P^2,\O(2))$&\hfill$n \leq 3$\hfill&\hfill $n\leq 4$\hfill& (3.4.0) \cr
\+ (2)&$(\S_3,\O(1))$&\hfill--------\hfill& \hfill$n= 4$\hfill& (3.3)\cr
\+ (3)&$(\S_2,\O(1))$&\hfill$n \leq 4$\hfill&\hfill $n\leq 4$\hfill&(3.2.1)
and (3.5.5)\cr
\+ (4)&$(\F_2,C_0+3f)$&\hfill$n = 3$\hfill& \hfill$n= 3$\hfill& \cr
\+ (5)&$(\F_1,C_0+2f)$&\hfill$n = 3$\hfill& \hfill$n\leq 4$\hfill&
(3.4.1), (3.4.3)\cr
\+ (6)&$(\F_0,C_0+2f)$&\hfill--------\hfill&\hfill $n=3$\hfill& \cr
\+ (7)&$(\F_0,C_0+f)$&\hfill$n \leq 4$\hfill& \hfill$n\leq 5$\hfill& (3.1)
and other\cr
\medskip
\hrule\bigskip
}}
\endinsert

\noindent (4.10).
Now we discuss the question of reducible fibers.

\proclaim Lemma (4.10.1).
The cone $\S_3$ can not be a component of a reducible 2-dimensional fiber
if $n=4$. Also $\S_2$ can not be a component of a reducible fiber if $n=3$
and in dimension 4 it meets one component at most. No three, respectively,
four components of $F$ can meet in a single point if $n=3$ or 4, respectively.

\proof 
We can take a general $X'\in |L|$ which passes through the point $x$
which denotes the vertex of the cone or the meeting point of
components. Such $X'$ is then smooth because there exist lines passing
through $x$ which are not contained in $X'$. Therefore, in view of
(1.5.2) at most 2, or respectively, 3 components of $X'\cup F$ should
meet in $x$ if $n=3$ or $4$, respectively.  So the last statement of
the lemma is clear.  Since the cones meet other components along
lines in the ruling, which contain $x$, this implies the other claims
too.

\medskip

\proclaim Lemma (4.10.2). Let $F$ be an isolated 2 dimensional fiber 
of a good contraction $\f$ of a 4 fold. Let $S_1$ and $S_2$ be two
(different) irreducible components of $F$ which meet along 
a line $l$. Then at least one of them is $\P^2$ and if $\f$
is birational then both of them are $\P^2$. Moreover, if one of the
components is $\F_1$ then the meeting line is $C_0$.

\proof
In the union $S=S_1\cup S_2$ we construct a connected curve $C=C_1\cup C_2$,
with irreducible components $C_i\subset S_i$ passing through the generic point
of each $S_i$ and having a common point at the line $l=S_1\cap S_2$. 
We assume that $S_1\ne \P^2$. In Table IV we indicate the choice of $C_1
\subset S_1$ which will depend on the position of $l$ in $S_1$. 
The curve $C_2\subset S_2$ will be always a general member
of $|L_{S_2}|$ passing through $C_1\cap l$. We note that
if $d_2=L\cdot C_2= L^2\cdot S_2$ then dim$|L_{S_2}|=d_2+1$.
Then we compute the degree of 
$C$ and the space of its deformations inside $S$ and the lemma follows by
comparing these two numbers. Namely, it follows that the curve $C$ moves
out of the fiber $F$, so that its degree has to be bounded by 1 or
2, if the contraction is birational or of fiber type, respectively.

\midinsert

{\baselineskip 20 pt
{\bf Table IV.}
{\sl 
\settabs\+ number& ${\F_2}_{C_0}\cup_{C_0}\F_2$\ \ \ \ \ \ \ \ \ \ \ \
&curve curve $C_1\subset S_1$ \ \ &
$L\cdot C=r+1+d_2$ \ \ \ &dim$_C Hilb(S)$\ \ &
\cr
\+ N$^o$& $S_1$ and $l$&  curve $C_1\subset S_1$ \ \ &$L\cdot C$
&dim$_{[C]} Hilb(S)$
\cr
\medskip
\hrule
\medskip
\+(A)&$\S_2$,  $l=$\ line & line &$1 + d_2$&$1+d_2$
\cr
\+(B)&$\F_r$, $r\geq 0$, $l=C_0$ & $f$ & $1+d_2$& $1+d_2$ 
\cr
\+(C)&$\F_r$, $r\geq 0$, $l=f$ & $C_1\in |rf+C_0|$&  $r+1+d_2$& $r+1+d_2$ 
\cr
\medskip\hrule\medskip
}}
\endinsert

\medskip
A similar argument is made in the proof of
\proclaim Lemma (4.10.3). In the above situation no three components can 
meet along a line.

\proof
Assume the contrary, i.e.~let $S_1$, $S_2$ and $S_3$ meet along
a line $l$.
From what we have just proved we know that two of them are $\P^2$.
Then, argueing exactly like in the previous lemma we conclude that
the third must be also $\P^2$. Namely, we consider $S_1=\F_r$ and $C_1$
as in Table IV, while $C_2\subset \P^2\cup\P^2$ consists of two lines,
one in each component.
Now we claim that the conormal
of $l$ in $X$ is (at least generically) spanned by the conormals of $l$
in each of $S_i$. Indeed, if we consider a general linear section $X'\in |L|$
then the tangent vectors of lines $S_i\cap X'$ will spann
the tangent space at $X'\cap l$. This however implies that we have an
embedding 
$$0\raa N_{l/S_1}\oplus N_{l/S_2}\oplus N_{l/S_3}=\O(1)^{\oplus 3}\raa
N_{l/X}\raa T$$
where $T$ is a torsion sheaf on $l$. This is absurd, since,
because of adjunction, $c_1(N_{l/X})=-1$.

\medskip
Now we are ready to prove the classification of isolated 2 dimensional fibers 
in dimension 4:
\proclaim Proposition (4.11).
Let $\f: X\ra Z$ be a good contraction
of a smooth 4-fold
and $F$ is an isolated 2-dimensional fiber.
If $\f$ is birational then the pair
$(F,L_F)$ is isomorphic to one of the following:
$$\matrix{(\P^2,\O(1)),& (\F_0,C_0+f),& (\S_2,\O(1)),& (\P^2\cup \P^2,
\O(1))}.$$
If $\f$ is of fiber type and $F$ is irreducible then $(F,L_{F})$ is
one of the following:
$$\matrix{(\P^2,\O(1)),&(\P^2,\O(2)),&(\S_3,\O(1)),&(\S_2,\O(1)),&
(\F_1,C_0+2f),&(\F_0,C_0+f).}$$ If $\f$ is of fiber type and $F$ is
reducible then it has at most three components and it is one of the
following (we supress the description of $L$ which is obvious):
$$\matrix{\P^2\cup\P^2,& \P^2\cup\F_0,& \P^2\cup_{C_0}\F_1,&
\P^2\cup\S_2,& \P^2\cup\P^2\cup\P^2,&\P^2\cup_f{(\F_0)}_{C_0}\cup\P^2,}$$
where any two components intersect along a line (indicated by a
subscript, when needed), and the exceptional case of 
$\P^2\bullet\P^2$ when the two components intersect at an isolated point.

\proof The birational case is an immediate consequence of Table III
as well as (4.10.2), (4.10.3) and
(4.7.1). Indeed, if a fiber is reducible then all its components are $\P^2$.
If it had more than two components then, since three of them do not 
meet along a common 
line, two of the components would meet along a common isolated point, 
contradictory to (4.7.1).

Now we pass to the fiber type case. The description of irreducible fibers and 
fibers of two components is already known due to Table III and (4.10.2).
Thus let us pass to fiber which contain at least three components, 
call them $S_i$. From (4.10.2) we know that from a pair of 
two components meeting along
a line, at least one is $\P^2$.
Moreover, no three components make
a cycle. That is, it is not possible that
 $S_1$ meets $S_2$ along $l_{1,2}$, $S_2$ meets $S_3$ along $l_{2,3}$
and $S_3$ meets $S_1$ along $l_{3,1}$ and all three lines are different,
because then the linear section of $S_1\cup S_2\cup S_3$ would
contain a cycle of rational curves. 

We recall that according to (4.7.1) if two components meet at
a single point then both are $\P^2$.  Thus, 
 if $S_2$ meets $S_1$ and $S_3$ along two different lines which intersect
then $S_1=S_3=\P^2$. According to Table III and lemmata (4.10.1) and (4.10.2) 
the central component $S_2$ is either $\P^2$ or $\F_0$. 
If the central component $S_2$ meets the other two components 
(equal to $\P^2$) along two 
non-intersecting lines then it has to be $\F_0$. In this case, however,
we can take a conneted curve $C$ of degree 3 consisting of three lines,
each in one component, and we get contradiction to the deformation
principle because dim$_{[C]}Hilb(F)=3$.

The case of more than 3 components is ruled out similarly. First we note that
all component would be $\P^2$ and then we can apply argument similar to
the one in (4.7.1) to get contradiction.

\bigskip

\remark. In Section 3 we have examples of appropriate 2 dimensional fibers
except cases $\P^2\cup\S_2$,
$\P^2\cup\P^2\cup\P^2$ and $\P^2\cup_f{(\F_0)}_{C_0}\cup\P^2$. A list of possible exceptional
fibers of a fiber type contraction of a 4-fold was obtained independently
by T.~Kachi [Kac], whom
we owe our thanks for pointing out missing cases in our preliminary
list.

\bigskip \noindent
(4.12) The last part of this section is dedicated to the case $n = 3$;
this case is  known (and proved by S.~Mori in the case of elementary
contractions) but the subsequent discussion may be
 a good test for our approach.
Our goal is to describe all possible two dimensional isolated fibers of a good
contraction of threefolds. A particular feature of this dimension is the
fact that these fibers are actually divisors.

The list of possible components of a fiber is set up in Table III.
Also,  because of (4.7.1) and (4.10.1) neither $(\P^2,\O(2))$ nor $\S_2$ 
is  a component of a reducible fiber
and --- the fact which will be used constantly in our arguments ---
no three components meet at one point.
\par
Let us discuss first the question of reducible fibers: assume therefore that
$F$ contains at least two intersecting components, $F_1$ and
$F_2$. Then $F_1 \cap F_2 := R$ is a line relative to $L$ and the 
two components intersect in $R$ tranversally (see Lemma (4.4) and (1.5))
\par
We use now the intersection theory of divisors to prove a
useful formula. Let $N_{R/X}$ denotes the normal bundle of $R$ in $X$;
we have  $N_{R/X} = N_1\oplus N_2$, where
$N_i$ denotes the normal bundle of $R$ in $F_i$, and
$deg(det(N_{R/X})) = -1$, by adjunction formula.
It follows that
$$-1 = deg (det(N_{R/X})) = deg N_1 + deg N_2 = F_1\cdot R + F_2\cdot R.$$
Since $F_i\cdot R = (K_{F_i}-K_X)\cdot R = K_{F_i}\cdot R +1$,
we obtain
$K_{F_1}\cdot R + K_{F_2}\cdot R = -3$.

\par
The curve $R$ is a line relatively to $L$; therefore, modulo numerical
equivalence, it is
unique in $\P^2$ while it can be either $f$ or $C_0$ in $\F_r$,
the last case only if $L_{\F_r} = C_0 + (r+1)f$.
This implies that $K_{F_i}\cdot R = -3$, if $F_i \iso \P^2$, and
$K_{F_i}\cdot R = (r-2)$ or $-2$, if $F_i \iso \F_r$.
This observation added to the above formula gives only two possibilities
(up to possible renumeration of components),
namely the following:
\par
\item{(i)} $F_1 = \P^2$ and $F_2 = \F_2$,  $R$ is a line in $\P^2$ and the
section $C_0$ in $\F_2$;
\par
\item{(ii)} $F_1 = \F_1$ and $F_2 =\F_r$,
$R$ is the section $C_0$ in $\F_1$ and a fiber $f$ in $\F_r$.

If we are in the case (ii) we can apply the argument of deforming
of a rational curve of degree 3 which is the union of
$f \subset \F_1$ and of $C_0 +f \subset \F_r$ (none of these curves is 
contained in another component because our ``no three meet'' rule).
Since $\hbox{dim}_{[C]}Hilb(F) \leq 2$ for $r>0$,
we obtain a contradiction unless $r=0$.
Then $F = \F_0 \cup \F_1$.  If the fiber $F$
has at least three components, by the above,
we must have that two of them are $\F_1$ and one is $\F_0$,
the $\F_1$ are disjoint and they intersect with the $\F_0$
along two disjoint fibers of its ruling. Deforming
a degree 3 rational curve in $F$ which is the union of
the section $C_0$ in $\F_0$ and of two fibers, one in each $\F_1$,
we obtain a contradiction, as above.

Now we know all possible fibers and it remains to distinguish the
fibers of birational contractions from these of fiber type
contractions.  To this end we note that if the fiber in question has
ample conormal bundle (i.e.~$\O_F(-F)$ is ample) then it is an
isolated positive dimensional fiber of the contraction (see also
(6.1)). Indeed, the
limit of 1-dimensional fibers approaching the fiber $F$ would produce
a curve in $F$ whose intersection with $F$ would be zero. Since
$\O_F(-F)$ can be computed easily by adjunction, we find out all the
pairs $(F,L_F)$ which may be non-isolated positive dimensional
fibers; these are the following:
$$\matrix{(\F_0,C_0+2f),&(\F_1,C_0+2f),&(\F_0\cup_{C_0}\F_1,\  
L_{\F_0}=C_0+f,\ L_{\F_1}=C_0+2f).}$$
Moreover, only in the case $(\F_1,C_0,+2f)$ the normal $F_F=C_0+f$
has zero intersection with the unique curve $C_0$ while in the remaining
two cases the locus of curves which have intersection 0 with $F$ coincides
with $F$. This distinguishes the birational and fiber case. Thus we have 
proved

\proclaim Proposition (4.13). Let $\f$ be a good contraction of a smooth
$3$-fold with an isolated two dimensional fiber $F$; as usually $L = -K_X$.
 If $\f$ is birational then $(F, L_F)$ is one of the following pairs
$$\matrix{(\P^2, \O(1)),&(\P^2,\O(2)),&(\S_2,\O_{\S_2}(1)),
&(\F_0 , C_0 + f),&(\F_1, C_0 + 2f)},$$
or  $\P^2 \cup_{C_0} \F_2$ and $L_{|\P^2} = \O(1)$,
$L_{|F_2} = C_0 + 3f$.
In case when $F=\F_1$ the contraction $\f$ contracts also a
smooth divisor $E \subset X$ to a smooth curve 
$E \cap F_1 = C_0 \subset F_1$. In the remaining birational cases
the fiber $F$ is an isolated poistive dimensional fiber of $\f$.
If $\f$ is of fiber type with generic fiber
of dimension 1 (a conic fibration) and $L$ is $\f$-spanned then the pair
$(F,L_F)$ is either
 $(\F_0, C_0 + 2f)$ or
 $F=\F_0 \cup_{C_0}\F_1$ and $L_{|\F_0} = C_0+f$, 
$L_{|\F_1} = C_0 + 2f$.

Let us note that all the above cases exist. Indeed, the construction
of elementary contractions with appropriate fibers is done by Mori
in [Mo1]. All cases which are not in the Mori's
list  are obtained from Mori's contractions.
Namely, the cases $(\F_1,C_0+2f)$ and $\P^2\cup\F_2$
are obtained from $(\P^2,\O(1))$ by a blow-up of a curve.
The case $(\F_0,C_0+2f)$ is just a product of $\P^1$ and a simple
blow-up of a surface while
the case of $\F_0\cup\F_1$ is obtained by blowing a line $C_0$
in the previous one.
Alternatively, the last case is the blow-up of the quadric cone singularity 
of the scroll 
$$\P^2\times\C^2\supset\{([t_0,t_1,t_2],(z_1,z_2): t_1z_1+t_2z_2=0\}
\raa \C^2$$
c.f.~Example (3.4.3) and (3.2) in [B-W].

\beginsection 5. Scheme theoretic fiber.

(5.0) Let $\f : X \ra Z$ be a good contraction of a smooth $n$-fold $X$.
As in (4.0) we fix a relative
ample line bundle $L:= -K_X$ and we choose $F=\f^{-1}(z)$
to be a positive dimensional fiber of $\f$.
\par
Suppose that we already know the geometric (reduced) structure of the
fiber $F$ of the map $\f$.
Then the next step in understanding of the map $\f$ is
the conormal sheaf (bundle) of $F$ in $X$, which is defined as the quotient
$N^*_{F/X} =: \I_F / \I_F^2$, where $\I_F$ is the ideal of $F$ with the
reduced structure.
\par
If $F$ is a locally complete intersection then $N^*_{F/X}$ is locally
free over $F$ (a vector bundle) and its dual is the normal bundle
$N_{F/X}$.  We note that almost all among possible fibers which we
have described in Section 4 are locally complete intersections
(exceptions are $\S_3$ and most of the reducible fibers) and moreover
the blow-up of $X$ along $F$ has always terminal singularities.
\par
If $F$ is locally complete intersection then we have  the adjunction formula:
$K_F=(K_X)_{|F} + det N_{F/X}$.
In particular, if $F$ and $L_{|F} =- K_{X{|F}}$
are fixed then the first Chern class of the normal bundle is given.
If $F$ has codimension $1$ the normal bundle is thus fixed.
\par

\par
It turns out that in many instances
the knowledge of the fiber and its conormal allows to describe the contraction
map. A typical example of such reasoning is the Castelnouvo theorem which
we recalled in (2.4).
The formal function theorem which was used in the proof of the Castelnouvo
theorem
allowed us to compare "asymptotic" behaviour of the ideals $\I_F$ and $m_z$.
This leads to a natural question of comparison of ideals $\I_F$ and
$\f^{-1}m_z\cdot\O_X$,
that is of understanding the fiber structure on $F$.

Let us begin with the following general observation.

\proclaim Lemma (5.1).
Let $\f: X\ra Z$ be a projective morphism of normal varieties such
that $\f_*\O_X=\O_Z$.
For a point $z\in Z$ let $F=\f^{-1}(z)$ be the geometric
fiber of $\f$ with the ideal (sheaf of ideals) $\I_F$.
Then $\f_*\I_F\subset\O_Z$ is the maximal ideal $m_z$ of the point $z$ and
the scheme theoretic structure on $F$ is defined by the ideal
$\I_{\tilde F}$ which is the image of the evaluation
$\f^{-1}\f_*(\I_F)\ra\I_F$.

Therefore, in order to understand the structure of the map $\f$ we will
analyse the behaviour of pull-backs and push-forwards.

A very useful example is the contraction
to the vertex of a cone.
\par

\example (5.1.1)
Let $Y$ be a smooth variety and let ${\cal L}$ be
an ample line bundle over $Y$.
Let $X:=Spec_Y(\bigoplus_{k\geq 0} k\L)$ be the total space of the dual
bundle $\cal L^*$
with a zero section $Y_0$. Consider the collapsing $\f:X\ra Z$ of $Y_0$
to the vertex $z$ of a cone $Z$. That is,
$Z=Spec(\bigoplus_{k\geq 0} H^0(Y,k\L))$
and the map $\f$ is associated to the evaluation of $m\L$.
If $Y$ is Fano and $-K_Y-{\cal L}$ is ample then the contraction
$\f$ is good. The maximal ideal of the vertex is
$m_z=\bigoplus_{m>0}H^0(Y,k\L)$. Let us note that $m_z=\f_*\O_X(-Y_0)$.

The fiber structure of the fiber $Y_0$ is defined by sections of
bundles $k\L$. More precisely, at a point $y\in Y_0$ the fiber structure ideal
is generated by functions $s_0^k\cdot s_k$ where $s_0$ is a local
generator of the reduced ideal of $Y_0$ in $X$ (zero section of $\L^*$)
and $s_k\in H^0(Y,k\L)$. Therefore the scheme theoretic structure 
of the fiber coincides with
the geometric structure at $y$
if and only if the exists $s_1\in H^0(Y,\L)$ which does not vanish
at $y$ (or $\L$ is generated at $y$).
In particular, the multiplicity of the fiber is 1,
or the fiber structure coincides with the geometric structure
at a general point if and only if the bundle $\L$ has a non-zero section.
In such a case the fiber structure has embeded components
at base-points of the line bundle $\cal L$.
In this situation, the natural property
$H^0(\tilde F,\O_{\tilde F})=\C$ is not true.

Let us note that the gradation in the ring
$\bigoplus_{k\geq 0}H^0(Y,k\L)$ may not coincide with the gradation of
the maximal ideal $m_z$ of the vertex. In fact these two coincide if and only
if  $\L$ is projectively normal i.e. the map $S^k(H^0(Y,\L))\ra H^0(Y,k\L)$
is surjective for all $k\geq 0$.

Then $Z$ is an affine cone over the $Y$
embedded in a projective space by $|\L|$.
In such a case $X$ is the blow-up of $Z$ at $z$ and
$m/m^2\iso H^0(Y,\L)$. In fact, we have the following criterion
relating the contraction to the vertex and the blow-up of the vertex.

\proclaim Lemma (5.2). (c.f.~[EGA], (8.8.3)) Suppose that $H^0(Y,\L)\ne 0$. 
Then $X$ is the blow-up of $Z$ at $z$ i.e.~$X=Proj_Z(\bigoplus_k m_z^k)$ 
if and only if $\L$ is spanned over $Y$.

\proof
If $X\ra Z$ is the blow-up of $m_z$ then $\f^{-1}m_z$ is invertible on $X$
and defines the fiber structure over $Y_0$. Therefore $\f^{-1}m_z=\O_X(-kY_0)$
for some $k\geq 1$.
Since $H^0(Y,\L)\ne 0$, the ideal is generically reduced and hence
$\f^{-1}m_z=\O_X(-Y_0)$. Thus $\L=\O_{Y_0}(-Y_0)$ is spanned.
Conversely, suppose that $\L$ is spanned. Then, according to what we have said
above, $\f^{-1}m_z=\O(-Y_0)$ and thus by the universal property of the blow-up
(see [EGA] or [Ha] (7.14)) we have a map $X\ra Proj_Z(\bigoplus_k m_z^k)$
which is an isomorphism.

\medskip

A similar situation occurs for semiample vector bundles.

\example (5.1.2) Let $\E$ be a rank $r$
semiample bundle on a smooth variety $Y$. 
That is $\O_{\P(\E)}(1)$ is semiample, or equivalently,
the symmetric power $S^k(\E)$ is generated
by global sections for $m\gg 0$
(see for instance Example (3.1)). Then similarly as above we consider
$X=Spec_Y(\bigoplus_{k\geq 0} S^k\E)$, the total space of the dual
bundle $\E^*$
with a zero section $Y_0$, and the collapsing 
$\f:X\ra Z=Spec(\bigoplus_{k\geq 0} H^0(Y,S^k\E))$ of $Y_0$
to the special point $z\in Z$.
If $Y$ is Fano and $-K_Y-det\E$ is ample then the contraction
$\f$ is good. Moreover $\f$ is birational if and only if
its top Segre class is positive. If $\E$ is spanned then the 
scheme theoretic fiber $\tilde Y_0$ is reduced and $\f$ factors through
the blow up of $Z$ at $z$.

\medskip

The good properties of the above examples can be extended to arbitrary
good contractions; the assumption which is needed is nefness of
the conormal of the fiber. Then the map behaves similarly as the
contraction to the cone because the sections of the conormal of the fiber
extend.

\medskip\noindent{\bf Lemma (5.3).}
{\sl  Let $\f: X\ra Z$ be a good or crepant contraction of a
smooth variety
with a fiber $F=\f^{-1}(z)$. Assume that $F$ is locally complete intersection
and the blow up $\beta:\hat X\ra X$ of $X$
along $F$ has log terminal singularities. By $\hat F$ we denote the
exceptional divisor
of the blow-up. Let $\L$ be an a line bundle on $X$ such that $-K_X+\L$ is
$\f$-big and nef.
If the conormal bundle $N^*_{F/X}$ is nef
then:
\par
\item{(a)} The line bundle $\O_{\hat X/X}(1)=-\hat F$ is $\f\circ\beta$-nef.
\item{(b)} Any section of $N^*_F$ extends to a function in $\Gamma(X,\O_X)$
vanishing along $F$, that is the natural map
$$\matrix{\f^!:m_z&\raa& H^0(F,N^*_F)=H^0(\hat F,\O_{\hat F}(-\hat F))\cr
m_z\ni f&\mapsto&(x\mapsto [\f\circ f]\in (\I_F/\I_F^2)_x)}$$
 is surjective.
\item{(c)} For $i>0$ and $t\geq 0$ we have vanishing
$H^i(F, N^*_{F/X}\otimes\L) = 0 $.
\item{(d)} Some positive multiple
$\O_{\hat X/X}(k)=-k\hat F$ is $\f\circ\beta$-spanned
and it defines a good contraction over $Z$:
$$\hat\f:\hat X\raa\hat Z=
Proj_Z(\bigoplus_k (\f\circ\beta)_*\O_{\hat X}(-k\hat F);$$
the scheme $\hat Z$ is a blow-up of $Z$ along some ideal of a scheme
supported at $z$.

}\medskip

\proof The nefness of $-\hat F$ is clear.
To prove (b) we consider a sequence
$$0\raa\O_{\hat X}(-2\hat F)\raa\O_{\hat X}(-\hat F)\raa\O_{\hat F}(-\hat
F)\raa 0.$$
Since, by assumption,
$-2\hat F-K_{\hat X}=-(dimX-dimF+1)\hat F-\beta^*K_X$ is
$\f \circ \beta$-big and nef,
and moreover $\hat X$ has good singularities it follows that
$H^1(\hat X,\O_{\hat X}(-2\hat F))=0$ and sections of $\O_{\hat F}(-\hat F)$
extends to $\hat X$. This implies (b).

To prove (c) we note that $H^i(F, N^*_{F/X}\otimes\L) =
H^i(\hat F,(-\hat F+ \beta ^*\L)_{|\hat F})$.
Moreover, as above, we note that
the line bundle $-s\hat F +  \beta^*\L -K_{\hat X}$
is $\f\circ\beta$-nef and big for $s\geq 0$.
Now we can apply the Kawamata-Viehweg vanishing to the cohomology of the
sequence
$$0 \raa (-2\hat F +\beta^*\L) \raa (-\hat F + \beta^*\L)
\raa (-\hat F + \beta^*\L)_{|\hat F} \raa 0$$
to get (c).

The claim (d) follows from (a) by Kawamata-Shokurov base-point-free
theorem. Indeed,
since $K_{\hat X}=(dimX-dimF-1)\hat F+\beta^*K_X$
then $-\hat F-K_{\hat X}=-(dimX-dimF)\hat F-\beta^*K_X$ is
$\f\circ\beta$-nef and big and some multiplicity of $-\hat F$ is spanned
by global sections on $\hat X$.
The statement on the blow-up of the ideal is general --- see [Ha], (7.14).
\medskip

\proclaim Corollary (5.3.1). In the situation of the previous Lemma
the exceptional set $\hat G$ of the blow-up $\hat Z\ra Z$ is
equal to $Proj(\bigoplus_k H^0(F,S^k(N^*_F)))$ and the map
$\hat\f_{\hat F}:\hat F\ra\hat G$ is defined by the
evaluation
$\bigoplus_k H^0(F,S^k(N^*_F))\ra \bigoplus_k S^k(N^*_F)$.
If $\hat G$ is irreducible and
$H^0(F,N^*_F)\ne 0$ then $\hat G$ is a Cartier divisor in $\hat Z$ and
$\O_{\hat G}(\hat G)\iso \O(-1)$, where the latter bundle is defined naturally
on the $Proj$ of the graded ring. If moreover $\hat G$ is smooth then
$\hat Z$ is smooth along $\hat G$.

\proof The connected part of the Stein factorization of the
 map $\hat\f_{|\hat F}$ is the evaluation map
$\hat F\ra Proj(\bigoplus_k H^0(F,S^k(N^*_F)))$ defined above. Then the map
$Proj(\bigoplus_k H^0(F,S^k(N^*_F)))$ to $\hat Z$ is associated to
the restriction
$$(\f\circ\beta)_*\O_{\hat X}(-k\hat F)\raa H^*(\hat F,\O_{\hat F}(-\hat F))
\iso H^0(F,N^*_F)$$
which is surjective because of Kawamata-Viehweg vanishing theorem.

Therefore the scheme
$Proj(\bigoplus_k H^0(F,S^k(N^*_F)))$ embeds in $\hat Z$ as
the exceptional set $\hat G$ of $\alpha:\hat Z\ra Z$.

By our construction $\hat\f^*(\O_{\hat Z}(1)=\O_{\hat X}(-\hat F)$.
If $s\in H^0(F,N^*_F)$ is a non-zero section then it extends to a section
of $\O_{\hat X}(-\hat F)$ over $\hat X$
and it descends to a divisor $D\in\O_{\hat Z}(1)$
which does not contain $\hat G$. Now, using the embedding
$\O_{\hat X}(-\hat F)\subset \O_{\hat X}\iso \O_Z$ we find a global function
$f$ on $Z$ such that $\alpha^*(f)=D+a\hat G$. Since the multiplicity of
$\beta^*(f)$ along (at least one of the components of) $\hat F$ is 1
it follows that $f$ vanishes with multiplicity 1
along $\hat G$ and thus $\hat G$
is Cartier and $\O_{\hat G}(-\hat G)\iso\O_{\hat G}(1)$.

Finally to get the last claim of the Corollary we note that
if a Cartier divisor is smooth then the ambient variety is smooth along
the divisor in question --- this is a general fact.

\medskip\noindent {\bf Proposition (5.4). }
{\sl Let $\f: X\ra Z$ be a good  or crepant contraction
of a smooth variety.
Assume that $F=\f^{-1}(z)$, a geometric fiber of $\f$, is locally complete
intersection
with the conormal bundle $N^*_{F/X}=\I_F/\I^2_F$.
Suppose moreover that the blow up $\beta:\hat X\ra X$ of $X$
along $F$ has log terminal singularities.
Then the following conditions are equivalent:
\par
\item{(a)} the bundle $N^*_{F/X}$ is generated by global sections on $F$,
\item{(b)} the invertible sheaf
$\O_{\hat X}(-\hat F)$ is generated by global sections
at any point of $\hat F$.
\item{(c)} $\f^{-1}m_z\cdot\O_X=\I_F$ or, equivalently,
the scheme-theoretic fiber structure of $F$ is reduced and contains
no embedded components, i.e. $\tilde F=F$.
\item{(d)} there exists a good contraction
$\hat\f : \hat X \ra \hat Z=Proj_Z(\bigoplus_k m_z^k)$
onto a blow-up of $Z$ at the maximal ideal of $z$,
and $\f^*(\O_{\hat Z}(1)) =
\O_{\hat X}(1)$.
\item{}
}

\proof (a) implies (b) because of the previous Lemma, part (b).
Claims (b) and (c) are equivalent because $\beta_*\O_{\hat X}(-\hat F)=
\I_F$ and $\beta^{-1}(\I_F)=\O_{\hat X}(-\hat F)$.
The implication (b)$\Rightarrow$(d) follows by the universal property of
the blow-up, since by (b) $(\f\circ\beta)^{-1}m_z\cdot\O_X=\I_{\hat F}$.
The implication (d)$\Rightarrow$(a) is clear since $\O_{\hat Z/Z}(1)$
is spanned over $Z$.
\bigskip

Before stating the last result of this subsection let us recall
that for a local ring $\O_{Z,z}$ with the maximal ideal $m_z$ one defines
the graded $\C$-algebra $gr(\O_{Z,z}):=\bigoplus_k m_z^k/m_z^{k+1}$.
The knowledge of the ring $gr(\O_{Z,z})$ allows sometimes to describe the
completion ring $\hat\O_{Z,z}$, like in the Castelnuovo theorem (2.4).
Also, we will say that a spanned vector bundle $\E$ on a projective variety
$Y$ is p.n.-spanned (p.n. stands for projectively normal) if
for any $k> 0$ the natural morphims $S^kH^0(Y,\E)\ra H^0(Y,S^k\E)$
is surjective. As we noted while discussing the contraction to the vertex,
projective normality allows us to compare gradings of rings ``upstairs''
and ``downstairs''.

\proclaim Proposition (5.5). (c.f. [Mo1], 3.32])
Let $\f:X\ra Z$ be a contraction as in the previous Proposition.
Suppose moreover that $N^*_F$ is p.n.-spanned. Then
$\f_*(\I^k_F)=m^k_z$,\ \  $\f^{-1}(m^k_z)\cdot\O_X=\I^k_F$
and there is a natural
isomorphism of graded $\C$-algebras:
$$gr(\O_{Z,z})\iso\bigoplus_kH^0(F,S^k(N^*_F)).$$

\proof See [Mo1], p.164.

\bigskip

\noindent {\bf (5.6). The normal bundle of a 1-dimensional fiber.}
\par\noindent
We first recall the case in which $F$ is a fiber of dimension 1. This is
well known and we will give it as a warm up
before discussing the two dimensional case.
Let $C$ be an irreducible component of $F$.
As we have seen in (4.1) $C$ is a rational curve and it can be
either a line or a conic with the respect to
$L$; in the last case $\f$ is of fiber type.
\par
Let $\I$ be the ideal
of $C\subset X$ (with the reduced structure) and consider the exact sequence
$$0 \raa \I/\I^2 \raa \O_X /\I^2
 \raa \O_X /\I \raa 0.$$
In the long cohomology sequence associated
the map of global sections $H^0(\O_X /\I^2)\ra H^0({\O_X /\I})$
is surjective; moreover, by (1.2.1), $H^1({\O_X /\I^2})=0$.
Therefore $H^1({\I/\I^2}) = 0$ which gives a bound
on $N^*_{C/X}=\I/\I^2$. Namely if $N_{C/X} = \oplus \O(a_i)$ then
$a_i < 2$.  On the other hand, by adjunction,
$det(N_{C/X})=  \Sigma a_i =\O(-2-K_X.C)$ and thus the list of possible values
of $(a_1,...,a_{n-1})$ are finite.
\par
If $\f$ is a good birational
contraction then we have even a better bound because,
similarly as above and using (1.2.1), we actually
get $H^1(N^*_{C/X}\otimes\O(K_X.C))=0.$
Therefore, since $K_X.C =1$,
there is only one possibility, namely $N_{C/X}=\O(-1) \oplus\O^{(n-2)} $.
\par
If $\f$ is of fibre type then the estimate coming
from this technique is not sufficient and one has to use other
arguments. More precisely, one has to use a scheme associated
to a double structure on $C$ --- see [An] --- or one may use arguments
coming from the deformation theory as it follows.
Namely the possibilities which can occur from the above vanishing,
if $n=3$, are the following:
$$\matrix{\O\oplus\O, & \O\oplus\O(-1), & \O(1)\oplus\O(-2), &
\O(1)\oplus\O(-1).}$$
We will show that the last possibility does not occur. The argument which
we apply is related to the deformation technique (1.4.1) and it will
 be used later to deal with 2 dimensional fibers too.

\proclaim Lemma (5.6.1). The normal bundle $N_{C/X}$ cannot be
$\O(1)\oplus\O(-1)$.

\proof Assume the contrary and let $\psi : \hat X \ra X$ be the blow-up of
$X$ along $C$;
let $E:=\P(\O(1)\oplus\O(-1))$ be the exceptional divisor.
Let $C_0$ be the curve contained in $E$ which is the section of the ruled
surface
$E \ra C$ corresponding to the line bundle $\O(-1)$. We have immediately that
$E\cdot C_0 = 1$ and that $\psi_{C_0}$ is a $1-1$ map from
$C_0$ to $C$; therefore $K_{\hat X}\cdot C_0 = K_X\cdot C +E\cdot C_0 = -1$. In
particular this
implies that $C_0$ moves at least in a $1$-dimensional family on $\hat X$
(see (1.4.1)); since it
does not move on $E$ it means that it goes out of $E$. Since $C_0$ is
contracted by $\f \circ \psi$ it implies that $E.C_0 = 0$,
but this is a contradiction since $E.C_0 = 1$.

\medskip

For crepant contractions we have the following useful results.

\proclaim Proposition (5.6.2). ({\rm see [C-K-M, (16.6)]})
Let $\f: X\ra Z$ be a crepant contraction
of a smooth $3$- fold $X$ with an irreducible (reduced)
1-dimensional fiber $f:=\f^{-1}(z)_{{\rm red}}$.
Then $f\iso \P^1$ and the conormal bundle
$N^*_{f/X}=I_{f}/I^2_{f}$
is isomorphic to either
$\O(1)\oplus\O(1)$ or to $\O\oplus\O(2)$
or to $\O(-1)\oplus\O(3)$.
In the first two cases the ideal $\I_{f}$ is spanned by
global functions on $X$. More precisely
$\I_{f}=\f^{-1}m_z\cdot \O_X$.

\proof The description of $N^*_{f/X}$ follows from the
above arguments using the vanishing in (1.2.1) and the fact
that $K_X \cdot f = 0$. The last part of the proposition follows
from (5.4) (see also the proof of the next proposition).

\proclaim Proposition (5.6.3). Let $\f: X\ra Z$ be a crepant contraction
of a smooth $3$- fold $X$ with a (reduced)
1-dimensional fiber $f$ consisting of two components,
$f=f_1\cup f_2$. Then each of the components
is a smooth $\P^1$, they meet transversally at one point and the conormal
$N^*_{f/X}=I_{f}/I^2_{f}$ is locally free and restricted to
the component $f_i$ it is isomorphic to either $\O\oplus\O(1)$ or
$\O(-1)\oplus\O(2)$.
If the restriction to both components is $\O\oplus\O(1)$
then  $\I_{f}=\f^{-1}m_z\cdot \O_X$.

\proof The description of the  structure of $f$ is clear and follows
from (1.5.1). The fiber then is locally complete
intersection in $X$ and the conormal is locally free.
Then, as above
we get $H^1(f,\I_f/\I^2_f)=0$. This, because of the restriction
to $f_i\subset f$, yields the vanishing of
$H^1(f_i,(N^*_{f/X})_{|f_i})$ and thus it implies the splitting type of
$N^*_f$ on each $f_i$.

The global structure of $N^*_f$ over $f$ is determined by its
restriction to each of $f_i$ and also by the glueing above
the common point $x_0:=f_1\cap f_2$.
Geometrical meaning of this is as follows:
if $F_i=\P((N^*_{f/X})_{|f_i})$ is the appropriate ruled surface
over $f_i$ then $F:=\P(N^*_f)$ is obtained by glueing
of these two surfaces along the fiber over $x_0$. There are two
possibilities: either the two negative sections of each of them
meet together over $x_0$  or they do not meet.
For example, if the restriction to each of the component is
$\O\oplus\O(1)$ then either the bundle is globally decomposable
or not. The two cases can be cohomologically distinguished when we twist
the bundle by a line bundle which on each of the components
is $\O(-1)$. Then the decomposable bundle has both cohomology of
dimension 1 while both cohomology vanish if the glueing yields a
non-decomposable bundle.
We  note that in both cases the bundle $N^*_f$  is spanned by global
sections.

The rest of the proposition follows then from (5.4) noting that
the variety $\hat X$ has a unique singular point
$\hat x_0$ over $x_0$ and that the singularity of $\hat X$ at $\hat x_0$
is of the quadric cone, in particular it is terminal.

\remark (5.6.4).
In the previous proof we have pointed out 6 types of the normal
bundle $N_f$ --- depending on the splitting type on each of the components
and on the glueing over the common point. Among these, two types
with a ``double'' splitting $\O\oplus\O(-1)$ are similar to
the types $\O(-1)\oplus\O(-1)$ and $\O\oplus\O(-2)$ in the irreducible case.
Let us also note that the cohomology yields non-existence of reducible
fiber and decomposable normal bundle with double splitting type
$\O(1)\oplus\O(-2)$
(because then the conormal has nonzero 1-st cohomology).
To the authors' knowledge the relation between the splitting type and the singularities
of the general plane section is not yet
completely understood  (see [Re], [Ka-Mo], [Ka2]).

\medskip
From (2.5) we get the following more general

\proclaim Corollary (5.6.5). Let $C$ be a reduced
1 dimensional fiber of a good or crepant contraction of a smooth manifold
$X$. Then the following properties are equivalent:
\item{(a)} $N^*_{C/X}$ is nef,
\item{(b)} $N^*_{C/X}$ is spanned at a generic point of any component
of $C$,
\item{(c)} $N^*_{C/X}$ is spanned everywhere on $C$.

\bigskip \noindent
{\bf (5.7)  The normal bundle of a two dimensional fiber.}
\par\noindent
In order to understand higher dimensional fibers of good contractions we will
slice them down. Thus we will need some kind of "ascending property".

Suppose that $\f: X\ra Z $ is a good contraction of a smooth variety,
$\L$ is an ample line bundle on $X$ such that $-K_X-\L$ is $\f$-big and nef.
Let $F=\f^{-1}(z)$ be a (geometric) fiber of $\f$. Suppose that
$F$ is locally complete intersection. Let $X'\in |\L|$ be a
normal divisor which does not contain any component of $F$.
Then the restriction of $\f$ to $X'$, call it $\f'$,
is a contraction, either good or crepant (see (1.3.2) and (4.4)).
The intersection $F'=X'\cap F$ is then a fiber of $\f'$.
The regular sequence of local generators $(g_1,\dots,g_r)$
of the ideal of the fiber $F$ in $X$ descends
to a regular sequence in the local ring of $X'$ which defines
a  subscheme $F\cdot X'$ supported on $F'=F\cap X'$,
call it $\bar F'$. Let us note that 
if the divisor $X'$ has multipicity 1 along each of the components of
$F$ then, since a locally complete intersection has no
embedded components, we get $\bar F'=F'$.

\proclaim Lemma (5.7.1).  The scheme $\bar F'$
is locally complete intersection in $X'$ and
$$N^*_{\bar F'/X'}\otimes_{\O_{\bar F'}}\O_{F'}\iso (N^*_{F/X})_{|F'}.$$
If moreover $X'$ is smooth, $\L$ is spanned
and $dimF'=1$ then
$H^1(F',(N^*_{F/X})_{|F'})=0$.

\proof The first part of the lemma follows from the preceding discussion
so it is enough to prove the vanishing.
Let $\J$ be the ideal of $\bar F'$ in $X'$.
From (1.2.1) we know that $H^1(\bar F',\O_{X'}/\J^2)=0$ and since
we have an exact sequence
$$0\raa \J/\J^2=N^*_{\bar F'/X'}\raa
\O_{X'}/J^2\raa \O_{X'}/\J=\O_{\bar F'}\raa 0$$
then we will be done if we show
$H^0(\bar F',\O_{\bar F'})=\C$.
Since $H^1(\bar F',\O_{\bar F'})=0$ then this is equivalent to
$\chi(\O_{\bar F'})=1$.
The equality $H^0(\bar F',\O_{\bar F'})=\C$ is clear
if $\bar F'$ is reduced. But since $\L$ is spanned and $F$ is locally complete
intersection
then there exists a flat deformation of $\bar F'$ to another
intersection $F\cdot X''$ which is reduced. This is what we need,
because flat deformation
preserves Euler characteristic.
\medskip
Now let us consider the following ascending property.
Let us consider a point $x\in F'$. Suppose that the ideal of $F'$,
or equivalently $N^*_{F'/X'}$, is generated by global functions from
$X'$. That is, there exist global functions $g'_1, \dots g'_r\in
\Gamma(X',\O_{X'})$
which define $F'$ at $x$. Then, since $H^1(X,-\L)=0$ these functions extend
to $g_1, \dots g_r\in \Gamma(X,\O_{X})$ which define $F$. Thus passing from
the ideal $\I$ to its quotient $\I/\I^2$ we get the first part of

\proclaim Lemma (5.7.2).
If $N^*_{F'/X'}$ is spanned by global functions from $\Gamma(X',\O_{X'})$ at
a point $x\in F'$ then $N^*_{F/X}$ is spanned at $x$ by functions from
$\Gamma(X,\O_{X})$.
If $N^*_{F'/X'}$ is spanned by global functions from $\Gamma(X',\O_{X'})$
everywhere on $F'$ then $N^*_{F/X}$ is nef.

\proof We are only to proof the second claim of the Lemma. Since
$F'\subset F$ is an ample section then the set where $N^*_{F/X}$ is not
generate by global sections is finite in $F$. Therefore the restriction
$(N^*_{F/X})_{|C}$ is spanned generically for any curve $C\subset F$ and
consequenly it is nef.

If the fiber is of dimension 2 then we have a better extension property.

\proclaim Lemma (5.7.3). Let $\f:X\ra Z$ be a good
birational contraction of a smooth variety
with an isolated 2-dimensional
fiber $F$ which is a locally complete intersection.
As usually $L=-K_X$ is a $\f$-ample line bundle which can be
assumed $\f$-very ample (see (1.3.4)).
Then the following conditions are equivalent:
\item{(a)}$N^*_{F/X}$ is generated by global sections at any point of $F$
\item{(b)} for a generic (smooth) divisor $X'\in |L|$ the bundle
$N^*_{F'/X'}$ is generated by global sections at a generic point of
any component of $F'$

\proof The implication (a)$\Rightarrow$(b) is clear. To prove the converse
we assume the contrary. Let $S$ denote the set of points on $F$ where
$N^*_{F/X}$ is not spanned. Because of the extension property (5.7.2)
and Corollary (5.6.5)
the set does not contain $F'$ and thus it is finite. Now we choose another
smooth section
$X'_1\in |L|$ which meets $F$ along a (reduced)
curve $F'_1$ containing a point of $S$.
(We can do it because $L$ is $\f$-very ample.)
The bundle $N^*_{F'_1/X'_1}$
is generated on a generic point of $F'_1$ so it is generated everywhere but
this, because
of the extension property, implies that $N^*_{F/X}$ is generated at some
point of $S$,
a contradiction.

\proclaim Lemma (5.7.4).
Let $\f:X\ra Z$ be a good birational contraction of a smooth 4-fold with an
isolated
2-dimensional fiber $F=\f^{-1}(z)$. As usually $L=-K_X$ is a $\f$-ample
line bundle which
may be assumed to be $\f$-very ample.
Then the fiber structure $\tilde F$ coincides with the geometric structure $F$
 unless one of the following occurs:
\par{\sl
\item{(a)} the fiber $F$ is irreducible and the restriction of
$N_F$ to any smooth curve $C\in |L_{|F}|$ is isomorphic to
$\O(-3)\oplus\O(1)$,
\item{(b)} $F=\P^2\cup\P^2$ and the restriction of
$N_F$ to any line
in one of the components is isomorphic to
$\O(-2)\oplus\O(1)$.
}

\proof  Let us consider an arbitrary curve $C\in |L_{|F}|$.
Since $L$ is $\f$-very ample we can take a smooth $X'\in |L|$ such that
$\f'=\f_{|X}$ is a crepant contraction and $C=F\cap X'$. Then considering the
embeddings $C=F\cap X'\subset F\subset X$ and $C=F\cap X'\subset X'\subset X$
$$N_{C/X}=N_{C/X'}\oplus L_C=(N_{F/X})_{|C}\oplus L_C$$
and therefore $N_{C/X'}=(N_{F/X})_{|C}$. Now we apply the propositions (5.6.2)
and (5.6.3) to describe $N_{C/X'}$.
In particular it follows that if neither (a) nor (b)
occurs then the fiber structure of the contraction $\f'$ is reduced.
Thus, using our ascending proposition (5.7.3) and
the equivalence in (5.4), we conclude that $\tilde F=F$.

\medskip

Now, let us discuss the possible exceptions described in the above lemma.
If $F=\P^2$ then, because of the theorem of Van de Ven (see (2.7)), $N_{F/X}$
would be decomposable and in particular $h^0(N_{F/X})-h^1(N_{F/X})>0$.
Because of Lemma (2.11) the possible exception over an irreducible quadric would
satisfy the same inequality (we use the vanishing in (5.7.1)).
This, however, because of the theory of
deformation see, e.g. [Ko] would imply that $F$ moves in $X$ which
contradicts our assumption that $F$ is an isolated 2-dimensional fiber.
A similar argument works for a reducible quadric. Namely, because of (5.7.4.(b))
we can apply
the theorem of Van de Ven to claim that for a component $F'$ of $F$ we have
$(N_{F/X})_{|F'}=\O(-2)\oplus\O(1)$. Thus, in view of (2.3) if we blow-up
the other
component and consider the strict transform of $F'$ then its normal would be
$\O(-3)\oplus\O$. Now we see that the strict transform would move which is
clearly
impossible. Therefore we have proved

\proclaim Theorem (5.7.5).
Let $\f:X\ra Z$ be a good birational contraction of a smooth 4-fold with an
isolated
2-dimensional fiber $F=\f^{-1}(z)$.
Then the fiber structure $\tilde F$ coincides with the geometric structure $F$
and the conormal $N^*_{F/X}$ is spanned by global sections.

Now we can verify which one among spanned vector bundles with the
appropriate $c_1$
is actually the conormal of an isolated 2-dimensional fiber.
The result is the following

\proclaim Theorem (5.7.6).
Let $\f:X\ra Z$ be a good birational contraction of a smooth 4-fold with an
isolated
2-dimensional fiber $F=\f^{-1}(z)$. If $F=\P^2$ then $N^*_{F/X}$ is either
$\O(1)\oplus\O(1)$ or $T(-1)\oplus\O(1)/\O)$, or
$\O^{\oplus 4}/\O(-1)^{\oplus2})$.
If is a quadric (possibly singular or even reducible) then $N^*_{F/X}$ is
the spinor
bundle ${\cal S}(1)$.

\proof
We use the classification results in (2.6) and (2.8)
together with the following observations.
Because of the deformation argument we know that
$h^0(N_{F/X})-h^1(N_{F/X})\leq 0$
and thus $N^*_{F/X}$ can not be decomposable with a trivial factor. On the
other hand
none of the bundles with $s_2=c_1^2-c_2=0$ can occur as $N^*_{F/X}$ because
of the following

\proclaim Lemma (5.7.7). Let $\f$ be a good birational contraction from a smooth
$4$-fold and let $F$ be an isolated two dimensional fiber.
Then $s_2(N^*_{F/X}) > 0$.

\proof Assume by contradiction that $s_2 = 0$.
Then from the classification of such bundles, (2.6) and (2.8), we see
that
$H^0(S^n(N^*_{F/X})) = S^n(H^0(N^*_{F/X}))$ is of dimension
$ {\big(\ }^{n+3}_{\ \ 2\ \ }\big)$.
Thus, by (2.4) the contraction should be to a a 3-dimensional smooth point
contrary to the fact that $dimZ=4$.

\remark The cases with $s_2 = 0$ do occur in the fiber type contractions;
see examples (3.5.5).

\bigskip

In some respect the above results about the fiber structure of a
2-dimensional fiber are nicer than one may expect. Namely, there is
no multiple fiber structure, the conormal is nef and the normal of the
geometric isolated fiber has no section. Thus the situation is better than
for 1-dimensional isolated fibers in dimension 2 and 3: the
fundamental cycle of a Du Val $D-E$ surface singularity is
non-reduced and in dimension 3 one may contract an isolated $\P^1$
with the normal $\O(1)\oplus\O(-3)$.  On the other hand, using the
double covering construction (see the Examples section, (3.5))
in dimension 5 one may contract a
quadric fibration over a smooth 3-dimensional base with an isolated
fiber equal to $\P^2$, scheme theoretically the fiber is a double
$\P^2$. Using the sequence of normal bundles and the deformation
of lines argument, one may verify that in this case $N_F\iso \O(1)\oplus\E^*$
where $\E$ is a rank 2 spanned vector bundle with $c_1=2$, $c_2=4$,
so that $dimH^1(\E^*)=-\chi(\E^*)=3$.

\smallskip
Let us also note that for a  divisorial fibers we have the following:

If $F=\bigcup F_i$ is a divisorial fiber of a surjective map
$X\ra Y$, where $X$ is smooth and $dimY\geq 2$ then for some $k>0$
the line bundle $\O_{F_i}(-kF_i)$ has non-trivial section and thus
no multiple of $\O_{F_i}(F_i)$ has a section.
In particular, if $rank(Pic(X/Y))=1$ then $\O_F(-F)$ is ample.

\smallskip
One can then try to conjecture that if $F$ is an isolated fiber
of a (good) contraction
which is locally complete intersection and with ``small'' codimension
then $H^0(F,N_F)=0$.

\bigskip \noindent
(5.8)
The above result on contractions of 4-folds can be generalized for
the adjunction mappings of an $n$-fold. Namely, suppose that
$\f: X\ra Y$ is a good contraction of a smooth $n$-fold $X$ supported
by a divisor $K_X+(n-3)H$, where $H$ is a $\f$-ample divisor on $X$.
Since we are interested in the local desription of $\f$ around a non-trivial
fiber $F=\f^{-1}(z)$, we may assume that the variety $Z$ is affine.

\medskip
\proclaim Corollary (5.8.1).
Let us assume that $\f$ is birational and that $F$
is an isolated fiber of dimension $n-2$. If $n\geq 5$ then
the contraction $\f$ is small and
$F$ is an isolated non-trivial fiber. More precisely
$F\iso\P^{n-2}$ and $N_{F/X}\iso\O(-1)\oplus\O(-1)$, and there exists
a flip of $\f$ (see [A-B-W]).

\proof
The proof follows by easy slicing method:
by [A-W], (5.1),  the line bundle $H$ is $\f$-spanned,  so there exists a
smooth
hypeplane section
$\in |H|$ which has the same properties as $X$.
Indeed, the restriction of $\f$ to the hyperplane section is a good
contraction (see [A-W], (2.6))
supported by $K+(n-4)H$ and, because of Bertini theorem,
the hyperplane section of $F$ is an isolated fiber of dimension $n-3$.
Thus we can arrive to the 4-dimensional linear section $X'$ of $X$ which we
know by (5.7.6). In particular we know the surface section $F'$ of the fiber
$F$ and its normal $N_{F'/X'}$. As in (5.7.1) we note that $N_{F'/X'}\iso
(N_{F/X})_{F'}$. Therefore it is enough to verify which of the pairs
$(F',N_{F'/X})$ is a plane section of a higher-dimensional pair.
Also, we note that because of Lemma (5.7.2) $N^*_{F'/X'}$ extends to a nef
vector bundle.

If $F'\iso\P^2$ then it must be a section of the projective space.
But then, among the bundles occuring in (5.7.6), only the bundle $\O(-1)\oplus
\O(-1)$ and $T\P^2(-1)\oplus\O/\O(-1)$ extend to nef vector bundles on
$\P^3$ (see [Sz-W2]).
The bundle $T\P^2(-1)\oplus\O/\O(-1)$ extends to the
null-correlation bundle on $\P^3$  and this one can occur as the normal of
an isolated 3-dimensional fiber of a fiber type contraction. Thus,
by a similar argument as in the proof of (5.7.7)
this case can not occur
in the case of a birational contraction.
Indeed, because of the nefness we have
the vanishing of $H^1(F,S^k(N^*_{F/X}))$ and
since $F$ with such a normal
is a fiber of a good contraction to a smooth
4-dimensional point it follows that we can apply (2.4).
The existence of flip and other results are proved in the
Lemma (6.1) of the next section.

A similar argument works if $F'$ is a quadric.
First we note that $F$ can not
be reducible quadric since $T\P^2$ does not extend.
A quadric $\Q^3$ with the normal ${\cal S}$ can be a fiber
of a fiber type contraction to a smooth point and again, by cohomological
argument, which actually depends only on the nefness and
Chern classes of the bundles in question,
it can not be a fiber of a birational contraction.

\bigskip

\noindent (5.9) {\bf The case of a conic fibration.}
\par\noindent
Note that the preceding arguments, which led to the classification
of the birational 4-dimensional case, depend on the isomorphism
$\f'_*\O_{X'}\iso\O_Z\iso\f_*\O_X$. This fails to be true if
$\f$ is of fiber type. Namely, let $\f:X\ra Z$ be a conic fibration,
i.e.~a good contraction such that $dimZ=dimX-1$. As usually,
we will assume that $F$ is an isolated
 2 dimensional fiber of $\f$ and $L=-K_X$ is 
$\f$-spanned. Then the restriction of $\f$ to
a general section $X'\in |L|$ is generically $2:1$ 
covering of $Z$. Let us assume that $X'$ is connected. The push-forward
of the divisorial sequence for $X'$ yields the exact sequence
$$0\raa\f_*\O_X=\O_Z\raa\f_*\O_{X'}\raa R^1\f_*\O_X(K_X)\raa 0\eqno(5.9.1)$$
of $\O_Z$-modules. Moreover, by (1.2.4) $R^1\f_*\O_X(K_X)=\omega_Z$.
If $\f$ is an elementary contraction then $Z$ is algebraically factorial
hence $\omega_Z=K_Z$ is Cartier and $Z$ is Gorenstein. 
Thus $\f_*\O_{X'}$ is locally free $\O_Z$-module.
\par
Let $\f = \pi\circ\f'$ be the Stein factorization of $\f$, where the map 
$\pi :Z' \ra Z$ is a $2:1$ cover and $\f'$ a crepant contraction 
(see (4.4)) so that $Z'$ is Gorenstein.
Since $\pi_*\O_{Z'}=\f_*\O_{X'}$ is locally free, it follows that
$\pi$ is flat.
We consider the {\sl trace map}, $tr: \pi_*\O_{Z'} \ra \O_Z$,
which over an open $U\subset Z$ is defined as it follows 
(for details see [A-K], p.~123).
The sheaf $\pi_* \O_{Z'}$ is free over $\O_Z$ and every element $f$ in
$\pi_* \O_{Z'}(U)$ defines a $\O_Z(U)$-homomorphism of the free 
$\O_Z(U)$-module $\pi_* \O_{Z'}(U)$; we define $tr(f)$ to be the 
trace of this homomorphism. The map $tr$ splits (5.9.1).
\par
Another proof of the splitting of (5.9.1), independent on
the assumption that $\f$ is elementary, is provided by 
[Ko2, Cor.~2.25]. We note that also in this case, since $\f_*\O_{X'}$
is reflexive (because $Z$ is normal and $\f_{|X'}$ contracts no divisor)
thus, by the splitting, $R^1\f_*\O_X(K_X)$ is reflexive, hence invertible.
Therefore, as before, $Z$ is Gorenstein and $\pi$ is flat.
\par
So $\pi$ is a double cyclic cover. 
The kernel of its trace map is equal to $\pi_*\I_R=\I_B$, where
$R$ and $B$ are, respectively, ramification and branching 
divisors of $\pi$, and both are Cartier.
\par
Now suppose that $F$ is as in (5.0), i.e.~it is locally complete intersection 
and the blow-up of $X$ at $F$ has log terminal singularities. Let, as
in (5.7), $F'=F\cap X'$ and set $N^*=N^*_{F/X}$. 
We note that Lemma (5.7.1) remains true also in
this case and in particular we have the restriction
$res: H^0(F,N^*)\ra H^0(F',N^*_{|F'})=H^0(F',N^*_{F'/X'})$. 
Suppose that $N^*_{F'/X'}$ is nef. Then, acccording to (5.3.b)
any section of $N^*_{F'/X'}$ extends to a section of $X'$ and we can
apply $tr$ to this extension. As the result we obtain a well defined
map $tr: H^0(F',N^*_{F'/X'})\ra H^0(F, N^*_{F/X})$
such that $res\circ tr=id$. Since $R$ is Cartier and it gives the kernel
of $tr: \pi_*\O_{Z'}=\pi_*(\f'_*\O_{X'})\ra\O_Z$ it follows

\proclaim Lemma (5.9.2).
The kernel of $tr: H^0(F',N^*_{F'/X'})\ra H^0(F,N^*_{F/X})$
is of dimension one at most and it is generated by the class
of $f_R\circ\f'$ where the function $f_R$ generates the ideal of $R$ at 
$z'=\f'(F')$.

Now, using Lemma (5.9.2) we investigate the structure of the 
conic type contraction in dimension 4. 
We will concentrate on the case when $F$ is a projective space or
an irreducible quadric.

\proclaim Lemma (5.9.3).
Let $\f:X\ra Z$ be a conic fibration of a smooth 4-fold.
Suppose that $F$ is an irreducible isolated
2-dimensional fiber of $\f$ such that $(F,L_F)$ is either $(\P^2,\O(1))$, or $(\F_0,C_0+f)$, or $(\S_2,\O(1))$. 
Then either $N^*_{F/X}$ is nef or it is one of the
exceptional bundles: either (2.7.0) or $\O(2,0)\oplus\O(-1,1)$, or
(2.11.0), respectively.

\proof We use the notation introduced above, i.e.~$X'\subset X$ is a 
linear section of $X$, $F'=X'\cap F$ and so on. In view of the vanishing
(1.2.2) and Lemma (4.7.1), the first cohomology of the 
restriction of $N^*:=N^*_{F/X}$ to $F'\subset F$ vanishes. 
In view of the results of Section 2 (2.7--2.11), 
either the general splitting type of $N^*$ is $(1,1)$ or $N^*$
is among some very special exceptions. 
\par
We claim that if for a general smooth $F'\subset F$ there is
$N^*_{|F'}\iso\O(1)\oplus\O(1)$ then
the bundle $N^*$ is nef.
Let  consider
the image of trace $\Lambda:=im(tr)\subset H^0(F,N^*)$
which gives rise to (at least) 2 dimensional linear system $|\Lambda|\subset |H^0(\P(N^*,\O_{\P(N^*)}(1)))|$. 
Since $N^*_{|F'}\iso\O(1)\oplus\O(1)$ thus
over $F'\subset F$ the system $|\Lambda|$ has at most one base point.
Therefore, if a curve $l_0\subset\P(N^*)$ has negative intersection with
$\O_{\P(N^*)}(1)$ then it descends to a line $l\subset F$.  Because
of the vanishing (1.2.2) we have $N^*_{|l}\iso\O(-1)\oplus\O(3)$ or
$\O(-1)\oplus \O(2)$, depending on whether $F$ is $\P^2$ or a
quadric.  Therefore $l_0\cdot \O_{\P(N^*)}(1)=-1$. Now the argument
is similar to the one in the proof of (5.6.1). Namely, let us
consider a blow-up $\hat X$ of $X$ along $F$. We compute that
$l_0\cdot K_X=0$ and thus, because of (1.4.1), $l_0$ moves in one
dimensional family in $\hat X$ (if $F=\S^2$ then $\hat X$ is
singular and one has to use the full strength of [Ko, II.1.14]).  This
implies that deformations of $l_0$ sweep out a divisor in $\P(N^*)$
which dominates
$F$ and contradicts the fact that the base locus of $|\Lambda|$ over
$F'$ has one point only.
\par
Now it remains to verify which of the bundles listed in (2.7.1) and (2.11)
are good candidates for $N^*$ in case when it is not nef. The decomposable 
bundles on $\P^2$ as well as some on the quadrics do not verify
$h^0(N)-h^1(N)\leq 0$ condition. Moreover, the bundles presented
in (2.11) (a-ii) and (b-ii) contain only a finite number of lines
with splitting type $(-1,2)$ and thus are excluded by the above 
deformation argument. Thus we are left with the exceptions which
are listed in the present lemma.
\medskip 

A similar argument proves the following
\proclaim Lemma (5.9.4). In the situation of (5.9.3) if $(F, L_F)=
(\P^2,\O(2))$ then $N^*\iso T\P^2(-1)$.

\proof In view of (2.7) it is enough to prove that the splitting type
of $N^*$ on any line is $(0,1)$. (It is easy to check that if   $N^*\iso\O(1)\oplus\O$ then $F$ is a non-isolated 2-dimensional fiber.)
If it is not the case, then as above
we consider the blow-up $\hat X$ and in $\hat X$ over a line $l\subset X$
we have a rational curve $l_0$ such that $l_0\cdot K_{\hat X}=-1$.
Now $l_0$ moves in a 2 dimensional family which would imply that
all lines in $F$ have splitting type $(-1,2)$ and thus 
$N^*\iso\O(-1)\oplus\O(2)$ which contradicts $h^0(N)-h^1(N)\leq 0$ condition.

\medskip

Now we deal with the exceptional bundles which are singled out in (5.9.3).
First, however, let us make an observation concerning its
proof. We use the notation 
from the proof of (5.9.3), i.e.~we consider a smooth 
section $F'\subset F$,
and we assume that $N^*_{|F'}\iso\O\oplus\O(2)$.
Let us consider a 2-dimensional linear system 
$|\Lambda|\subset|H^0(\P(N^*),\O(1))|$ which is the image of $tr$.
Then the argument, 
which shows the nefness of $N^*$, fails only if
the base point locus of $|\Lambda|$ over $F'$ is the exceptional $(-2)$-curve 
$C_0\subset\P(N^*_{|F'})\iso\F_2$.  
This is possible only if the
kernel of $tr: H^0(F',N^*_{|F'})=H^0(\P(N^*_{|F'}),\O(1))
\ra H^0(F,N^*)$ contains a section, zero locus of which
does not meet $C_0$. 
In other words, suppose that we take $\beta':\hat X'\ra X'$,
the blow-up of 
$X'$ along $F'$ with the exceptional divisor $\hat F'$,
and we take the strict transform
$\hat R$ of the ramification divisor
$R=\{f_R=0\}$ then $\hat R\equiv (\beta'\circ\f')^*(R)-\hat F'$
and $\hat R\cap \hat F'$ does not meet $C_0$
--- see Lemma (5.9.2).

Now,  $f_R^2$ is in $\O_Z$ and it defines the branch divisor
$B$. Let us note that $\f^*B\cdot X'=2\f'^*R$. We take the blow-up
$\beta:\hat X\ra X$ of $X$ along $F$ with the exceptional divisor
$\hat F$.  Then, the strict transform $\hat B$ is equivalent to 
$-2\hat F$ and over $F'$ it is equal to $2(R\cdot \hat F')$.
But --- as we noticed in the proof of (5.9.3) --- the base point locus
of $-2\hat F_{|\hat F}=\O_{\P(N^*)}(2)$ contains a divisor 
sweapt out by curves which have negative intersection with
$-\hat F_{|\hat F}$ and, as it follows from the description of 
the bundles in question, over $F'$ the base point locus contains
the curve $C_0$. This is in contradiction with our observation
that $\hat R\cap \hat F'\cap C_0=\emptyset$. 
Thus we have proved 

\proclaim Proposition (5.9.5).
Let $\f:X\ra Z$ be a conic fibration of a smooth 4-fold.
Suppose that $F$ is an irreducible isolated
2-dimensional fiber of $\f$ which is either a projective
plane or a quadric. Then the conormal bundle $N^*_{F/X}$
is nef.

Since $N^*_{F/X}$ is nef then one can use the results presented in the beginning
of the section to describe the contraction $\f$ around the fiber $F$. 
Indeed, by the results of [S-W1] (c.f.~(5.7.3)), nefness of such bundles implies
their spannedness and one can use their explicit classification (see~(2.6) and 
(2.8)).  In particular one gets the following result, the proof of which is
similar to that of (5.7.6) (see also Section 6):

\proclaim Theorem (5.9.6). 
Let $\f:X\ra Z$ be a conic fibration of a smooth 4-fold.
Suppose that $F=\f^{-1}(z)$ is an irreducible isolated
2-dimensional fiber. If $F\iso \P^2$ then
$N^*_{F/X}\iso \O^3/\O(-2)$. If $F$ is an irreducible quadric
then $N^*_{F/X}$ is the pullback of $T\P^2(-1)$ via some double covering
of $\P^2$. In both cases the fiber structure $\tilde F$ coincides with 
the geometric structure on $F$ and $Z$ is smooth at $z$.

The remaining cases of 4 dimensional conic fibrations with isolated 2 dimensional fibers,
including  these which are non-reducible or
non-locally complete intersection, will be treated separately in our 
forthcoming paper.

\beginsection 6. Geometry of a contraction of a 4-fold.

(6.0) In this section $\f: X \ra Z$ will be a birational contraction of a smooth
4-fold with
an isolated two dimensional fiber $F=\f^{-1}(z)$.
In the previous section we described the normal of $F$. We proved
that $N^*_{F/X}$ is spanned by global sections and thus we have a
contraction map over $Z$ of the blow-up $\hat X$ of $X$ at $\F$
to $\hat Z$, the blow-up of $Z$ at $z$ (see (5.4)).
Over the exceptional divisor $\hat F$ the map $\hat \f$ is
associated to the evaluation of sections of $N^*_{F/X}$.  Therefore,
knowing the list of possible normal bundles, we can deduce the
description of $\hat \f$. In particular, we note that $\hat \f$ is an
isomorphism if and only if $F=\P^2$ and $N^*_F=\O(1)\oplus\O(1)$
which is exactly the case when $\f$ is a small contraction. Indeed,
we have a more general

\proclaim Lemma (6.1). Let $\f: X\ra Z$ be a good or crepant contraction of a
smooth variety with a fiber $F=\f^{-1}(z)$.  Assume that $F$ is
locally complete intersection and the blow up $\beta:\hat X\ra X$ of
$X$ along $F$ has log terminal singularities.
If $N^*_{F/X}$ is ample then the exceptional
locus of $\f$ is equal to $F$. If $F$ is not a divisor and $\L$ is a
$\f$-ample line bundle then a flip of $\f$ is defined by a contraction
of $\hat X$ supported by $-\hat F-\tau\L$, where $\tau$ is a rational
number such that $-\hat F-\tau\L$ is nef but not ample.

\proof Outside of $F$ and, respectively,  $\hat F$ the exceptional loci
of $\f$ and $\hat \f$ coincide. On the other hand, because of the local
nature of our set-up, $\hat \f$ is an isomorphism on $\hat X$ if and only
if $N^*_F$ is ample. This proves the first statement.
The last part of the lemma is clear.

\medskip
The above lemma can not be inversed since in dimension 3 we have a
contraction of an isolated $\P^1$ with normal $\O\oplus\O(-2)$ or
$\O(1)\oplus\O(-3)$.  In our case, however, we can by using either the
result of Kawamata [Ka1] or applying a direct argument based on the
deformation theory, especially the inequality in [Wi] (Theorem (1.1)).
More precisely, using this inequality one proves:

\proclaim Lemma (6.2). Let $\f:X\ra Z$ be a good contraction of
a smooth variety with an irreducible fiber $F$ which is
locally complete intersection in $X$. Suppose that $\hat X$ is smooth,
$N^*_F$ is nef (therefore semiample) and that the map of $\P(N^*_F)$
associated to high multiple of $\O_{N^*_F}(1)$ is birational with
general non trivial fiber of dimension 1.
Then the exceptional set of $\f$ contains a divisor.

\proof By our assumptions the map $\hat\f$ supported by $-\hat F$
is birational with some 1-dimensional fibers and therefore, by dimension
estimate (see [Wi]) the exceptional set of $\hat\f$
contains a divisor which is not
$\hat F$ (because the map is birational on $\hat F$).

\proclaim Corollary (6.2.1).  Let $\f:X\ra Z$ be a birational
good contraction of a 4-fold with an isolated two dimensional fiber
$F$. Then $\f$ is small if and only if $F =\P^2$ and $N_{F/X} =
\O(-1)\oplus\O(-1)$.

\medskip

Now we assume that $n =4$ and that the contraction $\f:X\ra Z$
is divisorial with an exceptional set $E$ mapped to a
surface $S$ and the isolated 2-dimensional fiber $F=\f^{-1}(z)$.
Using the results from the previous section we can consider
$\beta:\hat X\ra X$,
the blow-up of $X$ along $F$ with the exceptional divisor $\hat F$,
$\alpha: \hat Z\ra Z$,
the blow-up of $Z$ at the maximal ideal of $z$ with the exceptional
set $\hat G$, and the map $\hat\f: \hat X\ra\hat Z$, which is a good
contraction supported by $-\hat F$. Let $\hat E$ and $\hat S$ denote,
respectively, the strict transforms of $E$ and $S$.
All the objects are presented on the following diagram.
$$\matrix{
&&&\hat X\supset (\hat E,\hat F)&&&\cr
&&\beta\swarrow&&\searrow\hat\f&&\cr
&&&&&&\cr
&X\supset (E\supset F)&&&&\hat Z\supset (\hat S,\hat G)&\cr
&&&&&&\cr
&&\f\searrow&&\swarrow\alpha&&\cr
&&&Z\supset (S\ni z)&&&}\leqno (*)$$

We will provide a description of all these objects and, in particular,
the singularities of $Z$ at $z$.
Actually, we can get the description of $Z$ at $z$ by using (2.4) and (5.5),
however, we find it interesting to provide a complete description of
all the objects
which occur in this blow-up-blow-down construction.
\bigskip

\noindent
{\bf Proposition (6.3).}
{\sl  Let $\f: X \ra Z$ be a divisorial contraction
of a $4$-fold $X$ with an isolated two dimensional fiber $F \iso \P^2$.
Then $\hat Z$ and $\hat S$ are smooth varieties and $\hat\f$
is the blow-up of $\hat Z$ along $\hat S$.
\item{(a)} If $N^*_{F/X}= T(-1)\oplus\O(1)/\O)$
then $\hat G$ is a smooth quadric
$\Q^3$ and $N_{\hat G/\hat Z}= \O_{\Q^3}(-1)$. The map $\alpha$ is
the contraction of $\Q^3$ to an isolated singular point $z\in Z$ which
is analytically isomorphic to the quadric cone singularity.
The surface $S\subset Z$ is smooth.
\item{(b)}
If $N^*_{F/X}= \O^{\oplus 4}/\O(-1)^{\oplus2}$ then $\hat G$ is equal to
$\P^3$ and $N_{\hat G/\hat Z} = \O_{\P^3}(-1)$. The
map $\alpha$ is a  blow-down of $\hat G$ to a smooth point $z\in Z$.
The surface $S$ has a singularity
in $z$ of the type of the vertex of a cone over a
twisted cubic.
\item{}
}

\smallskip \noindent
\proof The variety $\hat X$ is smooth and by Lemma (2.6)
the map $\hat\f_{|\hat F}$
is either the blow-up of a smooth three dimensional quadric along a line
or the blow-up of
$\P^3$ along a twisted cubic curve.
Consequently, the map $\hat\f$ is a divisorial contraction with all non
trivial fiber of dimension $1$. Therefore, by (6.1),
$\hat E$ is smooth and $\hat\f$ blows it down to a smooth $\hat S$
in smooth $\hat Z$.
\par If $m$ is the multiplicity of $E$ along
$F$, then we have the numerical equivalence
$$K_{\hat X}= \beta ^* K_X +\hat F = \beta^* E +\hat F = \hat E +
(m+1) \hat F$$
and also $K_{\hat X}= \hat\f^* K_{\hat Z} + \hat E$.
Therefore $\hat\f^* K_{\hat Z} = (m+1) \hat F$ and,
using the adjunction formula $(K_{\hat Z} +\hat G)_{|\hat G} = K_{\hat G}$,
we obtain
that $(m+2)\hat G_{|\hat G} = \O(-3)$ or $\O(-4)$ if
$\hat G$ is a quadric or $\P^3$, respectively.
This implies that the normal of $\hat G$ in $\hat Z$ is in both cases $\O
(-1)$ (which we know also from Corollary (5.3.1))
and also that $m = 1$, if $\hat G$ is a quadric, while $m = 2$ if $\hat G$
is $\P^3$
(this will be reproved later in (6.7)).
Since the above multiplicity coincides with the degree of the projection
$\hat E\cap\hat F\ra F$ (which we know because of the description of $N^*_F$)
it follows that $\hat E$ and $\hat F$ intersect transversaly.
The description of the singularity of $Z$ follows thus immediately
(see also (5.5)).
\par
Now we discuss the singularity of $S$ at $z$.
For this purpose we consider the
intersection curve $f:=\hat G\cap\hat S$. In both cases $f\iso\P^1$
and  the curve is either a line or a twisted cubic curve,
if $\hat G\iso\Q^3$ or $\P^3$,
$N_{f/\hat S}=\O(-1)$ or $\O(-3)$ if $G\iso\Q^3$ or $\P^3$, respectively.
This provides the description of the singularity of $S$.

\proclaim Proposition (6.4). Let $\f: X \ra Z$ be a divisorial contraction of
a $4$-fold $X$
with an isolated two dimensional fiber $F$ which is a, possibly singular
(even  reducible), quadric with conormal bundle equal to ${\cal S}(1)$.
Then $\hat Z$ is smooth, $\hat G = \P^3$ and $N_{\hat G/\hat Z} = \O(-1)$.
Thus $z\in Z$ is a smooth point and $\alpha$ is the blow-up of $z\in Z$.
The surface $S$ is smooth outside of $z$ and it 
has an isolated non-normal point at $z$.

\proof
Although we can make an argument similar to the one used in the previous
Proposition, it is definitely much more convenient to use Corollary (5.3.1).
First, however, we note that the blow-up of $X$ along $F$ has good
singularities so that we can apply that result. Now, we note that
$(Proj(\bigoplus_k H^0(F,S^k(N^*_F))),\O(1))\iso (\P^3,\O(1))$ ---
this follows  e.g.~from our examples (3.2) or can be verified directly.
Thus, because of (5.3.1), $\hat G\iso \P^3$, $\O_{\hat G}(-\hat G)\iso\O(1)$,
then $\hat Z$ is smooth along $\hat G$ and consequently, $Z$ is smooth at $z$.
The multiplicity of $E$ along $F$ can be computed similarly as in the case
of $F\iso\P^2$ --- provided that the fiber is ireducible. Indeed,
in such a case the singularities of $\hat X$ are $\Q$-factorial and one can
make a computation involving Cartier divisors to prove that the
multiplicity is 2. Also, in this case,
one can describe the singularities of $S$ in terms of $\hat S$: if
$F\iso \P^1\times\P^1$, then $\hat S$ is smooth (by the same argument as in
the case
of $F=\P^2$) and the intersection $\hat S\cap\hat G$
is transversal and consits of two disjoint lines. Therefore
the singularity of $S$ is of type of two tranversal planes in a 4-space.

If $F$ is a quadric cone, then the intersection $\hat S\cap\hat G$ consits
of a line along which $\hat S$ is not normal.
Indeed, we can take a hyperplane section $Y$
of the manifold $\hat Z$ and then its inverse image $\hat Y\subset \hat X$.
The variety $\hat Y$ is Gorenstein, it has an isolated singular points
(quadric cone singularities) above the vertex of $F$ and the contraction
$\hat Y\ra Y$ is a divisorial good contraction. Therefore by
Cutkosky [Cu] the contraction $\hat Y\ra Y$ is the blow-up of $Y$
along a curve $Y\cap \hat S$ which is locally intersection and has (planar)
singularities on $\hat G$. In other words,
the singularities of $\hat S$ are of the
double point type and the whole situation
is a degeneration of the two non-meeting
$(-1)$-lines which occur in the case of $F=\F_0$.

If $F\iso\P^2\cup\P^2$ then
$\hat X$ is not $\Q$ factorial and in fact the exceptional divisor
of the elementary contraction $\hat f:\hat X\ra\hat Z$ is reducible and
consists of the strict transform $\hat E$ and one of the components of
$\hat F$ which is the projectivisation of $T\P^2(-1)$.
The situation is better understood with an alternative approach based
on the study of a component of a Hilbert scheme of $X$. We will do it in the
remainder of this section, postponing the conclusion of the proof of the
Proposition until the end of the section.

\bigskip \noindent
(6.5) An alternative way of describing the geometry of
$\f$ around the special fiber is based on
the Hilbert scheme of lines in fibers of the good contraction. In
the situation of a birational contraction as in (6.0) a {\sl line} is a
proper rational curve over $Z$ whose intersection with $L$ is 1.
Suppose that the contraction $\f$ is divisorial and the exceptional
locus is an irreducible divisor $E$. The general fiber of $\f_{|E}$
is then a line with the normal $\O^{\oplus 2}\oplus\O(-1)$.  Let
${\cal H}$ be a component of the Hilbert scheme of $X$ over $Z$ which
parametrizes general fibers of $\f_{|E}$. By general properties of
Hilbert scheme the component ${\cal H}$ admits a morphism $\f_{{\cal
H}}: {\cal H}\ra Z$ such that $\f_{{\cal H}}([l])=\f(l)$ where $[l]$
represents a line $l\subset X$. We note that $\f_{{\cal H}}$ maps
${\cal H}$ birationally to the surface $\f(E)$. The exceptional set
of $\f_{{\cal H}}$ is over $z\in \f(E)$ and it parametrizes lines in
$F=\f^{-1}(z)$ which are limits of {\sl general lines}. We will say
that such lines can {\sl move out} of $F$. Since the incidence
variety of lines over ${\cal H}$ maps onto $E$ it follows that such
lines cover $F$ and thus $\f_{{\cal H}}^{-1}(z)$ is 1-dimensional.

\medskip

Suppose that $F_1\iso\P^2$ is a component of the fiber $F$. Lines in $F_1$
are parametrized by another component of the Hilbert scheme of $X$, we
call it $F_1^{\vee}$. The component $F_1^\vee$ is just a plane dual to
$F_1$, i.e. $F_1^\vee\iso\P^2$. The intersection $F_1^\vee\cap{\cal H}$
contains those lines which move out of $F_1$.

\proclaim Lemma (6.6).  Suppose that the conormal bundle $N^*_{F_1/X}$ is as
in (5.7.6). In particular, let us assume that $N^*_{F_1/X}$ is
generically spanned and its restriction to a generic line is
$\O(1)\oplus\O(1)$. A line $l\subset F_1$ moves out of $F_1$ if and only
if it is a jumping line of the bundle $N_{F_1/X}$, which means that the
restriction of the latter bundle to $l$ is either $\O\oplus\O(-2)$ or
$\O(1)\oplus\O(-3)$.

\proof
We have an exact sequence of vector bundles on a line $l\subset F_1$
$$0\raa N_{l/F_1}\iso\O(1)\raa N_{l/X}\raa (N_{F_1/X})_l\raa 0$$
and therefore $H^1(l,N_{l/X})\ne 0$ if and only if $l$ is a jumping of
$N_{F_1/X}$. If $l\in F_1^\vee\cap{\cal H}$ then the Hilbert scheme is
singular at $l$ and thus, by a general property of the Hilbert scheme, the
above cohomology does not vanish. Conversely, suppose that $l$ is a
jumping line of $N_{F_1/X}$. We may assume that the restriction of
$N_{F_1/X}$ to $l$ is $\O\oplus\O(-2)$ (because the jumping lines of the
other type are limits of these, if the conormal is spanned).
Let $\pi : \hat X\ra X$ be the blow-up of $X$ along $F_1$ with the
exceptional divisor $\hat F_1$.  Then over $l$ we have a unique curve $l_0$
such that $l_0.\hat F_1=0$. Therefore $K_{\hat X}.l_0=-1$ and computing
the deformation of $l_0$ on $\hat X$ we find out that it moves in
in a 2-dimensional family. So $l_0$ moves out of $\hat F_1$ and thus $l$
moves out of $F_1$.

\medskip

Let us note that the above argument describes the intersection $\hat E\cap
\hat F_1$ in the manifold $\hat X$, where $\hat E$ is the strict tranform
of the divisor $E$. Namely, the intersection is equal to the locus of all
sections $l_0$ which are described above. Knowing the list of possible
bundles (see (5.7.6)) we find out that
$\hat E_1:=\hat F_1\cap\hat E\subset \hat F_1$
is the unique irreducible divisor in $\hat F_1$ such that $\hat
E_1.l_0=-1$. On the other hand, since $l_0$ moves out of $\hat F_1$ and
becomes a fiber of contraction of $\hat E$ we seee that $\hat E.l_0=-1$ as
well. This implies however that the intersection of $\hat F_1$ and $\hat
E$ is transversal, or equivalently, $\O_{\hat F_1}(\hat E)=\O_{\hat
F_1}(\hat E_1)$.  Therefore we can count the multiplicity of $E$ along
$F_1$ using the jumping lines of $N_{F_1}$. The result is as it follows.

\proclaim Lemma (6.7).
$$mult_{F_1}E=c_2(N_{F_1/X})-1$$

\proof The left-hand-side of the above formula is equal to the
intersection number of the divisor $\hat E$ with a fiber of the blow-down
$\hat X\ra X$. The right-hand-side is equal to the degree of the curve of
jumping lines of $N_{F_1/X}$ in the dual plane $F_1^{\vee}$ which is the
same as the number of jumping lines passing through a generic point of
$F_1$ or, equivalently, the degree of the map $\hat E_1\ra F_1$. Since
$\O_{\hat F_1}(\hat E)=\O_{\hat F_1}(\hat E_1)$ then the above equality
follows.

\medskip

If $F$ is a quadric $\F_0$ or $\S_2$ then a similar argument can be applied
and the result is the following

\proclaim Lemma (6.8).
If the fiber $F$ is an irreducible quadric then $mult_F E=2$.

\proof
We use the notation from the proof of the previous lemma. If $F=F_1=\F_0$
then $\hat E_1$ consits of two components which are sections of
$\hat F_1\ra F$; consequently the argument can be applied to each of them.
If $F=F_1=\S_2$ then $\hat E_1$ is not Cartier but $2\hat E_1$ is Cartier
and we check that, as above, $2\hat E_1={\hat E}_{F_1}$.

\bigskip

Let $\hat{\cal H}$ be a component of the Hilbert scheme of
$\hat X$ which contains a (lift-up of) general line in $E$. Then
there is a map $\hat{\cal H}\ra {\cal H}$ which is identity outside of the
set parametrizing lines lying in $F$. Moreover, if $F$ is irreducible
then, by the proof of Lemma (6.6) the map $\hat{\cal H}\ra{\cal H}$ is
bijective over $F$ as well. If $F$ is irreducible then all fiber of
the contraction $\hat X\ra \hat Z$ are of dimension 1
(we know that all the fibers of the map $\hat
E\ra \hat S$ are 1-dimensional) and therefore the natural map
$\hat{\cal H}\ra\hat S$ is bijective.

\medskip \noindent
(6.9) Now let us discuss the case of $F=\P^2\cup\P^2$.
If $F$ is reducible then the Hilbert scheme parametrizing lines consists
of 3 irreducible components: ${\cal H}$ and $F_1^\vee$, $F_2^\vee$, the
latter two are just dual $\P^2$s. Because of Lemma (6.7),
after possible renumeration of them we
have $mult_{F_1}(E)=2$ and $mult_{F_2}(E)=3$.
According to the lemma (6.6)
the components meet along subsets parametrizing jumping lines. In
particular, the curve $F_1^\vee \cap {\cal H}$ is a reducible conic (two lines)
in $F_1^\vee$ and $F_2^\vee \cap {\cal H}$ is a rational cubic in $F_2^\vee$.
The components $F_1^\vee$ and $F_2^\vee$ meet at one point
which is the singularity of
each of the two preceeding curves. The point at which all the three
components meet is associated to the unique line $F_1\cap F_2$ with the
normal bundle $\O(1)\oplus\O(1)\oplus\O(-3)$.

The exceptional of divisor $\hat F$ consits of two components:
$\hat F_1=\P(\O\oplus\O(1))$ and $\hat F_2=\P(T\P^2(-1))$.
The map $\hat\f$ maps $F_1$ to $\hat G=\P^3$ contracting a section $\tilde F_1$
of $\hat F_1\ra F_1$ to a point $y_0\in \hat G$. The component $F_2$
is contracted by $\hat\f$ to a plane in $\hat G$ containing $y_0$.
The exceptional set of $\hat\f$ is equal to $\hat F_2\cup\hat E$.

Using the results on blow-ups (Lemma (2.2) and the subsequent discussion
which describes how to pass from $\P(N_{F_i/X})$ to $\P(N_{F_1\cup F_2/X})$
using vector bundle surgery)
and the above discussion on Hilbert
we can describe the intersection of $\hat E$ with $\hat F_1$ and $\hat F_2$.
The intersection $\hat F_1\cap\hat F_2$ is a ruled surface $\F_1$ (ruled over
$F_1\cap F_2$) with a $(-1)$ curve $l_0$. Then the locus of singular
points of $\hat X$ is a section $C_s$ of $\hat F_1\cap\hat F_2\ra F_1\cap F_2$
such that $C_s\cdot l_0=2$. In other words, if we contract $l_0$ in $\F_1$
then $C_s$ becomes a rational cubic in $\P^2$.
If we transform jumping lines of $N_{F_i/X}$ to $\hat F_i$ and use Lemma (2.2)
then we can describe the locus of lines from $\hat{\cal H}$. Namely,
in $\hat F_1$ these are the lines in $\tilde F_1$ which pass
through one of the points in $l_0\cap C_s$; in $\hat F_2$ these are
the fibers in $\hat\f^{-1}(\hat\f(C_s))$. Thus
$\hat F_1\cdot\hat E=2\tilde F_1$
and $\hat F_2\cdot \hat E=\hat\f^{-1}\f(C_s)$
(the multiplicity of the intersection can be verified with the result of
Lemma (6.7)).

From the above discussion it follows that the intersection $\hat S\cap\hat G$
is the rational cubic $\hat\f(C_s)\subset\hat\f(F_2)\subset\hat G$ which
is singular at $y_0$. Also we find out that the map $\hat{\cal H}\ra {\cal H}$
resolves the ``triple component'' singularity at the point which represents
the curve with normal $\O(1)^2\oplus\O(-3)$. That is, the lines in
$\hat{\cal H}$ over the point $z$ are parametrized by three rational
curves forming a simple chain. The map $\hat{\cal H}\ra \hat S$ contracts
the two exterior lines and loops the central one to the rational cubic
by glueing its meeting points. Therefore $\hat S$ is not normal at $y_0$.

\bigskip

\beginsection References.

\smalltype{

\item{[An]} Ando, T., On extremal rays of the higher
dimensional varieties, Invent. Math. {\bf 81} (1985), 347---357.

\item {[A-W]} Andreatta, M, Wi\'sniewski, J.A., A note on non vanishing
and applications, Duke Math. J., {\bf 72} (1993), 739-755.

\item{[A-B-W]} Andreatta,M., Ballico, E., Wi\'sniewski, J.A., Two theorems on
elementary contractions, Math. Ann. {\bf 297} (1993), 191-198.

\item {[B-P-V]} Barth,W., Peters, C. , Van de Ven, A., {\it Compact Complex
Surfaces},
Ergebnisse der Mathematik, {\bf 4} (1984), Springer-Verlag.

\item{[B-W]} Ballico, E., Wi\'sniewski,J.A., On B\v anic\v a sheaves and
Fano manifolds, to appear in Compositio Math..

\item{[Be]} Beltrametti, M., Contractions of non numerical effective extremal
rays in dimension 4, Proceedings of Conf. on Algebraic Geometry,
Berlin 1985, Teubner-Texte Math., {\bf 92} (1987), 24-37.

\item {[B-S]} Beltrametti, M., Sommese A.J.,
{\it The adjunction theory of complex projective varieties},
Expositions in Mathematics, Vol. 16 (1995), Walter De Gruyter

\item{[C-K-M]} Clemens,H., Koll\'ar,J., Mori,S., Higher dimensional
complex geometry, Asterisque {\bf 166} (1988).

\item{[CC]} CoCoA, Computation in Commutative Algebra, version 1.5,
by A. Giovini - G. Niesi.

\item{[Cu]} Cutkosky,S., Elementary Contractions of Gorenstein Threefolds,
Math. Ann., {\bf 280} (1988), 521-525.

\item{[Da]} Danilov, V. I., Birational Geometry of Toric 3-folds,
Math USSR Izvetia {\bf 21} n.2 (1983), 266-280.

\item{[Ei]} Ein, L., Varieties with small dual varieties, II,
Duke Math J., {\bf 52} (1985), 895-907.

\item{[E-W]} Esnault, H., Vieheweg, E., {\it Lectures on Vanishing Theorems},
DMV Seminar Band 20, Birkh\"auser Verlag (1992).

\item{[F-H-S]} Foster, O., Hirschowitz, A., Schneider, M., Type de scindage
g\'en\'eralis\'e pur les fibr\'es stables, Proc. of Vector Bundles and
Differential Equations, Nice (1979), Progress in Mathematics {\bf 7},
Birkh\"auser (1980), 65-82.

\item{[Fu]} Fujita, T., Classification theories of polarized varieties,
London Lect. Notes {\bf 115}, Cambridge Press 1990.

\item{[Fl]} Fulton,W., {\it Introduction to Toric Varieties}, Annals
of Math. Studies {\bf 131}, Princeton Univ. Press (1993).

\item{[EGA]} Grothendieck,A. - Dieudonne,J.,
{\it \'El\'ements de G\'eom\'etrie Alg\'ebrique,
vol. II}, Publ. Math. I.H.E.S., n.8 (1961).

\item{[Ha]} Hartshorne, R.,{\it Algebraic Geometry}, Springer-Verlag, 1977.

\item{[Io]} Ionescu, P., Generalized adjunction and applications, 
Math.~Proc.~Cambridge Math.~Soc.~{\bf 99} (1988), 457--472.

\item{[Kac]} Kachi, Y.,
Extremal contractions from 4-dimensional manifolds to 3-folds, preprint.

\item{[Ka1]} Kawamata, Y., Small contractions of four dimensional
algebraic manifolds, Math. Ann., {\bf 284} (1989), 595-600.

\item{[Ka2]} Kawamata, Y., General hyperplane sections of nonsingular
flops in dimension 3, preprint.

\item{[K-M-M]} Kawamata, Y., Matsuda, K., Matsuki, K.: Introduction to the
Minimal Model Program in {\it Algebraic Geometry, Sendai}, Adv. Studies
in Pure Math. {\bf 10}, Kinokuniya--North-Holland 1987, 283---360.

\item{[Ka-Mo]}  Katz, S., Morrison, D.R., Gorenstein threefold singularities
with small resolution via invariant theory of Weyl groups,
J. Alg. Geom. {\bf 1}, 449-530.

\item{[Ke]} Kempf, G.~et al, Toroidal embeddings I, 
Lecture Notes in Math.~{\bf 339} (1973), Springer Verlag.

\item{[Ko]} Koll\'ar, J., {\it Rational Curves on Algebraic Varieties},
Springer Verlag, Ergebnisse der Math. {\bf 32}, 1995.

\item{[Ko1, Ko2]} Koll\'ar, J., Higher direct images of dualizing sheaves I 
and II, Ann.~Math.~{\bf 123} (1986), 11--42 and {\bf 124} (1986), 171--202.

\item{[Ko-Mo]} Koll\'ar, J., Mori, S.,   Classification of three-dimensional
flips, Journal of the A.M.S., {\bf 5} (1992), 533-703.

\item{[Ma]} Maruyama,M., Moduli of stable sheaves I, J. Math. Kyoto Univ.,
{\bf 17}
(1977), 91-126.

\item{[Mo1]} Mori, S, Threefolds whose canonical bundles are not numerically
effective, Ann. Math., {\bf 116} (1982), 133-176.

\item{[Mo2]} Mori,S., Flip theorem and the existence of minimal models for
3-folds, Journal of the A. M. S., {\bf 1}, 1988, 117-224.

\item{[Od]} Oda, Convex bodies and algebraic geometry. {\it An introduction
to the theory of toric varieties}, Springer Verlag 1988.

\item{[O-S-S]} Okonek, C., Schneider, M., Spindler, H., Vector Bundles on
Complex Projective Spaces, Progress in Mathematics {\bf 3},
Birkh\"auser 1980.

\item{[Re]} Reid,M., Minimal Models of Canonical 3-folds, Proc. Algebraic
Varieties and Analytic Varieties-Tokyo 1981, Adv. Studies
in Pure Math. {\bf 1}, Kinokuniya--North-Holland 1982, 131-180.

\item{[S-W1]} Szurek, M, Wi\'sniewski,J.A., Fano Bundles of rank 2 on surfaces,
Compositio Mathematica, {\bf 76} (1990), 295-305.

\item{[S-W2]} Szurek, M, Wi\'sniewski,J.A., On Fano manifolds, which are
{\bf P}$^k$-
bundles over {\bf P}$^2$, Nagoya Math. J., {\bf 120} (1990), 89-101.

\item{[Wi]} Wi\'sniewski,J.A., On contractions of extremal rays of Fano
manifolds,
Jounal f\"ur die reine und angew. Mathematik, {\bf 417} (1991), 141-157.

\item{[Y-Z]} Ye, Y. and Zhang, Q., On ample vector bundles whose adjunction
bundles are not numerically effective, Duke Math. J. {\bf 60} (1990), 671---687.

}

\end